\documentclass[11pt]{article}
\usepackage{graphicx,amsmath,fullpage,latexsym,amssymb,tikz,pgfplots}
\usepackage{psfrag}
\usetikzlibrary{shapes.multipart}

\newcommand{\reals}{{\mbox{\bf R}}}
\newcommand{\BEAS}{\begin{eqnarray*}}
\newcommand{\EEAS}{\end{eqnarray*}}
\newcommand{\BEA}{\begin{eqnarray}}
\newcommand{\EEA}{\end{eqnarray}}
\newcommand{\BEQ}{\begin{equation}}
\newcommand{\EEQ}{\end{equation}}
\newcommand{\BIT}{\begin{itemize}}
\newcommand{\EIT}{\end{itemize}}

\newcommand{\ie}{{\it i.e.}}

\newcommand{\intr}{\mathop{\bf int}}
\newcommand{\Range}{\mbox{\textrm{range}}}
\newcommand{\diag}{\mathop{\bf diag}}
\newcommand{\Tr}{\mathop{\bf tr}}
\newcommand{\argmin}{\mathop{\rm argmin}}

\newcommand{\SV}[2]{\symm^{#1}_{#2}} 
\newcommand{\SVp}[2]{\symm^{#1}_{{#2},+}}
 
\newcommand{\SVc}[2]{\symm^{#1}_{{#2},\mathrm{c}}}

\newcommand{\symm}{{\mbox{\bf S}}}

\newcommand{\svec}{\mathrm{\bf vec}}

\newcommand{\prox}{\mathrm{prox}}

\newcommand{\prnt}{\mathop{\mathrm{pa}}}
\newcommand{\ch}{\mathop{\mathrm{ch}}}

\title{Decomposition in conic optimization with partially separable 
structure}
\author{
Yifan Sun\thanks{
    Electrical Engineering Department, University of California, 
    Los Angeles. Email: {\tt ysun01@ucla.edu}, 
   {\tt lieven.vandenberghe@ucla.edu}.
   Research supported by NSF Grants DMS-1115963 and
   ECCS-1128817.}
\and 
Martin S. Andersen\thanks{
    Department of Applied Mathematics and Computer Science,
    Technical University of Denmark, 
    Email: {\tt mskan@dtu.dk}.} 
  \and 
   Lieven Vandenberghe\footnotemark[1]}
\date{}

\begin{document}
\maketitle
\begin{abstract}
Decomposition techniques for linear programming are difficult to
extend to conic optimization problems with general non-polyhedral convex 
cones because the conic inequalities introduce an additional nonlinear 
coupling between the variables.
However in many applications the convex cones have a 
partially separable structure that allows them to be characterized
in terms of simpler lower-dimensional cones.
The most important example is sparse semidefinite programming with a
chordal sparsity pattern.  Here partial separability derives from 
the clique decomposition theorems that characterize positive 
semidefinite and positive-semidefinite-completable matrices 
with chordal sparsity patterns.
The paper describes a decomposition method that exploits
partial separability in conic linear optimization.
The method is based on Spingarn's method for equality constrained
convex optimization, combined with a fast interior-point method for 
evaluating proximal operators.
\end{abstract}

\section{Introduction} \label{s-intro}
We consider conic linear optimization problems (conic LPs)
\BEQ \label{e-conic-intro}
 \begin{array}{ll}
 \mbox{minimize} & c^T x \\
 \mbox{subject to} & Ax = b \\ & x\in\mathcal C
 \end{array}
\EEQ
in which the cone $\mathcal C$ is defined in terms of lower-dimensional 
convex cones $\mathcal C_k$ as
\BEQ \label{e-C-intro}
 \mathcal C = \{ x\in \reals^n \mid x_{\gamma_k} \in \mathcal C_k,
 \; k=1,\ldots, l\}.
\EEQ
The sets $\gamma_k$ are ordered subsets of $\{1,2,\ldots,n\}$ 
and $x_{\gamma_k}$ denotes the subvector of $x$ with entries indexed
by $\gamma_k$.
We refer to the structure in the cone $\mathcal C$ as \emph{partial 
separability}. 
The purpose of the paper is to describe a decomposition method 
that exploits partially separable structure.

In (standard) linear optimization, with $\mathcal C =\reals^n_+$,
the cone is separable, \ie, a product of one-dimensional cones, 
and the coupling of variables and constraints 
is entirely specified by the sparsity pattern of $A$.  
The term \emph{decomposition} in linear optimization usually refers
to techniques for exploiting \emph{angular} or \emph{dual-angular}
structure in the coefficient matrix $A$, \ie, a sparsity pattern
that is almost block-diagonal, except for a small number of dense rows or 
columns \cite{Las:02b,BeT:97a}.  The goal of a decomposition 
algorithm is to solve the problem iteratively, by solving a sequence
of separable problems, obtained by removing the complicating
variables or constraints.  
The decoupled  subproblems can be solved 
in parallel or sequentially (for example, to reduce memory usage).
Moreover, if the iterative coordinating process is simple enough to be 
decentralized, the decomposition method can be used as a distributed
algorithm.
By extension, decomposition methods can also be applied to more general
sparsity patterns for which removal of complicating variables and 
constraints makes the problem substantially easier to solve 
(even if it does not decompose into independent subproblems).

When the cone $\mathcal C$ in~(\ref{e-conic-intro}) is not 
separable or block-separable (a product of lower-dimensional cones), 
the formulation of decomposition algorithms is 
more complicated because the inequalities  introduce an 
additional coupling between the variables.   
However if the cone is partially separable, as defined 
in~(\ref{e-C-intro}), and the
overlap between the index sets $\gamma_k$ is small, one can still
formulate efficient decomposition algorithms.  The purpose
of this paper is to discuss such a decomposition method.
The method is based on Spingarn's method of partial inverses for convex 
optimization
problems with equality constraints
\cite{Spi:83,Spi:85}, combined with an interior-point 
method applied to sparse conic subproblems.  
The details of the method are described in 
sections~\ref{s-par-sep} and~\ref{s-conic}.

An important example of a partially separable cone are
the positive-semidefinite-completable sparse matrices with  
a chordal sparsity pattern.  Matrices in this cone are characterized by
the property that all their principal dense submatrices are positive
semidefinite \cite[theorem 7]{GJSW:84}.
This fundamental result has been applied in previous methods for 
sparse semidefinite optimization.
It is the basis of the  conversion methods 
used to reformulate sparse semidefinite programs (SDPs) in equivalent 
forms that are easier to handle  by interior-point 
algorithms  \cite{KKMY:11,FKMN:00} or more suitable for distributed
algorithms via the Alternating Direction Method of Multipliers (ADMM)
\cite{DZG:12}.
Partial separability also underlies the saddle-point 
mirror-prox algorithm for `well-structured' sparse SDPs proposed by Lu,
Nemirovski and Monteiro \cite{ZNM:07}.
We discuss the sparse semidefinite optimization application
of the decomposition method in detail in 
sections~\ref{s-sdp}--\ref{s-experiments}.

\paragraph{Notation} \label{s-notation}
If $\alpha$ is a subset of $\{1,2,\ldots,n\}$, then $E_\alpha$ 
will denote the $|\alpha| \times n$-matrix with entries
\[
 (E_\alpha)_{ij} = \left\{
 \begin{array}{ll} 1 & \alpha(i) = j \\
                   0 & \mbox{otherwise}. 
 \end{array}\right.
\]
Here $\alpha(i)$ is the $i$th element of $\alpha$, sorted 
using the natural ordering.
If not explicitly stated the column dimension $n$ of $E_\alpha$ 
will be clear from the context.
The result of multiplying an $n$-vector $x$ with 
$E_\alpha$ is the subvector of $x$ of length $|\alpha|$ with elements 
$(x_\alpha)_k = x_{\alpha(k)}$.  
The adjoint operation $x = E_\alpha^T y$ maps an $|\alpha|$-vector $y$
to an $n$-vector $x$ by copying the entries of $y$ to the positions
indicated by $\alpha$, \ie, by setting $x_{\alpha(k)} = y_k$ and 
$x_i = 0$ for $i\not\in \alpha$.
Therefore $E_\alpha E_\alpha^T$ is an identity matrix of order
$|\alpha|$ and $E_\alpha^T E_\alpha$ is a diagonal $0$-$1$ matrix of 
order $n$, with $i$th diagonal entry equal to one if and only if 
$i\in\alpha$.
The matrix $P_\alpha = E_{\alpha}^T E_{\alpha}$ represents projection 
in $\reals^n$ on the sparse $n$-vectors with support $\alpha$.

Similar notation will be used for principal submatrices in a symmetric 
matrix.
If $X\in\symm^p$ (the symmetric matrices of order $p$)
and $\alpha$ is a subset of $\{1,\ldots,p\}$, then 
\[
 \mathcal E_{\alpha}(X) = X_{\alpha\alpha} = E_\alpha X E_\alpha^T
\in \symm^{|\alpha|}.
\]
This is the submatrix of order $|\alpha|$ with 
$i,j$ entry $(X_{\alpha\alpha})_{ij} = X_{\alpha(i)\alpha(j)}$.
The adjoint operation $\mathcal E_\alpha^*$ copies a matrix 
$Y\in\symm^{|\alpha|}$ to an otherwise zero symmetric $p\times p$-matrix:
\[
\mathcal E_{\alpha}^*(Y) = E_\alpha^T Y E_\alpha.
\]
The projection of a matrix $X\in\symm^p$ on the matrices that are zero 
outside of a diagonal $\alpha\times \alpha$ block is denoted
\[
\mathcal P_{\alpha}(X) = P_\alpha X P_\alpha
= E_{\alpha}^TE_\alpha X E_\alpha^T E_\alpha.
\]

\section{Partially separable cones}  
\label{s-par-sep}

\subsection{Partially separable functions}
A function $f:\reals^n \rightarrow \reals$ is \emph{partially separable}
if it can be expressed as
\[
 f(x) = \sum_{k=1}^l f_k(A_kx),
\]
where each $A_k$ has a nontrivial nullspace, 
\ie, a rank substantially less than $n$.
This concept was introduced by Griewank and Toint
in the context of quasi-Newton algorithms 
\cite{GrT:82,GrT:84}\cite[section 7.4]{NoW:06}.  
Here we consider the simplest and most common example of partial
separability and assume that $A_k = E_{\gamma_k}$ for some
index set $\gamma_k\subset \{1,2,\ldots,n\}$.
This means that $f$ can be written as a sum of functions that 
depend only on subsets of the components of $x$:
\BEQ \label{e-part-sep}
 f(x) = \sum_{k=1}^l f_k(E_{\gamma_k}x)
 = \sum_{k=1}^l f_k(x_{\gamma_k}).
\EEQ
Partial separability generalizes \emph{separability} 
($l=n$, $\gamma_k = \{k\}$) and \emph{block-separability} 
(the sets $\gamma_k$ form a partition of $\{1,2,\ldots,n\}$).

\subsection{Partially separable cones} \label{s-parsep-cones}
We call a cone $\mathcal C \subset \reals^n$ 
\emph{partially separable} if it can be expressed as
\BEQ \label{e-C-primal}
 \mathcal C = \{x \mid E_{\gamma_k}x \in \mathcal C_k, \; k=1,\ldots,l\}
\EEQ
where $\mathcal C_k$ is a convex cone in $\reals^{|\gamma_k|}$. 
The terminology is motivated by the fact the indicator function 
$\delta_\mathcal C$ of $\mathcal C$ is a partially separable function:
\[
\delta_\mathcal C(x) = 
\sum_{k=1}^p \delta_{\mathcal C_k} (E_{\gamma_k}x)
\]
where $\delta_{\mathcal C_k}$ is the indicator function of $\mathcal C_k$.
The following assumptions will be made.
\BIT
\item The index sets $\gamma_k$ are \emph{distinct} and \emph{maximal},
 \ie, $\gamma_i \not\subseteq \gamma_j$ for $i\neq j$, and their
 union is equal to $\{1,2,\ldots,n\}$. 
 
\item The convex cones $\mathcal C_k$ are proper, \ie, closed, pointed,
 with nonempty interior. This implies that their dual cones 
 \[
 \mathcal C_k^* = \{v \in\reals^{|\gamma_k|} \mid
   u^T v \geq 0 \; \forall u\in\mathcal C_k\}
 \]
are proper convex cones and that $\mathcal C_k = \mathcal C_k^{**}$
\cite[page 53]{BoV:04}.

\item There exists a point $\bar x$ with $E_{\gamma_k}\bar x \in \intr
\mathcal C_k$ for $k=1,\ldots,l$.
\EIT
These assumptions imply that $\mathcal C$ is itself a proper cone.
Indeed, $\mathcal C$ is clearly convex, with nonempty interior
(the point $\bar x$ is in the interior).
It is closed because it can be expressed as an intersection of closed 
halfspaces:
\[
 \mathcal C = \{x \in\reals^n \mid
  v_k^T E_{\gamma_k} x \geq 0 \; \; \forall v_k \in \mathcal C_k^*,
 \; k=1,\ldots, l\}.
\]
Finally, $\mathcal C$ is pointed because $x \in \mathcal C$, 
$-x\in\mathcal C$
implies $E_{\gamma_k} x \in\mathcal C_k$ and $-E_{\gamma_k} x\in 
\mathcal C_k$ for all $k$.  Since the cones $\mathcal C_k$ are
pointed, this means $E_{\gamma_k}  x= 0$ for $k=1,\ldots,l$.
Since the index sets $\gamma_k$ cover $\{1,2,\ldots,n\}$,
this implies $x=0$.

It follows that the dual cone $\mathcal C^*$ is proper.
It can be verified that 
\BEQ \label{e-C-dual}
{\mathcal C}^* = 
\{\sum_{k=1}^l E_{\gamma_k}^T \tilde s_k \mid 
\tilde s_k \in \mathcal C_k^*, \; k=1,\ldots,l \}.
\EEQ

\paragraph{Example}
Take $n=6$ and
\BEQ \label{e-example-sets}
 \gamma_1 = \{1,2,6\}, \qquad
 \gamma_2 = \{2,5,6\}, \qquad
 \gamma_3 = \{3,5\}, \qquad
 \gamma_4 = \{4,6\}.
\EEQ
Let
$\mathcal C_1 \subset \reals^3$, $\mathcal C_2 \subset \reals^3$,
$\mathcal C_3 \subset \reals^2$, $\mathcal C_4 \subset \reals^2$ 
be proper convex cones. 
A vector $x\in\reals^6$ is in the cone $\mathcal C$ defined 
in~(\ref{e-C-primal}) if
\[
 (x_1,x_2,x_6) \in \mathcal C_1, \qquad
 (x_2,x_5,x_6) \in \mathcal C_2, \qquad
 (x_3,x_5) \in \mathcal C_3, \qquad
 (x_4,x_6) \in \mathcal C_4.
\]
A vector $s \in\reals^6$  is in the dual cone $\mathcal C^*$ if
\[
 s = \left[\begin{array}{c}
 \tilde s_{11} \\ \tilde s_{12} \\ 0 \\ 0 \\ 0 \\ \tilde s_{13} 
 \end{array}\right] + 
 \left[\begin{array}{c}
 0 \\ \tilde s_{21} \\ 0 \\ 0 \\ \tilde s_{22} \\ \tilde s_{23} 
 \end{array}\right] + 
 \left[\begin{array}{c}
 0 \\ 0 \\ \tilde s_{31} \\ 0 \\ \tilde s_{32} \\ 0
 \end{array}\right] + 
 \left[\begin{array}{c}
 0 \\ 0 \\ 0 \\ \tilde s_{41} \\ 0 \\ \tilde s_{42}
 \end{array}\right]
\]
for some
\[
 \tilde s_1 = 
 (\tilde s_{11}, \tilde s_{12}, \tilde s_{13}) \in \mathcal C_1^*, \quad 
 \tilde s_2 = 
 (\tilde s_{21}, \tilde s_{22}, \tilde s_{23}) \in \mathcal C_2^*, \quad 
 \tilde s_3 = 
 (\tilde s_{31}, \tilde s_{32}) \in \mathcal C_3^*, \quad 
 \tilde s_4 = (\tilde s_{41}, \tilde s_{42}) \in \mathcal C_4^*. 
\]

\subsection{Sparsity and intersection graph} \label{e-int-graph}
Two different undirected graphs can be associated with the partially 
separable structure defined by the index sets 
$\gamma_1$, \ldots, $\gamma_l$.
These graphs will be referred to as the \emph{sparsity graph}
and the \emph{intersection graph}.

\paragraph{Sparsity graph} 
The sparsity graph $\mathcal G$ 
has $n$ vertices, representing the $n$ variables.
There is an edge between two distinct vertices $i$ and $j$  
if $i$, $j \in \gamma_k$ for some $k$.
We call this the sparsity graph because it represents the sparsity pattern
of a matrix
\BEQ \label{e-sp-H}
 H = \sum_{k=1}^l E_{\gamma_k}^T H_k E_{\gamma_k}
\EEQ
where the matrices $H_k$ are dense symmetric matrices.
For example, if the component functions $f_k$ in~(\ref{e-part-sep})
are twice differentiable with dense Hessians, then the Hessian of $f$ 
has this structure.  The entries 
$(i,j) \not\in \cup_{k=1,\ldots,l} \, (\gamma_k \times \gamma_k)$ 
are the positions of the zeros in the sparsity pattern of $H$.

Each index set $\gamma_k$ thus defines a complete subgraph of 
the sparsity graph $\mathcal G$. 
Since the index sets are maximal (by assumption), these complete 
subgraphs are the cliques in $\mathcal G$. 

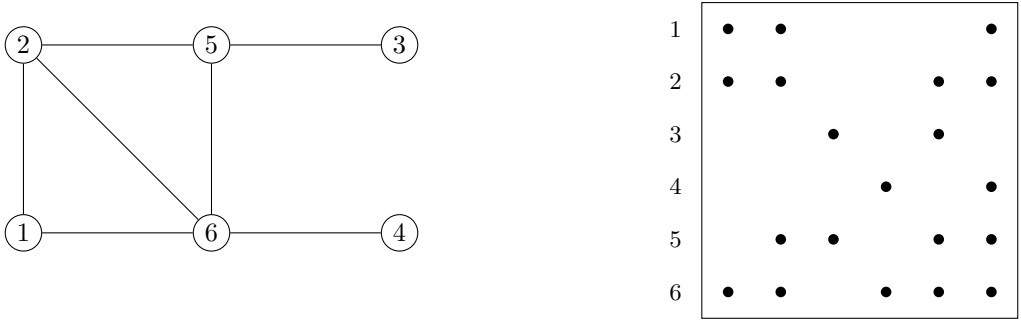
\begin{figure}
\hspace*{\fill}
\begin{minipage}{.4\linewidth}
\centering
\begin{tikzpicture}[scale=5]
   \tikzset{VertexStyle/.style = {shape = circle, draw, 
       minimum size = 12pt, inner sep = 2, fill=none}}

   \foreach \l/\k/\d/\x/\y in {
       1/1/1/0.0/0.0,
       2/2/2/0.0/0.5,
       3/3/3/1.0/0.5,
       4/4/2/1.0/0.0,
       5/5/3/0.5/0.5,
       6/6/2/0.5/0.0}
       \node[VertexStyle, font=\small](\k) at (\x,\y){\l};

   \foreach \i/\j in {1/2, 1/6, 2/5, 2/6, 3/5, 4/6, 5/6} \draw (\i)--(\j);
\end{tikzpicture}
\hspace*{\fill}
\end{minipage}
\hspace*{\fill}
\begin{minipage}{.4\linewidth}
\hspace*{\fill}
\begin{tikzpicture}[scale=0.7]
   \tikzset{VertexStyle/.style = {shape = circle, minimum size = 4pt, 
       inner sep = 0 pt, fill=black}}
   \tikzset{VertexStyleR/.style = {shape = circle, minimum size = 4pt, 
       inner sep = 0 pt, fill=red}}

   \foreach \i/\j in {
       1/0, 5/0, 
       4/1, 5/1, 
       4/2, 
       5/3, 
       5/4 
      }
      {
       \node[VertexStyle] at (\j,-\i){};
       \node[VertexStyle] at (\i,-\j){};
      }
   \foreach \k/\i in { 1/0, 2/1, 3/2, 4/3, 5/4, 6/5}
      {
      \node[VertexStyle] at (\i,-\i){};
      \node[font=\footnotesize] at (-1,-\i){\k};
      \node[font=\footnotesize] at (\i,1){\k};
      }

   \draw (-0.5,0.5) rectangle (5.5,-5.5);
\end{tikzpicture}
\hspace*{\fill}
\end{minipage}
\hspace*{\fill}

\caption{Sparsity graph and sparsity pattern for an example with 
$n=6$ and four index sets $\gamma_1 = \{1,2,6\}$, 
$\gamma_2 = \{2,5,6\}$, 
$\gamma_3 = \{3,5\}$, 
$\gamma_4 = \{4,6\}$.} \label{f-sparsity}
\end{figure}

Figure~\ref{f-sparsity} shows the sparsity graph and sparsity pattern 
for the index sets~(\ref{e-example-sets}).

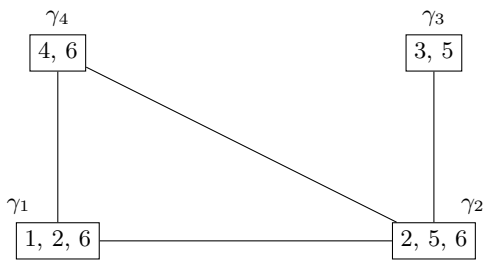
\begin{figure}
\hspace*{\fill}
\begin{minipage}{.5\linewidth}
\hspace{\fill}
\begin{tikzpicture}[scale=5, clique/.style = {rectangle, draw}],
 \footnotesize
 \node[clique][label = above left: $\gamma_1$](1) at (0,0){1, 2, 6};
 \node[clique][label = above right: $\gamma_2$](2) at (1,0){2, 5, 6};
 \node[clique][label = above: $\gamma_3$](3) at (1,0.5){3, 5};
 \node[clique][label = above: $\gamma_4$](4) at (0,0.5){4, 6};
 \draw (1) -- (2);
 \draw (1) -- (4);
 \draw (2) -- (3);
 \draw (2) -- (4);
\end{tikzpicture}
\hspace*{\fill}
\end{minipage}
\hspace*{\fill}

\caption{Intersection graph for the same example as in 
Figure~\ref{f-sparsity}.} \label{f-ex-intersection}
\end{figure}

\paragraph{Intersection graph}
The intersection graph has the index sets $\gamma_k$ as its vertices
and an edge between distinct vertices $i$ and $j$ if the sets 
$\gamma_i$ and $\gamma_j$ intersect.  
We place a weight $|\gamma_i \cap \gamma_j|$ on edge $\{i,j\}$.
The intersection graph is therefore identical to the clique graph 
of the sparsity graph~$\mathcal G$.
(The clique graph of an undirected graph has the cliques of the
graph as its vertices and undirected edges between cliques that intersect,
with edge weights equal to the sizes of the intersection.)
An example is shown in Figure~\ref{f-ex-intersection}.

\subsection{Chordal structure} \label{s-chordal}
An undirected graph is \emph{chordal} if for every cycle of length 
greater than three there is a chord 
(an edge connecting non-adjacent vertices in the cycle).
If the sparsity graph representing a partially separable
structure is chordal (as will be the case in the application to
semidefinite optimization discussed in the second half of the paper),
several additional useful properties hold.
In this section we summarize the most important of these properties. 
For more background and proofs we refer the reader to the 
survey paper \cite{BlP:93}.

\paragraph{Running intersection property}
A spanning tree of the intersection graph
(or, more accurately, a spanning forest,
since we do not assume the intersection graph is connected)
has the \emph{running intersection property} if
\[
 \gamma_i \cap \gamma_j \subseteq \gamma_k 
\]
whenever vertex $\gamma_k$ is on the path between vertices 
$\gamma_i$ and $\gamma_j$ in the tree.  
A fundamental theorem states that a spanning tree with the running 
intersection property exists if and only if the corresponding sparsity 
graph is chordal \cite{BlP:93}.  

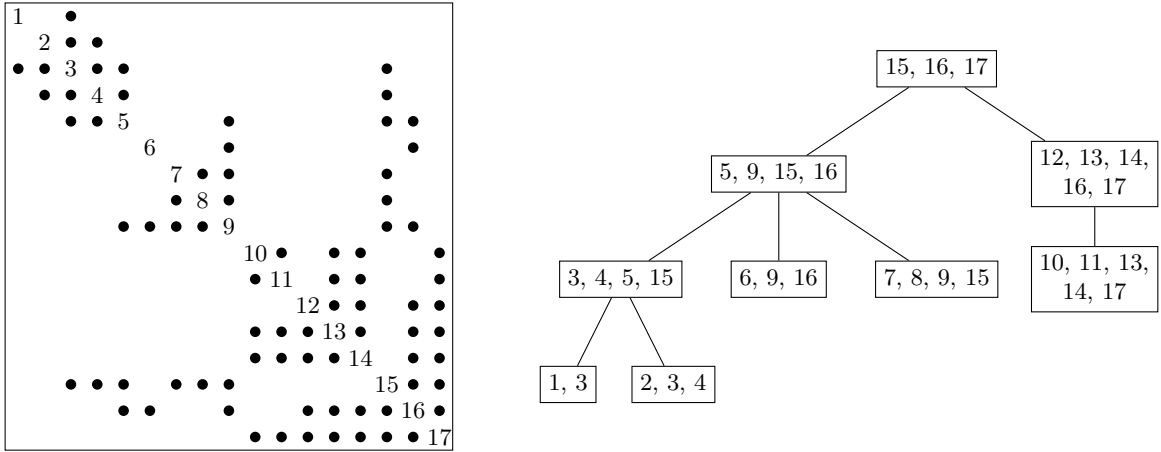
\begin{figure}
\hspace*{\fill}
\begin{minipage}{.40\linewidth}
\hspace*{\fill}
\begin{tikzpicture}[scale=0.35]
   \tikzset{VertexStyleA/.style = {shape = circle, minimum size = 4pt, 
       inner sep = 0pt, fill=black}}

   \foreach \i/\j in {
       2/0,
       2/1, 3/1,
       3/2, 4/2, 14/2,
       4/3, 14/3,
       8/4, 14/4, 15/4, 
       8/5, 15/5,
       7/6, 8/6, 14/6,  
       8/7, 14/7, 
       14/8, 15/8,
       10/9, 12/9, 13/9, 16/9,
       12/10, 13/10, 16/10,
       12/11, 13/11, 15/11, 16/11,
       13/12, 15/12, 16/12,
       15/13, 16/13,
       15/14, 16/14,
       16/15 
      }
      {
      \node[VertexStyleA] at (\j,-\i){};
      \node[VertexStyleA] at (\i, -\j){};
      }

   \foreach \k/\i in {
       4/3, 8/7, 9/8, 11/10, 13/12, 14/13, 17/16} 
      \node[font=\footnotesize] at (\i,-\i) {$\k$};

   \foreach \k/\i in {
       1/0, 2/1, 3/2, 5/4, 6/5, 7/6, 10/9, 12/11, 15/14, 16/15} 
      \node[font=\footnotesize] at (\i,-\i) {$\k$};

   \draw (-0.5,0.5) rectangle (16.5,-16.5);
\end{tikzpicture}
\hspace*{\fill}
\end{minipage}
\hspace*{\fill}
\begin{minipage}{.55\linewidth}
\hspace*{\fill}
\begin{tikzpicture}[scale=1.4, clique/.style = {rectangle, draw}]
 \footnotesize
 \node[clique](1) at (1,-1){ 1, 3 };
 \node[clique](2) at (2,-1){ 2, 3, 4 };
 \node[clique](3) at (1.5,0){ 3, 4, 5, 15};
 \node[clique](6) at (3,0){6, 9, 16};
 \node[clique](7) at (4.5,0){ 7, 8, 9, 15};
 \node[clique](5) at (3,1){ 5, 9, 15, 16};
 \node[clique](10) at (6,0){ \parbox{4.5em}{\centering 10, 11, 13, 14, 17}};
 \node[clique](12) at (6,1){ \parbox{4.5em}{\centering 12, 13, 
      14, 16, 17}};
 \node[clique](15) at (4.5,2){15, 16, 17}; 
 \draw (1) -- (3);
 \draw (2) -- (3);
 \draw (3) -- (5);
 \draw (6) -- (5);
 \draw (7) -- (5);
 \draw (5) -- (15);
 \draw (10) -- (12);
 \draw (12) -- (15);
\end{tikzpicture}
\hspace*{\fill}
\end{minipage}
\hspace*{\fill}

\caption{Spanning tree in an intersection graph for nine index sets
$\gamma_k$ and  $n = 17$ variables (right),
and the sparsity pattern for the corresponding sparsity graph (left).
The sparsity graph is chordal.  It can be verified that the spanning 
tree on the right-hand side has the running intersection property.}
\label{f-rip}
\end{figure}

The right-hand figure in Figure~\ref{f-rip} shows a spanning tree of 
the intersection graph of $l=9$ index sets $\gamma_k$ with $n=17$ 
variables.  On the left-hand side we represent the corresponding sparsity 
graph as a sparse matrix pattern 
(a dot in positions $i,j$ and $j,i$ indicates an edge $\{i,j\}$).
It can be verified that the 
tree satisfies the running intersection property.  

It is easy to see that a spanning tree with the running intersection
property is a maximum weight spanning tree (if the weight of 
edge $\{i,j\}$ is defined as the size of the intersection 
$\gamma_i \cap \gamma_j$.)
To show this, assume the spanning tree has the running intersection
property but is not a maximum weight spanning tree.  
Then there exists an edge $\{\gamma_i, \gamma_j\}$ in the intersection 
graph that is not an edge of the tree and that can be substituted for an 
edge on the path between $\{\gamma_i,\gamma_j\}$ in the tree, say edge 
$\{\gamma_s, \gamma_t\}$, to obtain 
a spanning tree with larger weight.
This means that the edge $\{\gamma_i,\gamma_j\}$ is heavier than the edge 
$\{\gamma_s,\gamma_t\}$, 
\ie, $|\gamma_s \cap \gamma_t| < |\gamma_i \cap \gamma_j|$.
However this contradicts the running intersection property,
which states that $\gamma_i \cap \gamma_j \subseteq \gamma_s$
and $\gamma_i \cap \gamma_j \subseteq \gamma_t$.

It is therefore not surprising that chordality of a sparsity 
graph can be tested efficiently using modifications of algorithms for
constructing maximum-weight spanning trees in graphs.
An example is the \emph{maximum-cardinality search} algorithm
\cite{TaY:84,BlP:93}.

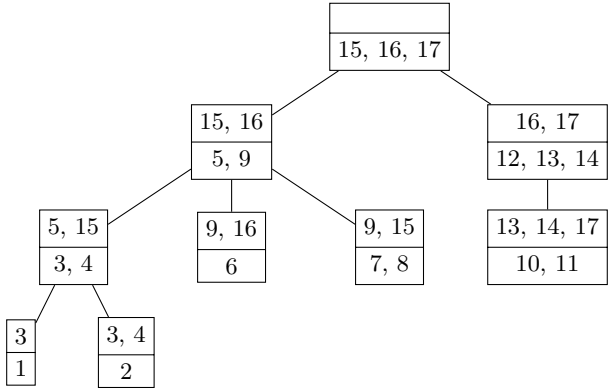
\begin{figure}
\hspace*{\fill}
\begin{minipage}{.65\linewidth}
\hspace*{\fill}
\begin{tikzpicture}[scale=1.4, 
clique/.style = {rectangle split, rectangle split parts = 2, draw}]
 \footnotesize
 \node[clique](1) at (1,-1){ 3 \nodepart{second} 1 };
 \node[clique](2) at (2,-1){ 3, 4 \nodepart{second} 2 };
 \node[clique](3) at (1.5,0){ 5, 15 \nodepart{second} 3, 4};
 \node[clique](6) at (3,0){9, 16 \nodepart{second} 6 };
 \node[clique](7) at (4.5,0){ 9, 15 \nodepart{second} 7, 8};
 \node[clique](5) at (3,1){ 15, 16 \nodepart{second} 5, 9};
 \node[clique](10) at (6,0){ 13, 14, 17 \nodepart{second} 10, 11};
 \node[clique](12) at (6,1){ 16, 17 \nodepart{second} 12, 13, 14   };
 \node[clique](15) at (4.5,2){ \nodepart{second} 15, 16, 17 }; 
 \draw (1) -- (3);
 \draw (2) -- (3);
 \draw (3) -- (5);
 \draw (6) -- (5);
 \draw (7) -- (5);
 \draw (5) -- (15);
 \draw (10) -- (12);
 \draw (12) -- (15);
\end{tikzpicture}
\hspace*{\fill}
\end{minipage}
\hspace*{\fill}

\caption{Each vertex $\gamma_k$ in the intersection tree of 
Figure~\ref{f-rip} is split in two sets $\alpha_k$ and 
$\gamma_k\setminus\alpha_k$ with $\alpha_k$ the intersection of 
$\gamma_k$ and its parent.
The indices listed in the top row of each vertex form $\alpha_k$.  
The indices in the bottom row form $\gamma_k \setminus \alpha_k$.}
\label{f-rip2}
\end{figure}

\paragraph{Properties} Suppose a spanning tree with the running 
intersection property exists.  
We partition each index set $\gamma_k$ in two 
sets $\alpha_k$ and $\gamma_k\setminus\alpha_k$, defined as follows.
If $\gamma_k$ is the root of the tree, then 
$\alpha_k = \emptyset$.
For the other vertices,
\[
 \alpha_k  = \gamma_k \cap \gamma_{\prnt(\gamma_k)}
\]
where $\prnt(\gamma_k)$ is the parent of $\gamma_k$ in the tree.
Note that $\alpha_k$ is a strict subset of $\gamma_k$ because 
$\alpha_k = \gamma_k$  would imply 
$\gamma_k = \alpha_k \subseteq \prnt{(\gamma_k)}$,
contrary to our assumption that $\gamma_k \not\subseteq \gamma_j$
for $j\neq k$.
The definition of the sets $\alpha_k$ is illustrated in 
Figure~\ref{f-rip2}.

The running intersection property has the following implications
\cite{LPP:89,PoS:90}.

\BIT
\item Every index $i \in\{1,2,\ldots,n\}$ belongs to at least one 
set $\gamma_k\setminus\alpha_k$.

This is easily seen by contradiction.
By assumption, each index $i$ belongs to at least one set $\gamma_k$.   
Suppose $i \in \alpha_j$ whenever $i\in \gamma_j$.
By definition of $\alpha_j$, this implies that 
$i\in \gamma_{\prnt(\gamma_j)}$ whenever $i\in\gamma_j$.
Therefore $i\in \gamma_r$ where $r$ is the root of the tree.
Since $\alpha_r = \emptyset$, this contradicts the assumption that
$i$ does not belong to any set $\gamma_k\setminus\alpha_k$. 

\item If an element $i\in \gamma_k\setminus\alpha_k$ is contained in 
$\gamma_j$, $j\neq k$,
then $\gamma_j$ is a descendant of $\gamma_k$.

Suppose $\gamma_j$ is not a descendant of $\gamma_k$.
Then the path connecting $\gamma_j$ and $\gamma_k$ in the spanning tree
includes the parent of $\gamma_k$ and by the running intersection
property, $i$ is an element of the parent of $\gamma_k$.
This implies $i\in\alpha_k$.

\item Every index $i\in\{1,2,\ldots,n\}$ 
belongs to at most one set $\gamma_k\setminus\alpha_k$.

This follows directly from the previous property: 
$i\in\gamma_j\setminus\alpha_j$ and $i\in\gamma_k\setminus\alpha_k$ 
implies that the vertex $\gamma_j$ is a descendant of vertex $\gamma_k$ 
in the tree and vice-versa, so $j=k$.
\EIT
Combining the first and third properties, we conclude that the sets 
$\gamma_k\setminus\alpha_k$, 
$k=1,\ldots,l$, form a partition of $\{1,2,\ldots,n\}$.
This is illustrated in Figure~\ref{f-rip2}:
the indices in the bottom rows of the vertices in the tree 
are the sets $\gamma_k\setminus\alpha_k$ 
and form a partition of $\{1,2,\ldots,17\}$.

\section{Conic optimization with partially separable cones} 
\label{s-conic}
We now consider a pair of conic linear optimization problems
\BEQ \label{e-conic}
 \begin{array}[t]{ll}
 \mbox{minimize} & c^T x \\
 \mbox{subject to} & Ax = b \\ & x\in\mathcal C
 \end{array} \qquad \qquad
 \begin{array}[t]{ll}
 \mbox{maximize} & b^Ty\\
 \mbox{subject to} & A^T y + s = c \\ & s\in\mathcal C^*
 \end{array} 
\EEQ
with respect to a partially separable cone~(\ref{e-C-primal})
and its dual~(\ref{e-C-dual}).
The variables are $x,s \in\reals^n$, $y\in\reals^m$.
In addition to the assumptions listed in Section~\ref{s-parsep-cones}
we assume that the sparsity graph associated with the index sets 
$\gamma_k$ is chordal and that a maximum weight spanning tree (or forest)
$T$ in the intersection graph is given.
We refer to $T$ as the \emph{intersection tree}
and use the notation $\prnt(\gamma_k)$ and $\ch(\gamma_k)$ for 
the parent and the children of vertex $\gamma_k$ in $T$.  

\subsection{Reformulation}
The decomposition algorithm developed in the following 
sections is based on a straightforward reformulation of the
conic LPs~(\ref{e-conic}).  
The primal and dual cones can be expressed  as 
\[
\mathcal C = \{x \mid Ex \in \tilde {\mathcal C}\}, \qquad
\mathcal C^* = \{E^T\tilde s \mid \tilde s \in 
\tilde {\mathcal C}^*\},
\]
where 
$\tilde {\mathcal C} = \mathcal C_1 \times \cdots \times \mathcal C_l$
and $\tilde {\mathcal C}^* = \mathcal C_1^* \times \cdots \times 
\mathcal C_l^*$, and $E$ is the $\tilde n \times n$ matrix
\[
 E = \left[\begin{array}{llll}
 E_{\gamma_1}^T & E_{\gamma_2}^T & \cdots & E_{\gamma_l}^T
 \end{array}\right]^T
\]
with $\tilde n = \sum_k |\gamma_k|$.
Define $\mathcal V = \Range(E)$. 
A change of variables $\tilde x = Ex$, $s= E^T \tilde s$
allows us to write the problems~(\ref{e-conic}) equivalently as
\BEQ \label{e-converted-V}
\begin{array}[t]{ll}
\mbox{minimize} & \tilde c^T \tilde x \\
\mbox{subject to} & \tilde A \tilde x  = b \\
 & \tilde x \in \mathcal V \\
 & \tilde x \in \tilde {\mathcal C} 
\end{array} \qquad\qquad
\begin{array}[t]{ll}
\mbox{maximize} & b^Ty \\
\mbox{subject to} & \tilde A^Ty + v + \tilde s = \tilde c \\
 & v  \in \mathcal V^\perp \\
 & \tilde s \in \tilde {\mathcal C}^*
\end{array} 
\EEQ
with variables 
$\tilde x = (\tilde x_1, \ldots, \tilde x_l) \in \reals^{\tilde n}$,
$y\in\reals^m$, 
$\tilde s = (\tilde s_1, \ldots, \tilde s_l) \in\reals^{\tilde n}$,
provided $\tilde A$ and $\tilde c$ are chosen to satisfy
\[
\tilde A E = \sum_{k=1}^l \tilde A_k E_{\gamma_k} = A, \qquad
E^T \tilde c = \sum_{k=1}^l E_{\gamma_k}^T \tilde c_k = c.
\]
Here $\tilde A_k$ and $\tilde c_k$ are blocks of size $|\gamma_k|$
in the partitioned matrix and vector
\BEQ \label{e-AJ}
\tilde A = \left[\begin{array}{cccc}
 \tilde A_1 & \tilde A_2 & \cdots & \tilde A_l
\end{array}\right], \qquad 
\tilde c^T = \left[\begin{array}{cccc}
 \tilde c_1^T & \tilde c_2^T & \cdots & \tilde c_l^T
\end{array}\right].
\EEQ
It is straightforward to find $\tilde A$ and $\tilde c$
that satisfy these conditions.
For example, one can take $\tilde A = AJ$, $\tilde c = J^T c$ with
$J$ equal to 
\BEQ \label{e-J}
 J = \left[\begin{array}{cccc}
  P_{\gamma_1\setminus\alpha_1} E_{\gamma_1}^T & 
  P_{\gamma_2\setminus\alpha_2} E_{\gamma_2}^T & 
  \cdots &
  P_{\gamma_l\setminus\alpha_l} E_{\gamma_l}^T 
 \end{array}\right]
\EEQ
or any other left-inverse of $E$.  
However we will see later that other choices of $\tilde A$ may offer
advantages.

\begin{figure}
\begin{minipage}{.4\linewidth}
\centering
\begin{tikzpicture}[scale=2.2,
    clique/.style = {rectangle split, rectangle split parts=2, draw} ]
 \small
 \node[clique][label = left :$\tilde x_j$](1) at (0,0){
     $\alpha_j$ \nodepart{second} $\gamma_j \setminus \alpha_j$ };
 \node[clique][label = left: $\tilde x_k$](2) at (1,1.0){
     $\alpha_k$ \nodepart{second} $\gamma_k \setminus \alpha_k$ };
 \node[clique][label = left: $\tilde x_i$](3) at (2,2){
     $\alpha_i$ \nodepart{second} $\gamma_i \setminus \alpha_i$ };
 \node(4) at (-0.4,-0.4){ };
 \node(5) at (0.4,-0.4){ };
 \node(6) at (0,-0.4){ };
 \node(7) at (1.4,0.6){ };
 \node(8) at (1,0.6){ };
 \node(9) at (2.4,1.6){ };
 \node(10) at (2,1.6){ };

 \node at (0.5,.5)[label = left: 
   $E_{\alpha_j} (E_{\gamma_j}^T \tilde x_j -
    E_{\gamma_k}^T \tilde x_k) {=} 0$]{};
 \node at (1.5,1.5)[label = left: 
   $E_{\alpha_k} (E_{\gamma_k}^T \tilde x_k -
    E_{\gamma_i}^T \tilde x_i) {=} 0$]{};
 \node at (-.2,-0.4)[label = center: $\cdots$]{};
 \node at (.2,-0.4)[label = center: $\cdots$]{};
 \node at (1.2,0.6)[label = center: $\cdots$]{};
 \node at (2.2,1.6)[label = center: $\cdots$]{};

 \draw (1) -- (2);
 \draw (1) -- (4);
 \draw (1) -- (5);
 \draw (1) -- (6);
 \draw (2) -- (3);
 \draw (2) -- (7);
 \draw (2) -- (8);
 \draw (3) -- (9);
 \draw (3) -- (10);
\end{tikzpicture}
\end{minipage}
\hspace*{\fill}
\begin{minipage}{.4\linewidth}
\centering
\begin{tikzpicture}[scale=2.2,
    clique/.style = {rectangle split, rectangle split parts=2, draw} ]
 \node[clique][label = left :$\tilde s_j$](1) at (0,0){
     $\alpha_j$ \nodepart{second} $\gamma_j \setminus \alpha_j$ };
 \node[clique][label = left: $\tilde s_k$](2) at (1,1.0){
     $\alpha_k$ \nodepart{second} $\gamma_k \setminus \alpha_k$ };
 \node[clique][label = left: $\tilde s_i$](3) at (2,2){
     $\alpha_i$ \nodepart{second} $\gamma_i \setminus \alpha_i$ };
 \node(4) at (-0.4,-0.4){ };
 \node(5) at (0.4,-0.4){ };
 \node(6) at (0,-0.4){ };
 \node(7) at (1.4,0.6){ };
 \node(8) at (1,0.6){ };
 \node(9) at (2.4,1.6){ };
 \node(10) at (2,1.6){ };

 \node at (0.5,.5)[label = left: $u_j$]{};
 \node at (1.5,1.5)[label = left: $u_k$]{};
 \node at (-.2,-0.4)[label = center: $\cdots$]{};
 \node at (.2,-0.4)[label = center: $\cdots$]{};
 \node at (1.2,0.6)[label = center: $\cdots$]{};
 \node at (2.2,1.6)[label = center: $\cdots$]{};

 \draw (1) -- (2);
 \draw (1) -- (4);
 \draw (1) -- (5);
 \draw (1) -- (6);
 \draw (2) -- (3);
 \draw (2) -- (7);
 \draw (2) -- (8);
 \draw (3) -- (9);
 \draw (3) -- (10);
\end{tikzpicture}
\end{minipage}
\hspace*{\fill}

\caption{\emph{The subspaces $\mathcal V$ and $\mathcal V^\perp$.}
The figures show three vertices of the intersection tree.
The left-hand figure illustrates $\mathcal V$.
We associate the subvector $\tilde x_k$ of 
$\tilde x = (\tilde x_1, \ldots, \tilde x_l)$ 
with vertex $\gamma_k$ in the tree and associate a consistency constraint
$E_{\alpha_j} (E_{\gamma_j}^T \tilde x_j - E_{\gamma_k}^T \tilde x_k) 
= 0$ with the edge between vertex~$\gamma_j$ and its parent~$\gamma_k$.
Then $(\tilde x_1, \ldots, \tilde x_l)$ is in $\mathcal V$
if and only if the consistency constraints are satisfied.
The right-hand figure illustrates $\mathcal V^\perp$.
Here we associate the subvector $\tilde s_k$ of 
$\tilde s = (\tilde s_1, \ldots, \tilde s_l)$ with 
vertex $\gamma_k$ in the tree 
and a vector $u_j \in \reals^{|\alpha_j|}$ with the edge
between vertex $\gamma_j$ and its parent $\gamma_k$. 
Then $\tilde s = (\tilde s_1, \ldots, \tilde s_l)$ is in 
$\mathcal V^\perp$ if and only there exist values of  $u_k$ 
such that $\tilde s_k = E_{\gamma_k} (E_{\alpha_k}^T u_k -
 \sum_{\gamma_j\in\ch(\gamma_k)} E_{\alpha_j}^T u_j )$.
\label{f-primal-conversion}
\label{f-dual-conversion}}
\end{figure}
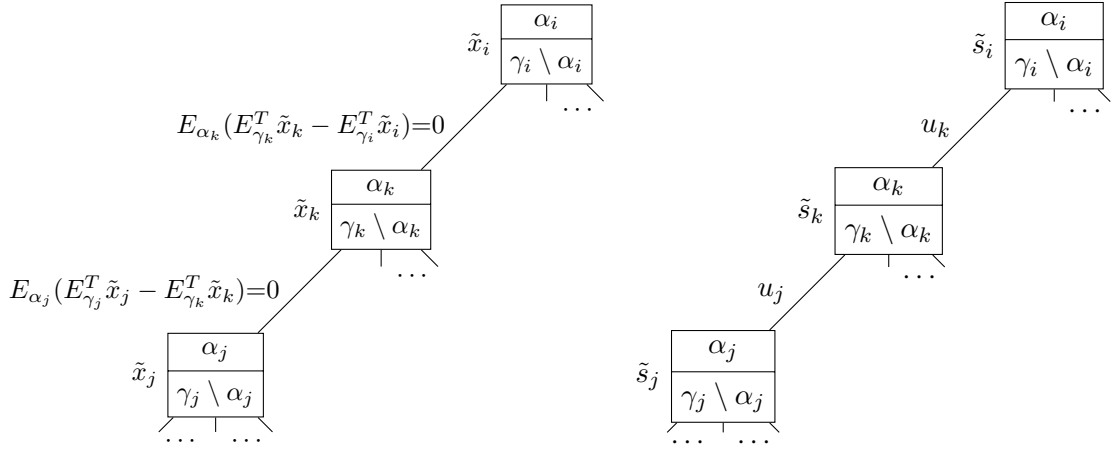

\paragraph{Consistency constraint}
The running intersection property of the intersection tree $T$ 
can be used to derive a simple representation of
the subspaces $\mathcal V$ and $\mathcal V^\perp$.
We first note that a vector 
$\tilde x = (\tilde x_1, \ldots, \tilde x_l)$ is in $\mathcal V$ 
if and only if
\BEQ \label{e-primal-cone-overlap}
E_{\alpha_j} (E_{\gamma_j}^T \tilde x_j - E_{\gamma_k}^T \tilde x_k) = 0, 
\quad k=1,\ldots,l, \quad \gamma_j\in\ch(\gamma_k).
\EEQ
This can be seen as follows.
Since $\alpha_j = \gamma_j \cap \gamma_{\prnt(\gamma_j)}$ by definition,
the equalities~(\ref{e-primal-cone-overlap}) mean that
\BEQ \label{e-primal-cone-overlap-2}
  E_{\gamma_j \cap \gamma_k} (E_{\gamma_j}^T \tilde x_j 
  -E_{\gamma_k}^T\tilde x_k) = 0 
\EEQ
for all $\gamma_k$ and all $\gamma_j \in \ch(\gamma_k)$.  
This is sufficient to guarantee that~(\ref{e-primal-cone-overlap-2})
holds for \emph{all} $j$ and $k$, because the running intersection 
property guarantees that if $\gamma_j$ and $\gamma_k$ intersect then 
their intersection is included in every index set on  the path between 
$\gamma_j$ and $\gamma_k$ in the tree.
The equations~(\ref{e-primal-cone-overlap}) therefore hold if and
only if there is an $x$ such that 
$\tilde x_k = E_{\gamma_k} x$ for $k=1,\ldots,l$, \ie,
$x\in\mathcal V$. 
We will refer to the constraint $\tilde x \in\mathcal V$ 
as the \emph{consistency constraint} in~(\ref{e-converted-V}).
It is needed to ensure that the variables $\tilde x_k$  
can be interpreted as  copies $\tilde x_k = E_{\gamma_k} x$ of
overlapping subvectors of some $x\in\reals^n$.
This is illustrated graphically in Figure~\ref{f-primal-conversion}.

Likewise, a vector $\tilde s = (\tilde s_1,\ldots, \tilde s_l)$
is in $\mathcal V^\perp$ if and only if
there exist $u_j\in\reals^{|\alpha_j|}$, $j=1,\ldots,l$, such that
\BEQ \label{e-Vperp}
 \tilde s_k = E_{\gamma_k} (E_{\alpha_k}^T u_k -
 \sum\limits_{\gamma_j\in\ch(\gamma_k)} E_{\alpha_j}^T u_j),
 \quad k=1,\ldots,l.
\EEQ
This is illustrated in the left-hand part of 
Figure~\ref{f-dual-conversion}.  

The equations~(\ref{e-primal-cone-overlap}) 
and~(\ref{e-Vperp}) can be written succinctly as
\[
   B \tilde x= 0, \qquad \tilde s = B^T u,
\]
where $u = (u_1,\ldots, u_l) \in\reals^{|\alpha_1|} \times \cdots
\times \reals^{|\alpha_l|}$ and the matrix $B$ is constructed as follows.
Define an $l\times l$ matrix $S$ with elements
\[
 S_{ij} = \left\{ \begin{array}{ll}
  1 & i = j \\
 -1 & \gamma_i = \prnt(\gamma_j) \\
  0 & \mbox{otherwise.}
 \end{array}\right.
\]
This is the transpose of the node-arc incidence matrix of 
the spanning tree $T$, if we direct the arcs from children to parents. 
Define
\[
 B = \left[\begin{array}{ccc}
 E_{\alpha_1} & \cdots & 0 \\
 \vdots & \ddots & \vdots \\
 0 & \cdots & E_{\alpha_l} \end{array}\right]
 (S^T \otimes I_n) 
\left[\begin{array}{ccc}
 E_{\gamma_1}^T & \cdots & 0 \\
 \vdots & \ddots & \vdots \\
 0 & \cdots & E_{\gamma_l}^T \end{array}\right] 
\]
where $S\otimes I_n$ is the Kronecker product.
The matrix $B$ is an $l\times l$ block matrix with
diagonal blocks $E_{\alpha_k} E_{\gamma_k}^T$, $k=1,\ldots, l$.
Block row $j$ of $B$ has a nonzero block 
$-E_{\alpha_j}E_{\gamma_k}^T$ in block column $k$, where 
$k=\prnt(\gamma_j)$.  The rest of the matrix $B$ is zero.
If the vertices of the intersection tree are numbered in a topological 
ordering (\ie, each vertex receives a lower number than its parent),
then the matrix $S$ is lower triangular.
Note that the sparsity pattern of $B^TB$ can be embedded
in the chordal sparsity pattern of the sparsity graph.
(If we ignore the sparsity within the blocks 
$E_{\alpha_j} E_{\gamma_k}^T$ and treat these blocks as dense, 
then the sparsity pattern of $B^TB$ is exactly the sparsity pattern
of the sparsity graph.)
The matrix $BB^T$, on the other hand, is not necessarily sparse.

With this notation the reformulated primal and dual 
problems~(\ref{e-converted-V}) are
\BEQ \label{e-converted-matrix} 
\begin{array}[t]{ll}
 \mbox{minimize}   & \tilde c^T \tilde x \\
 \mbox{subject to} & \tilde A \tilde x = b \\
                   & B \tilde x = 0 \\
                   & \tilde x \in \tilde {\mathcal C},
\end{array} \qquad\qquad
\begin{array}[t]{ll}
 \mbox{maximize} & b^T y \\
 \mbox{subject to} & \tilde A^T y + B^T u + \tilde s = \tilde c \\
                   & \tilde s \in \tilde{\mathcal  C}^*.
\end{array}
\EEQ

\subsection{Correlative sparsity} \label{s-corr-sparse}
The reformulated problems generalize the clique-tree conversion 
methods proposed for semidefinite programming 
in \cite{KKMY:11,FKMN:00}.
These conversion methods were proposed 
with the purpose of reformulating large, sparse SDPs in
an equivalent form that is easier to solve by interior-point methods.
In this section we discuss the benefits of the reformulation 
in the context of general conic optimization problems with partially 
separable cones.  
The application to semidefinite programming is discussed
in the next section.

The reformulated problems (\ref{e-converted-matrix})
are of particular interest if the sparsity of the matrix 
$\tilde A$ implies that a matrix of the form 
\BEQ \label{e-corr-sparse}
 \tilde A G \tilde A^T = \sum_{k=1}^l \tilde A_k G_k \tilde A_k^T,
\EEQ
where $G$ is block-diagonal, with arbitrary dense diagonal 
blocks $G_k$, is sparse.  
We call the sparsity pattern of $\tilde AG\tilde A^T$ the
\emph{correlative sparsity} pattern of the reformulated problem, 
after Kobayashi \emph{et al.}\ \cite{KKK:08}.
The correlative sparsity pattern can be determined as follows:
the $i,j$ entry of $\tilde AG\tilde A^T$ is zero if there are no
block columns $\tilde A_k$ in which the $i$th and $j$th row are both
nonzero.  
The correlative sparsity pattern clearly depends on the choice of 
$\tilde A$ as illustrated by the following example.

Consider a small conic LP with $m=4$, $n=6$, index sets 
$\gamma_k$ given in~(\ref{e-example-sets}), and a constraint matrix $A$
with zeros in the following positions:
\[
 A = \left[\begin{array}{cccccc}
  A_{11} & A_{12} & 0      & 0      & 0      & A_{16} \\
  0      & A_{22} & 0      & 0      & A_{25} & A_{26} \\
  0      & 0      & 0      & A_{34} & 0      & A_{36} \\
  0      & 0      & A_{43} & 0      & A_{45} & 0
 \end{array}\right].
\]
In other words, equality $i$ in $Ax=b$ 
involves only variables $x_k \in \gamma_i$.
\begin{figure}
\hspace*{\fill}
\begin{minipage}{.65\linewidth}
\hspace*{\fill}
\begin{tikzpicture}[scale=1.4, 
clique/.style = {rectangle split, rectangle split parts = 2, draw}]
 \footnotesize
 \node[clique][label = left: $\gamma_1$](1) at 
     (1,0){ 2, 6 \nodepart{second} 1 };
 \node[clique][label = left: $\gamma_3$](2) at 
     (3,0){ 6 \nodepart{second} 4 };
 \node[clique][label = left: $\gamma_2$](3) at 
     (2,1){ 5 \nodepart{second} 2, 6};
 \node[clique][label = left: $\gamma_4$](4) at 
     (2,2){  \nodepart{second} 3, 5 };
 \draw (1) -- (3);
 \draw (2) -- (3);
 \draw (3) -- (4);
\end{tikzpicture}
\hspace*{\fill}
\end{minipage}
\hspace*{\fill}
\caption{Spanning tree in the intersection graph of
Figure~\ref{f-ex-intersection}.}  \label{f-ex-tree}
\end{figure}
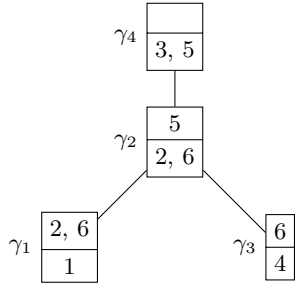
The primal reformulated problem has a variable
$\tilde x = (\tilde x_1, \tilde x_2, \tilde x_3, \tilde x_4)
\in \reals^3 \times \reals^3 \times \reals^2 \times \reals^2$.
If we use the intersection tree shown in Figure~\ref{f-ex-tree},
The consistency constraints are  $B\tilde x= 0$ with
\[
 B = \left[\begin{array}{ccc|ccc|cc|cc}
  0 & -1 & 0 & 1 & 0 & 0 & 0 & 0 & 0 & 0 \\ 
  0 &  0 & -1 & 0 & 0 & 1 & 0 & 0 & 0 & 0 \\ 
  0 &  0 &  0 & 0 & 0 & 1 & 0 & -1 & 0 & 0 \\ 
  0 &  0 &  0 & 0 & -1 & 0 & 0 & 0 & 0 & 1 
\end{array}\right].
\]
If we define $\tilde A$ and $\tilde c$ via~(\ref{e-AJ}) 
and~(\ref{e-J}), we obtain
\BEAS
\tilde A & = & \left[\begin{array}{ccc|ccc|cc|cc}
 A_{11} & 0 & 0 & A_{12} & 0 & A_{16} & 0      & 0 & 0      & 0 \\
 0      & 0 & 0 & A_{22} & 0 & A_{26} & 0      & 0 & 0      & A_{25} \\
 0      & 0 & 0 & 0      & 0 & A_{36} & A_{34} & 0 & 0      & 0 \\
 0      & 0 & 0 & 0      & 0 & 0      & 0      & 0 & A_{43} & A_{45} 
 \end{array}\right] \\
\tilde c & = & \left[\begin{array}{ccc|ccc|cc|cc}
 c_1 & 0 & 0 & c_2 & 0 & c_6 & c_4 & 0 & c_3 & c_5 \end{array}\right]^T.
\EEAS
With this choice the $4\times 4$ matrix~(\ref{e-corr-sparse}) is dense,
except for a zero in positions $(4,1)$, $(4,3)$, $(1,4)$, $(3,4)$.
On the other hand, 
if we choose 
\[
\tilde A = \left[\begin{array}{ccc|ccc|cc|cc}
 A_{11} & A_{12} & A_{16} & 0 & 0 & 0 & 0 & 0 & 0 & 0 \\
 0 & 0 & 0 & A_{22} & A_{25} & A_{26} & 0 & 0 & 0 & 0 \\
 0 & 0 & 0 & 0 & 0 & 0 & A_{34} & A_{36} & 0 & 0 \\
 0 & 0 & 0 & 0 & 0 & 0 & 0 & 0 & A_{43} & A_{45} 
 \end{array}\right],
\]
then the correlative sparsity pattern is diagonal.

\subsection{Interior-point methods} \label{s-decomp-ipm}
In this section we first compare the cost of 
interior-point methods applied to the reformulated and the original 
problems, for problems with correlative sparsity.

An interior-point method applied to the reformulated primal
and dual problems~(\ref{e-converted-matrix})
requires at each iteration the solution of a linear equation
(often called the 
\emph{Karush-Kuhn-Tucker (KKT)} equation)
\BEQ \label{e-kkt}
\left[\begin{array}{ccc}
H & \tilde A^T & B^T  \\ \tilde A & 0 & 0  \\ B &  0 & 0 
\end{array}\right] 
\left[\begin{array}{c} 
\Delta \tilde x \\ \Delta y \\ \Delta u \end{array}\right]
= \left[\begin{array}{c} r_{\tilde x} \\ r_y \\ r_u \end{array}\right]
\EEQ
where $H = \diag(H_1, \ldots, H_l)$ is a positive definite 
block-diagonal scaling matrix that depends on the algorithm used,
the cones $\mathcal C_k$, and the current primal and dual iterates 
in the algorithm.
Here we will assume that the blocks of $H$ are defined as
\[
H_k = \nabla^2\phi_k(w_k) 
\]
where $\phi_k$ is a logarithmic barrier function for $\mathcal C_k$
and $w_k$ is some point in $\intr\mathcal C_k$.
This assumption is sufficiently general to cover 
path-following methods based on primal scaling, 
dual scaling, and the Nesterov-Todd primal-dual scaling.
In most implementations, the KKT equation is solved by 
eliminating
$\Delta \tilde x$ and solving  a smaller system
\BEQ \label{e-converted-schur}
\left[\begin{array}{cc}
 \tilde A H^{-1} \tilde A^T  & \tilde AH^{-1} B^T \\
 B H^{-1}\tilde A^T  & B H^{-1} B^T
\end{array}\right] \left[\begin{array}{c} \Delta y \\ \Delta u 
 \end{array}\right] = \left[\begin{array}{c}
 \tilde A H^{-1}r_{\tilde x} - r_y \\ 
 B H^{-1}r_{\tilde x} - r_u \end{array}\right].
\EEQ
The coefficient matrix in~(\ref{e-converted-schur})
is called the \emph{Schur complement matrix}.
Note that the 1,1 block has the form~(\ref{e-corr-sparse}), 
so its sparsity pattern is the correlative sparsity pattern.  
The 2,2 block $BH^{-1}B^T$ on the other hand may be quite dense.  
The sparsity pattern of the coefficient matrix of~(\ref{e-converted-schur})
must be compared with the Schur complement matrix in an interior-point 
method applied to the original conic LPs~(\ref{e-conic}).  This matrix  
has the same sparsity pattern as the system obtained by
eliminating $\Delta u$ in~(\ref{e-converted-schur}), \ie,
\BEQ \label{e-schur-unconverted}
 \tilde A (H^{-1} - H^{-1}B^T(B H^{-1}B^T)^{-1} BH^{-1})\tilde A^T.
\EEQ
The matrix~(\ref{e-schur-unconverted})
is often dense (due to the $BH^{-1}B^T$ term), even 
for problems with correlative sparsity.

An interior-point method for the reformulated problem can exploit
correlative sparsity by solving~(\ref{e-converted-schur}) using
a sparse Cholesky factorization method.  
If $\tilde AH^{-1}\tilde A^T$ is nonsingular,
one can also explicitly eliminate $\Delta y$  and reduce it to a 
dense linear equation in $\Delta u$ with coefficient matrix
\[
B(H^{-1}- H^{-1}\tilde A^T 
(\tilde A H^{-1}\tilde A^T)^{-1} \tilde A H^{-1})B^T.
\]
To form this matrix one can take advantage of correlative sparsity.
(This is the approach taken in \cite{KKK:08}.)
Whichever method is used for solving~(\ref{e-converted-schur}),
the advantage of the enhanced sparsity resulting from the  
sparse 1,1 block $\tilde AH^{-1}\tilde A^T$ must be weighed against 
the increased size of the reformulated problem.
This is especially important for semidefinite programming, where 
the extra variables $\Delta u$ are vectorized matrices,
so the difference in size of the two Schur complement systems is very
substantial.

\subsection{Spingarn's method} \label{s-spingarn}
Motivated by the high cost of solving the KKT 
equations~(\ref{e-converted-schur}) of the converted problem
we now examine the alternative  of using a first-order splitting method
to exploit correlative sparsity.
The converted primal problem~(\ref{e-converted-V}) 
can be written as 
\BEQ \label{e-dr}
 \begin{array}{ll}
 \mbox{minimize} & f(\tilde x) \\
 \mbox{subject to} & \tilde x \in \mathcal V 
 \end{array}
\EEQ
where 
the cost function $f$ is defined as
\BEQ \label{e-f-def}
f(\tilde x) = \tilde c^T \tilde x + \delta(\tilde A\tilde x- b)
+ \delta_{\tilde{\mathcal C}}(\tilde x),
\EEQ
with $\delta$ and $\delta_{\tilde{\mathcal C}}$ the 
indicator functions for $\{0\}$ and $\tilde{\mathcal C}$, respectively.
Spingarn's \emph{method of partial inverses} \cite{Spi:83,Spi:85}
is a decomposition method that exploits separable structure in equality
constrained convex problems of the form~(\ref{e-dr}).
The method is known to be equivalent to the Douglas-Rachford 
splitting method  \cite{LiM:79,EcB:92} applied to 
\[
 \begin{array}{ll}
 \mbox{minimize} & f(\tilde x)  + \delta_{\mathcal V}(\tilde x).
 \end{array}
\]
Starting at some $z^{(0)}$, the following three steps are repeated:
\BEAS
  \tilde x^{(k)} & = & \prox_{f/\sigma} (z^{(k-1)}) \\
  w^{(k)} & = & P_{\mathcal V} (2\tilde x^{(k)} - z^{(k-1)}) \\
  z^{(k)} & = & z^{(k-1)} + \rho_k (w^{(k)} - \tilde x^{(k)}).
\EEAS
Here $P_\mathcal V$ denotes Euclidean projection on $\mathcal V$
and $\prox_{f/\sigma}$ is the \emph{proximal operator} of $f$,
defined as 
\[
\prox_{f/\sigma}(z) = 
\argmin_{\tilde x} \left(f(\tilde x) + 
 \frac{\sigma}{2} \|\tilde x-z\|_2^2\right).
\]
It can be shown that $\prox_{f/\sigma}(z)$
exists and is unique for all $z$ \cite{Mor:65,BaC:11}.  
The value $\tilde x = \prox_{f/\sigma}(z)$ 
of the prox-operator of the function~(\ref{e-f-def}) 
is the primal optimal solution in the pair of
conic quadratic optimization problems
\BEQ \label{e-proxop}
\begin{array}[t]{ll}
\mbox{minimize}
&  \displaystyle
c^T \tilde x  + \frac{\sigma}{2} \|\tilde x - z \|_2^2 \\
\mbox{subject to}
& \tilde A\tilde x =b \\
& \tilde x \in \tilde{\mathcal C}
\end{array}
\qquad\qquad
\begin{array}[t]{ll}
\mbox{maximize} & \displaystyle b^T y - 
\frac{1}{2\sigma} \|c - \tilde A^T y - \sigma z - \tilde s\|_2^2\\
\mbox{subject to} & \tilde s \in \tilde{\mathcal C}^*
\end{array}
\EEQ
with primal variables $\tilde x$ and dual variables $y$, $\tilde s$.
Equivalently, $\tilde x$, $y$, $\tilde s$ satisfy the 
optimality conditions
\BEQ \label{e-proxop-oc}
 \tilde A\tilde x = b, \qquad
 \tilde A^Ty + \tilde s + \sigma (z-\tilde x)  = c, \qquad
 \tilde x \in \mathcal C, \qquad
 \tilde s \in {\tilde{\mathcal C}}^*, \qquad
 \tilde x^T \tilde s = 0.
\EEQ
In the following discussion we assume that the prox-operator 
of $f$ is computed exactly, \ie, we do not explore the possibility of 
speeding up the algorithm by using inexact prox-evaluations.
This is justified if an interior-point method 
is used for solving~(\ref{e-proxop}), since interior-point methods
achieve a high accuracy and offer only a modest gain in efficiency 
if inaccurate solutions are acceptable.

The algorithm depends on two algorithm parameters: 
a positive constant $\sigma$ 
(we will refer to $1/\sigma$ as the \emph{steplength})
and a \emph{relaxation parameter} $\rho_k$, which can change at
each iteration but must remain in an interval
$(\rho_{\mathrm{min}}, \rho_{\mathrm{max}})$ with
$0 < \rho_{\mathrm{min}} < \rho_{\mathrm{max}} < 2$.
More details on the Douglas-Rachford method and its applications 
can be found in \cite{Eck:94,CoP:07,BaC:11,PaB:12}.

The complexity of the first two steps (the evaluation of the 
prox-operator and the projection on $\mathcal V$) will be discussed
later, after we make some general comments about the interpretation of
the method and stopping criteria.

\paragraph{Interpretation as fixed-point iteration}
The three steps in the Spingarn iteration can be combined 
into a single update
\BEQ
 z^{(k)} = z^{(k-1)} - \rho_k G(z^{(k-1)}) \label{e-dr-step-z}
\EEQ
with the operator $G$ defined as
\[
G(z) = \prox_{f/\sigma}(z) - P_\mathcal V(2\prox_{f/\sigma}(z) - z).
\]
For $\rho_k=1$ this is a fixed-point iteration $z^{(k)} = 
z^{(k-1)} - G(z^{(k-1)})$ for solving $G(z) = 0$;
for other values of $\rho_k$ it is a fixed-point
iteration with relaxation (underrelaxation for $\rho_k < 1$,
overrelaxation for $\rho_k > 1$).

Zeros of $G$ are related to the solutions of~(\ref{e-dr}) as
follows. If $z$ is a zero of $G$ then 
$x = \prox_{f/\sigma}(z)$ and $v = \sigma(z-x)$
satisfy the optimality conditions for~(\ref{e-dr}), which are
\BEQ \label{e-dr-oc}
 \tilde x \in\mathcal V, \qquad
 v \in\mathcal V^\perp, \qquad
 v\in \partial f(\tilde x).
\EEQ
Conversely, if $x$, $v$ satisfy these  optimality conditions,
then $z = x + (1/\sigma)v$ is a zero of $G$.
To see this, first assume $G(z) = 0$ and 
define $x = \prox_{f/\sigma}(z)$, $v=\sigma(z-x)$.  
By definition of the prox-operator, $v \in \partial f(x)$.
Moreover, $G(z) = 0$ gives 
$x = P_\mathcal V(x) + (1/\sigma) P_\mathcal V(v)$. 
Therefore $x\in\mathcal V$ and $v \in \mathcal V^\perp$.
Conversely suppose $x$, $v$ satisfy these  optimality conditions.
Define $z = x + (1/\sigma)v$.  Then it can be verified that
$x = \prox_{f/\sigma}(z)$ and
$G(z) = x - P_\mathcal V(x - (1/\sigma)v) = 0$.

\paragraph{Primal and dual residuals}
From step~1 in the algorithm and the definition of the proximal operator 
we see that the vector $v^{(k)} = \sigma( z^{(k-1)} - \tilde x^{(k)})$
satisfies $v^{(k)} \in \partial f(\tilde x^{(k)})$.
If we define 
\[
r_\mathrm{p}^{(k)} =  P_\mathcal V(\tilde x^{(k)}) - \tilde x^{(k)},
\qquad
r_\mathrm{d}^{(k)} =  -P_\mathcal V(v^{(k)})
\]
then
\BEQ \label{e-pd-oc}
 \tilde x^{(k)} + r_{\mathrm p}^{(k)} \in \mathcal V, \qquad 
 v^{(k)} + r_{\mathrm d}^{(k)} \in \mathcal V^\perp, \qquad
 v^{(k)} \in \partial f(\tilde x^{(k)}).
\EEQ
The vectors $r_{\mathrm p}^{(k)}$ and
$r_{\mathrm d}^{(k)}$ can be interpreted as primal and dual residuals
in the optimality conditions~(\ref{e-dr-oc}), 
evaluated at the approximate primal and dual solution 
$\tilde x^{(k)}$, $v^{(k)}$.

More specifically, using the optimality 
conditions~(\ref{e-proxop-oc}) that characterize 
$\tilde x^{(k)} = \prox_{f/\sigma} (z^{(k-1)})$,
we see that $\tilde x^{(k)}$, $\tilde z^{(k-1)}$ 
satisfy all the optimality conditions for the conic
LPs~(\ref{e-converted-V}),
except two conditions: in general, 
$\tilde x^{(k)} \not\in\mathcal V$  and $v^{(k)}\not\in\mathcal V^\perp$.
The primal and dual residuals measure the errors in these equations.

\paragraph{Stopping condition}
One can also note (from line 2 in the algorithm) that the 
step $G(z^{(k-1)}) = \tilde x^{(k)} - w^{(k)}$ in~(\ref{e-dr-step-z})
can be decomposed as
$G(z^{(k-1)}) 
= -r_{\mathrm p}^{(k)} - (1/\sigma) r_{\mathrm d}^{(k)}$
and since the two terms on the right-hand side are orthogonal,
\BEQ \label{e-error-sum}
\|G(z^{(k-1)})\|_2^2 = 
\|r_{\mathrm p}^{(k)}\|_2^2 + \frac{1}{\sigma^2} 
\|r_{\mathrm d}^{(k)}\|_2^2.
\EEQ
A simple stopping criterion is to terminate when
\BEQ \label{e-stopping}
\frac{\|r_{\mathrm p}^{(k)}\|_2}
{\max\{1.0, \|\tilde x^{(k)}\|_2\}} \leq \epsilon_\mathrm p 
\qquad \mbox{and} \qquad
\frac{\|r_{\mathrm d}^{(k)}\|_2}{\max\{1.0, \|v^{(k)}\|_2\}} 
\leq \epsilon_\mathrm d
\EEQ
for some relative tolerances $\epsilon_\mathrm p$ and 
$\epsilon_\mathrm d$.

\paragraph{Choice of steplength}
In the standard convergence analysis of the Douglas-Rachford algorithm 
the parameter $\sigma$ is assumed to be an arbitrary positive constant  
\cite{EcB:92}. 
However the efficiency in practice is greatly influenced by the 
steplength choice  and several strategies have been proposed for 
varying $\sigma$ during the algorithm \cite{HYW:00,WaL:01,HLW:03}.
As a guideline, it is often observed that the convergence is slow
if one of the two components of $\|G(z^{(k-1)})\|_2$
in~(\ref{e-error-sum}) decreases much more rapidly than the other,
and that adjusting $\sigma$ can help control the balance between 
the primal and dual residuals.
A simple strategy is to take
\BEQ \label{e-adapt-t}
\sigma_{k+1} = \left\{ \begin{array}{ll}
 \sigma_k\tau_k & t_k > \mu \\
 \sigma_k/\tau_k & t_k < 1/\mu \\
 \sigma_k  & \mbox{otherwise,}
 \end{array} \right.
\EEQ
where $t_k$ is the ratio of relative primal and dual residuals,
\[
 t_k = 
 \frac{\|r_{\mathrm p}^{(k)}\|_2} {\|\tilde x^{(k)}\|_2}
 \cdot
 \frac{\|v^{(k)}\|_2} {\|r_{\mathrm d}^{(k)}\|_2},
\]
and $\tau_k$ and $\mu$ are parameters greater than one.
This is further discussed in section~\ref{s-step}.

\paragraph{Projection}
The subspace $\mathcal V$ contains the vectors
$\tilde x = (\tilde x_1, \ldots, \tilde x_l)$ that can
be expressed as $\tilde x_k = E_{\gamma_k}x$ for some $x\in\reals^n$.
The Euclidean projection of a vector $\tilde x$ on 
$\mathcal V$ is therefore easy to compute.
For each $i \in \{1,2,\ldots,n\}$, define
$M(i) = \{k \mid i \in \gamma_k\}$.
The vertices of $T$ indexed by $M(i)$ define a subtree 
(this is a consequence of the running intersection property).
The projection $P_{\mathcal V}(\hat x)$ of 
$\tilde x$ on $\mathcal V$ is the vector 
\[
 P_{\mathcal V}(\tilde x) =
(E_{\gamma_1}(\bar x), E_{\gamma_2}(\bar x), \ldots, 
E_{\gamma_l}(\bar x))
\]
where $\bar x$ is the $n$-vector with components
\[
 \bar x_i = 
 \frac{(\sum_{k\in M(i)} E_{\gamma_k}^T \tilde x_k)_i} {|M(i)|},
 \quad i=1\ldots,n.
\]
In other words, component $i$ of $\bar x$ is a simple average of
the corresponding components of $\tilde x_k$, for 
the sets $\gamma_k$ that contain $i$.

\paragraph{Proximal operator}
The value of the proximal operator 
$\tilde x = \prox_{f/\sigma}(z)$ of $f$, 
applied to a vector $z = (z_1,\ldots, z_l)$ is the solution of the
conic quadratic optimization problem~(\ref{e-proxop}).
An interior-point method applied to this problem
requires the solution of KKT systems
\[
 \left[\begin{array}{cc}
  \sigma I + H & \tilde A^T \\ \tilde A &  0 
 \end{array}\right] \left[\begin{array}{c} 
 \Delta \tilde x \\ \Delta y \end{array}\right]
 = \left[\begin{array}{c} r_{\tilde x} \\ r_y \end{array}\right] 
\]
where $H$ is a block-diagonal positive definite scaling matrix.
As before, we assume that the diagonal blocks of $H$ are of the form
$H_k= \nabla^2\phi_k(w_k)$ where  $\phi_k$ is a logarithmic
barrier of $\mathcal C_k$.  
The cost per iteration of evaluating the proximal operator is dominated
by the cost of assembling the coefficient matrix 
\BEQ \label{e-prox-schur-compl}
\tilde A ( \sigma I + H)^{-1} \tilde A^T 
 = \sum_{k=1}^l \tilde A_k (\sigma I + H_k)^{-1} \tilde A_k^T 
\EEQ
in the Schur complement equation
\[
\tilde A ( \sigma I + H)^{-1} \tilde A^T \Delta y
 = \tilde A(\sigma I+H)^{-1} r_x - r_y   
\]
and the cost of solving the Schur complement system.
For many types of conic LPs the extra term $\sigma I$ 
in (\ref{e-prox-schur-compl})
can be handled by simple changes in the interior-point algorithm.
This is true in particular when $H_k$ is diagonal or 
diagonal-plus-low-rank, as is the case when $\mathcal C_k$
is a nonnegative orthant or second-order cone.
For positive semidefinite cones the modifications are more involved  
and will be discussed in section~\ref{s-prox-sdp}.
In general, it is therefore fair to assume that in most applications
the cost of assembling 
the Schur complement matrix in~(\ref{e-prox-schur-compl}) is roughly the 
same as the cost of computing $\tilde A^T H^{-1}\tilde A^T$.
Since the Schur complement matrix in~(\ref{e-prox-schur-compl})
is sparse (under our assumption of correlative sparsity), 
it can be factored at a smaller cost than its 
counterpart~(\ref{e-converted-schur}) for the reformulated conic LPs.
Depending on the amount of correlative sparsity,
the cost of one evaluation of the proximal operator $\prox_{f/\sigma}$ 
via an interior-point method can therefore be substantially less
than the cost of solving the reformulated problems directly by an
interior-point method. 

\section{Sparse semidefinite optimization} \label{s-sdp}
In the rest of the paper we discuss the application to sparse
semidefinite optimization.  In this section we first explain 
why sparse SDPs with a chordal sparsity pattern can be viewed as
examples of partially separable structure.
In section~\ref{s-decmp-sdp} we then apply the decomposition
method described in section~\ref{s-spingarn}.

We formally define a symmetric sparsity pattern of order $p$
as a set of index pairs
\[
   V\subseteq \{1,2,\ldots, p\} \times \{1,2,\ldots,p\} 
\]
with the property that $(i,j) \in V$ whenever $(j,i)\in V$. 
We say a symmetric matrix $X$ of order $p$ has sparsity pattern $V$ 
if $X_{ij} =0$ when $(i,j) \not\in V$.
The entries $X_{ij}$ for $(i,j)\in V$ are referred to as the 
\emph{nonzero entries} of $X$, even though they may be numerically zero.
The set of symmetric matrices of order $p$ with sparsity pattern $V$
is denoted $\SV{p}{V}$.

\subsection{Nonsymmetric formulation}
Consider a semidefinite program (SDP) in the standard form and 
its dual:
\BEQ \label{e-sdps}
\begin{array}[t]{ll}
\mbox{minimize}   & \Tr(CX) \\
\mbox{subject to} & \Tr(F_iX) = b_i, \quad i=1,\ldots,m \\
& X\succeq 0
\end{array} \qquad\qquad
\begin{array}[t]{ll}
\mbox{maximize} & b^Ty \\
\mbox{subject to} & 
\sum\limits_{i=1}^m y_i F_i + S = C \\
& S\succeq 0.
\end{array}
\EEQ
The primal variable is a symmetric matrix $X\in\symm^p$; the dual 
variables are $y\in\reals^m$ and the slack matrix $S\in\symm^p$.
The problem data are the vector $b\in\reals^m$
and the matrices $C$, $F_i\in\symm^p$.

The \emph{aggregate} sparsity pattern is the union of the sparsity
patterns of $C$, $F_1$, \ldots, $F_m$.
If $V$ is the aggregate sparsity pattern, then we can take
$C$, $F_i \in \SV{p}{V}$.
The dual variable $S$ is then necessarily 
sparse at any dual feasible point, with the same sparsity pattern $V$.
The primal variable $X$, on the other hand, is dense in general, 
but one can note that the cost function and the equality constraints
only depend on the entries of $X$ in the positions of the 
nonzeros of the sparsity pattern $V$.  The other entries of $X$ are 
arbitrary, as long as the matrix is positive semidefinite.
The primal and dual problems can therefore be viewed alternatively as 
conic linear optimization problems with respect to a pair of 
non-self-dual cones:
\BEQ \label{e-matrix-cone-LPs}
\begin{array}[t]{ll}
\mbox{minimize}   & \Tr(CX) \\
\mbox{subject to} & \Tr(F_iX) = b_i, \quad i=1,\ldots,m \\
& X\in\SVc{p}{V}
\end{array} \qquad\qquad
\begin{array}[t]{ll}
\mbox{maximize} & b^Ty \\
\mbox{subject to} & 
\sum\limits_{i=1}^m y_i F_i + S = C \\
& S\in\SVp{p}{V}.
\end{array}
\EEQ
Here the variables $X$ and $S$, as well as the coefficient matrices
$C$, $F_i$, are matrices in $\SV{p}{V}$. 
The primal cone $\SVc{p}{V}$ is the set of matrices in $\SV{p}{V}$
that have a positive semidefinite completion, \ie, the projection
of the cone of positive semidefinite matrices of order $p$ 
on the subspace $\SV{p}{V}$.
We will refer to $\SVc{p}{V}$ as the sparse 
\emph{p.s.d.-completable cone}.
The dual cone $\SVp{p}{V}$ is the set of positive semidefinite
matrices in $\SV{p}{V}$,
\ie, the intersection of the cone of positive semidefinite matrices
of order $p$ with the subspace $\SV{p}{V}$.
This cone will be referred to as the sparse p.s.d.\ cone.
It can be shown that the two cones form a dual pair of proper 
convex cones, provided the nonzero positions in the sparsity 
pattern $V$ include the diagonal entries
(a condition that naturally holds in semidefinite programming).

\paragraph{Vector notation}
It is often convenient to use vector notation for the matrix variables
in~(\ref{e-matrix-cone-LPs}).  For this purpose we introduce an operator
$x = \svec_V(X)$ that maps the lower-triangular nonzeros of 
a matrix $X\in\SV{p}{V}$ to a vector $x$ of length $n = (|V|+p)/2$, 
using a format that preserves inner products, \ie, 
$\Tr(XY) = \svec_V(X)^T \svec_V(Y)$ for all $X$, $Y$.
For example, one can copy the nonzero lower-triangular entries 
of $X$ in column-major order to $x$, scaling the strictly 
lower-triangular entries by $\sqrt 2$.
A similar notation $x=\svec(X)$ (without subscript)
will be used for a packed vector representation of a dense matrix:
if $X\in\symm^p$, then $x=\svec(X)$ is a vector of length $p(p+1)/2$
containing the lower-triangular entries of $X$ in a storage format
that preserves the inner products.
Using this notation, the matrix cones $\SVc{p}{V}$ and $\SVp{p}{V}$ 
can be `vectorized' to define two cones 
\[
\mathcal C = \{\svec_V(X) \mid X \in \SVc{p}{V}\}, \qquad
\mathcal C^* = \{\svec_V(S) \mid S \in \SVp{p}{V}\}.
\]
These cones form a dual pair of proper convex cones in $\reals^n$ 
with $n = (|V|+p)/2$. 
The conic linear optimization
problems~(\ref{e-matrix-cone-LPs}) can then be written as~(\ref{e-conic})
with variables $x=\svec_V(X)$, $s= \svec_V(S)$, $y$, 
and problem parameters 
\[
c = \svec_V(C),  \qquad
 A= \left[\begin{array}{cccc}
 \svec_V(F_1) & \svec_V(F_2) & \cdots & \svec_V(F_m) \end{array}\right]^T.
\]

\subsection{Clique decomposition of chordal sparse matrix cones}
\label{s-clique-decmp}
The nonsymmetric conic optimization or matrix completion approach 
to sparse semidefinite programming, based on the 
formulation~(\ref{e-matrix-cone-LPs}), was first proposed
by Fukuda \emph{et al.} \cite{FKMN:00} and further developed in 
\cite{NFFKM:03,Bur:03,SrV:04a,ADV:10,KKMY:11}.
The various techniques described in these papers all assume that the 
sparsity pattern $V$ is \emph{chordal}.  
In this section we review some key results concerning positive 
semidefinite matrices with chordal sparsity patterns.

With each sparsity pattern $V$ one associates an undirected graph
$\mathcal G_V$ with $p$ vertices and edges $\{i,j\}$  between pairs of 
vertices $(i,j) \in V$ with $i>j$.
A clique in $\mathcal G_V$ is a maximal complete subgraph, \ie, a 
maximal set  $\beta \subseteq \{1,2,\ldots,p\}$ such that 
$\beta \times \beta \subseteq V$.
Each clique defines a maximal dense principal submatrix in any matrix
with sparsity pattern $V$.
If the cliques in the graph $\mathcal G_V$ are $\beta_k$, $k=1,\ldots,l$,
then the sparsity pattern $V$ can be expressed as 
$V = \bigcup_{k=1,\ldots,l} \beta_k \times \beta_k$.
A sparsity pattern $V$ is called chordal if the graph 
$\mathcal G_V$ is chordal. 

In the remainder of the paper we assume that $V$ is a chordal sparsity
pattern that contains all the diagonal entries ($(i,i) \in V$ for
$i=1,\ldots,p$).
We denote by $\beta_k$, $k=1,\ldots,l$, the cliques of $\mathcal G_V$ and 
define $V_k = \beta_k \times \beta_k$.
We will make use of two classical theorems that characterize
the matrix cones $\SVc{p}{V}$ and $\SVp{p}{V}$ for chordal patterns $V$.
These theorems are discussed in the next two paragraphs.

\paragraph{Decomposition of sparse positive semidefinite cone}
The first theorem \cite[theorem 2.3]{AHMR:88}
states that the sparse p.s.d.\ cone $\SVp{p}{V}$ 
is a sum of positive semidefinite cones
with simple sparsity patterns:
\BEQ 
\SVp{p}{V} = \sum_{k=1}^l \SVp{p}{V_k} 
= \{
\sum_{k=1}^l \mathcal E^*_{\beta_k}(\tilde S_k)  \mid
 \tilde S_k \in \symm^{|\beta_k|}_+\}
\label{e-clique-decomp-1}
\EEQ
where $\symm^{|\beta_k|}_+$ is the positive semidefinite cone of order 
$|\beta_k|$.
The operator $\mathcal E^*_{\beta_k}$ copies a dense matrix of
order $|\beta_k|$ to the principal submatrix indexed by $\beta_k$
in a symmetric matrix of order $p$; see section~\ref{s-notation}.
According to the decomposition result~(\ref{e-clique-decomp-1}),
every positive semidefinite matrix $X$ with sparsity pattern $V$ 
can be decomposed as a sum of positive semidefinite matrices, 
each with a sparsity pattern consisting of a single principal dense 
block $V_k = \beta_k\times\beta_k$.  
If $X$ is positive definite, a decomposition of this form is easily 
calculated via a zero-fill Cholesky factorization.

\paragraph{Decomposition of positive-semidefinite-completable cone}
The second theorem characterizes the p.s.d.-completable cone $\SVc{p}{V}$ 
\cite[theorem 7]{GJSW:84}:
\BEA 
\SVc{p}{V} & = & \{X \in\SV{p}{V} \mid X_{\beta_k\beta_k}
\succeq 0, \; k=1,\ldots,l\} \nonumber \\
& = & \{X \in\SV{p}{V} \mid \mathcal E_{\beta_k}(X)
\in \symm^{|\beta_k|}_+, \; k=1,\ldots,l\}.
\label{e-clique-decomp-2}
\EEA
The operator $\mathcal E_{\beta_k}$ extracts from its argument
the dense principal submatrix indexed by $\beta_k$. 
(This is the adjoint operation of $\mathcal E_{\beta_k}^*$;
see section~\ref{s-notation}.)
In other words, a matrix in $\SV{p}{V}$ has a positive semidefinite 
completion if and only if all its maximal dense principal submatrices 
$X_{\beta_k\beta_k}$ are positive semidefinite.
This result can be derived from the characterization of 
$\SVp{p}{V}$ in~(\ref{e-clique-decomp-1}) and the fact that the cones
$\SVc{p}{V}$ and $\SVp{p}{V}$ are duals. 

\paragraph{Clique decomposition in vector notation}
We now express the clique decomposition 
formulas~(\ref{e-clique-decomp-1}) 
and~(\ref{e-clique-decomp-2})  in vector notation.
For each clique $\beta_k$, define an index set 
$\gamma_k \subseteq\{1,2,\ldots,n\}$ via the identity 
\BEQ \label{e-gamma-beta}
 E_{\gamma_k} \svec_V(Z) = \svec(Z_{\beta_k\beta_k})
 \quad\forall Z\in\SV{p}{V}.
\EEQ
The index set $\gamma_k$ has length 
$|\gamma_k| = |\beta_k| (|\beta_k|+1)/2$ and its elements 
indicate the positions of the entries of the 
$\beta_k\times \beta_k$ submatrix of $Z$
in the vectorized matrix $\svec_V(Z)$.
Using this notation, the cone $\mathcal C$ can be expressed 
as~(\ref{e-C-primal}) where 
$\mathcal C_k = \{\svec(W) \mid W\in\symm^{|\beta_k|}_+\}$
is the vectorized dense positive semidefinite matrix cone of order 
$|\beta_k|$.
The clique decomposition~(\ref{e-clique-decomp-2}) 
of the p.s.d.\ cone can be expressed in vector notation 
as~(\ref{e-C-dual}).
(Note that $\mathcal C_k$ is self-dual, so here
$\mathcal C_k = \mathcal C_k^*$.)
The decomposition result~(\ref{e-C-primal})
shows that the p.s.d.-completable cone associated with
a chordal sparsity pattern $V$ is partially separable. 

\paragraph{Clique tree}
The cliques $\beta_k$ of $V$ can be arranged in a clique tree that 
satisfies the running intersection property 
($\beta_i \cap \beta_j \subseteq
\beta_k$ if clique $k$ is on the path between cliques $\beta_i$
and $\beta_j$ in the tree); see \cite{BlP:93}.
We denote by $\eta_k$ the intersection of the clique $\beta_k$ 
with its parent in the clique tree.  

Since there is a one-to-one relation between the index sets 
$\gamma_k$ defined in~(\ref{e-gamma-beta}) and the cliques $\beta_k$ 
of $\mathcal G_V$, we can identify the 
clique graph of $\mathcal G_V$ (which has vertices $\beta_k$)
with the intersection graph for the index sets $\gamma_k$.
Similarly, we do not have to distinguish between a clique tree $T$ for 
$\mathcal G_V$ and a spanning tree with the running intersection 
property in the intersection graph of the sets $\gamma_k$.
The sets $\alpha_k = \gamma_k \cap \prnt(\gamma_k)$
are in a one-to-one relation to the 
sets $\eta_k = \beta_k \cap \prnt(\beta_k)$
via the identity $E_{\alpha_k} (\svec_V(Z)) = \svec(Z_{\eta_k\eta_k})$
for arbitrary $Z\in\SV{p}{V}$.

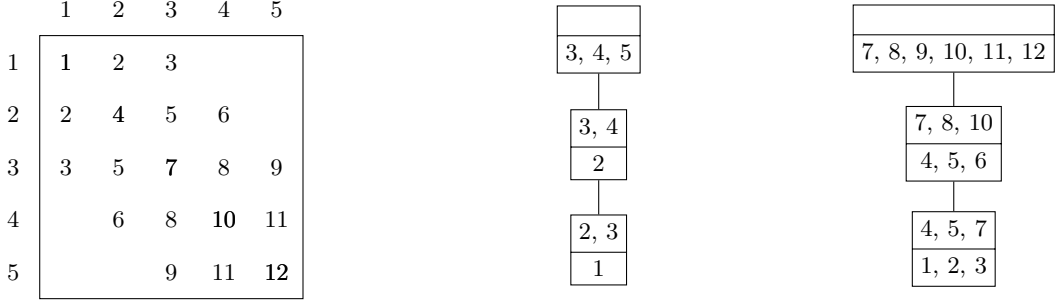
\begin{figure}
\begin{minipage}{.4\linewidth}
\hspace*{\fill}
\begin{tikzpicture}[scale=0.7]
   \tikzset{VertexStyle/.style = {shape = circle, minimum size = 4pt, 
       inner sep = 0 pt, fill=black}}
   \tikzset{VertexStyleR/.style = {shape = circle, minimum size = 4pt, 
       inner sep = 0 pt, fill=red}}

   \foreach \i/\j/\k in {
       0/0/1, 1/0/2, 2/0/3, 
       1/1/4 , 2/1/5, 3/1/6, 
       2/2/7, 3/2/8, 4/2/9,
       3/3/10, 4/3/11,
       4/4/12
      }
      {
      \node[font=\footnotesize] at (\j,-\i){\k};
      \node[font=\footnotesize] at (\i,-\j){\k};
      }
   \foreach \k/\i in { 1/0, 2/1, 3/2, 4/3, 5/4}
      {
      \node[font=\footnotesize] at (-1,-\i){\k};
      \node[font=\footnotesize] at (\i,1){\k};
      }

   \draw (-0.5,0.5) rectangle (4.5,-4.5);
\end{tikzpicture}
\hspace*{\fill}
\end{minipage}
\hspace*{\fill}
\begin{minipage}{.25\linewidth}
\hspace*{\fill}
\begin{tikzpicture}[scale=1.4, 
clique/.style = {rectangle split, rectangle split parts = 2, draw}]
 \footnotesize
 \node[clique](1) at (1,0){ 2, 3 \nodepart{second} 1 };
 \node[clique](2) at (1,1){ 3, 4 \nodepart{second} 2 };
 \node[clique](3) at (1,2){ \nodepart{second} 3, 4, 5 };
 \draw (1) -- (2);
 \draw (2) -- (3);
\end{tikzpicture}
\hspace*{\fill}
\end{minipage}
\hspace*{\fill}
\begin{minipage}{.25\linewidth}
\hspace*{\fill}
\begin{tikzpicture}[scale=1.4, 
clique/.style = {rectangle split, rectangle split parts = 2, draw}]
 \footnotesize
 \node[clique](1) at (1,0){ 4, 5, 7 \nodepart{second} 1, 2, 3 };
 \node[clique](2) at (1,1){ 7, 8, 10 \nodepart{second} 4, 5, 6 };
 \node[clique](3) at (1,2){ \nodepart{second} 7, 8, 9, 10, 11, 12 };
 \draw (1) -- (2);
 \draw (2) -- (3);
\end{tikzpicture}
\hspace*{\fill}
\end{minipage}
\hspace*{\fill}
\caption{A $5\times 5$ chordal sparsity pattern with 12 nonzero
entries in the lower triangular part.
The numbers in the matrix are the indices of the entries in the 
vectorized matrix.
The center of the figure shows a clique tree. 
The right-hand part of the figure shows the corresponding spanning tree
in the intersection graph. }
\label{f-beta-gamma}
\end{figure}

The notation is illustrated in Figure~\ref{f-beta-gamma} for 
a simple example. 
There are three cliques
\[
\beta_1 = \{1,2,3\}, \qquad \beta_2 = \{2,3,4\}, \qquad 
\beta_3 = \{3,4,5\}.
\]
If we use the column-major order for the nonzero entries in the 
vectorized matrix, these cliques correspond to the index sets
\[
\gamma_1 = \{1,2,3,4,5,7\}, \qquad 
\gamma_2 = \{4,5,6,7,8,10\}, \qquad
\gamma_3 = \{7,8,9,10,11,12\}.
\]
The sets $\eta_k = \beta_k \cap \prnt(\beta_k)$ and 
$\alpha_k = \gamma_k \cap \prnt(\gamma_k)$ are 
\[
\eta_1 = \{2,3\}, \qquad 
\eta_2 = \{3,4\}, \qquad 
\eta_3 = \{\},  \qquad
\alpha_1 = \{4,5,7\}, \qquad
\alpha_2 = \{7,8,10\}, \qquad
\alpha_3 = \{\}. 
\]

\section{Decomposition in semidefinite programming} \label{s-decmp-sdp}
We now work out the details of the decomposition method when 
applied to sparse semidefinite programming. 
In particular, we describe an efficient method for solving the
quadratic conic optimization problem~(\ref{e-proxop}), needed for
the evaluation of the proximal operator, when the cone
$\tilde{\mathcal C}$ is a product of positive semidefinite matrix cones.

\subsection{Converted problems}
We first express the reformulated problems~(\ref{e-converted-V}) 
and~(\ref{e-converted-matrix}) for SDPs in matrix notation.
The reformulated primal problem can be written as
\BEQ \label{e-sdp-converted-primal}
\begin{array}{ll}
\mbox{minimize} & 
\sum\limits_{k=1}^l \Tr(\tilde C_k \tilde X_k) \\*[1ex]
\mbox{subject to} & 
\sum\limits_{k=1}^l \Tr(\tilde F_{ik} \tilde X_k)
= b_i, \quad i=1,\ldots,m \\*[1ex]
 & \mathcal E_{\eta_j} (\mathcal E_{\beta_k}^*(\tilde X_k)
   - \mathcal E_{\beta_j}^*(\tilde X_j)) = 0, \quad
   k=1,\ldots, l, \quad \beta_j\in\ch(\beta_k) \\*[1ex]
 & \tilde X_k \succeq 0, \quad k=1,\ldots, l
\end{array}
\EEQ
with variables $\tilde X_k \in \symm^{|\beta_k|}$, $k=1,\ldots,l$.
The coefficient matrices $\tilde C_k$ and $\tilde F_{ik}$ are chosen
so that
\BEQ \label{e-converted-CF}
 \Tr(CZ) = \sum_{k=1}^l \Tr(\tilde C_k Z_{\beta_k\beta_k}), \qquad
 \Tr(F_iZ) = \sum_{k=1}^l \Tr(\tilde F_{ik} Z_{\beta_k\beta_k}) 
 \qquad \forall Z\in \SV{p}{V}.  
\EEQ
One possible choice is
\BEQ \label{e-Fik}
 \tilde C_k = \mathcal E_{\beta_k}(C - \mathcal P_{\eta_k}(C)),
 \qquad
 \tilde F_{ik} = \mathcal E_{\beta_k}(F_i - \mathcal P_{\eta_k}(F_i)).
\EEQ
The converted dual problem is
\BEQ \label{e-sdp-converted-dual}
\begin{array}{ll}
\mbox{maximize} & b^T y \\
\mbox{subject to} & \sum\limits_{i=1}^m y_i \tilde F_{ik}
  + \mathcal E_{\beta_k} (\mathcal E_{\eta_k}^*(U_k) -
 \sum\limits_{\beta_j\in\ch(\beta_k)} \mathcal E_{\eta_j}^*(U_j))
  + \tilde S_k = \tilde C_k, \quad k=1,\ldots, l \\
& \tilde S_k \succeq 0,\quad k=1,\ldots,l.
\end{array}
\EEQ
with variables $y$, $\tilde S_k\in\symm^{|\beta_k|}$,
and $U_k\in \symm^{|\eta_k|}$, $k=1,\ldots,l$.
The reformulations~(\ref{e-sdp-converted-primal}) 
and~(\ref{e-sdp-converted-dual})
also follow from the clique-tree 
conversion methods proposed in \cite{KKMY:11,FKMN:00}. 

The variables $\tilde X_k$ in~(\ref{e-sdp-converted-primal})
are interpreted as copies 
of the dense submatrices $X_{\beta_k\beta_k}$.
The second set of equality constraints in~(\ref{e-sdp-converted-primal})
are the consistency constraints that ensure that the
entries of $\tilde X_k$ agree when they refer to the same entry of $X$.  
The consistency equations can be written in simpler form
if we assume that the indices are sorted so that
indices in $\beta_k\setminus\eta_k$ precede those in $\eta_k$.
(This is the case if the indices are sorted using a perfect elimination
ordering for the Cholesky factorization with chordal sparsity pattern $V$.)
If we partition $\tilde X_k$ and $X_{\beta_k\beta_k}$ conformably as
\[
\tilde X_k 
= \left[\begin{array}{cc}
 \tilde X_{k,11} & \tilde X_{k,21}^T \\
 \tilde X_{k,21} & \tilde X_{k,22} \end{array}\right], \qquad
X_{\beta_k\beta_k} 
= \left[\begin{array}{cc}
 X_{\beta_k\setminus\eta_k,\beta_k\setminus\eta_k} & 
 X_{\beta_k\setminus\eta_k,\eta_k} \\
 X_{\eta_k,\beta_k\setminus\eta_k} & X_{\eta_k\eta_k} 
 \end{array}\right]
\]
then the consistency equations reduce to
\[
\tilde X_{j,22} - 
\mathcal E_{\eta_j} (\mathcal E_{\beta_k}^*(\tilde X_k)) = 0,
\quad k=1,\ldots,l, \quad \beta_j\in\ch(\beta_k).
\]
Similarly, the definitions~(\ref{e-converted-CF}) simplify as
\[
 \tilde C_k 
 = \left[\begin{array}{cc}
   C_{\beta_k\setminus\eta_k,\beta_k\setminus\eta_k} & 
   C_{\beta_k\setminus\eta_k,\eta_k} \\
   C_{\eta_k,\beta_k\setminus\eta_k} & 0 \end{array}\right], \qquad
 \tilde F_{ik} = 
  \left[\begin{array}{cc}
   (F_i)_{\beta_k\setminus\eta_k,\beta_k\setminus\eta_k} & 
   (F_i)_{\beta_k\setminus\eta_k,\eta_k} \\
   (F_i)_{\eta_k,\beta_k\setminus\eta_k} & 0 \end{array}\right].
\]
We can also note that the matrices $U_k$ in the dual problem 
play an identical role as  the update matrices in a 
multifrontal supernodal Cholesky factorization \cite{ADV:12}.

\subsection{Proximal operator} \label{s-prox-sdp}
In the clique tree conversion methods of \cite{NFFKM:03,KKMY:11}
the converted SDP~(\ref{e-sdp-converted-primal}) 
is solved by an interior-point method.
A limitation to this approach is the large number of equality constraints
added in the primal problem or, equivalently, the large dimension
of the auxiliary variables $U_k$ in the dual problem.
In section~\ref{s-spingarn} we proposed an operator-splitting
method to address this problem. 
The key step in each iteration of the splitting method is the evaluation 
of a proximal operator, by solving the quadratic conic optimization
problem (QP)
\BEQ \label{e-sdp-prox-qp}
\begin{array}{ll}
\mbox{minimize} & \sum\limits_{k=1}^l \Tr(\tilde C_k \tilde X_k) 
 + (\sigma/2) \sum\limits_{k=1}^l \|\tilde X_k - Z_k\|_F^2 \\
\mbox{subject to} & 
\sum\limits_{k=1}^l \Tr(\tilde F_{ik} \tilde X_k)
= b_i, \quad i=1,\ldots,m \\*[2ex]
 & \tilde X_k \succeq 0, \quad k=1,\ldots, l.
\end{array}
\EEQ
Solving this problem by a general-purpose solver can be quite expensive 
and most solvers require a reformulation to remove the quadratic term 
in the objective by adding second-order cone constraints.  
However the problem can be solved efficiently via a 
customized interior-point solver, as we now describe.
A similar technique was used for handling variable bounds 
in SDPs in~\cite{NWV:08,TTT:07}.

The Newton equation or KKT system that must be solved in each iteration 
of an interior-point method for the conic QP~(\ref{e-sdp-prox-qp}) has
the form
\BEA
\label{e-prox-kkt-1}
\sigma  \Delta\tilde X_k + W_k \Delta\tilde X_k W_k + 
 \sum\limits_{i=1}^m \Delta y_i \tilde F_{ik} & = & R_k, 
 \quad k=1,\ldots, l \\
\sum_{k=1}^l \Tr(\tilde F_{ik}\,\Delta X_k) & = & r_i, \quad 
i=1,\ldots,m,
\label{e-prox-kkt-2}
\EEA
with variables $\Delta \tilde X_k$, $\Delta y$,
where $W_k$ is a positive definite scaling matrix.
The first term $\sigma \Delta \tilde X_k$ results from the quadratic
term in the objective.  Without this term it is straightforward
to eliminate the variable $\Delta \tilde X_k$ from first equation, 
to obtain an equation in the variable $\Delta y$.
To achieve the same goal at a similar cost with a customized solver we
first compute eigenvalue decompositions 
$W_k = Q_k \diag(\lambda_k) Q_k^T$ of the $l$ scaling matrices,
and define $l$ matrices $S_k \in\symm^{|\beta_k|}$ with entries
\[
  (S_k)_{ij} = \frac{1}{\sigma + \lambda_{ki}\lambda_{kj}}, \quad
 i,j =1, \ldots, |\beta_k|.
\] 
We can now use the first equation in~(\ref{e-prox-kkt-1}) to 
express $\Delta \tilde X_k$ in terms of $\Delta  y$:
\[
\Delta\tilde X_k
= Q_k \,
(S \circ (\widehat R_k  - \sum_{i=1}^m \Delta y_i \widehat F_{ik})) \,
Q_k^T
\]
with $\widehat R_k = Q_k^T R_k Q_k$, 
$\widehat F_{ik} = Q_k^T \tilde F_{ik} Q_k$, 
and where $\circ$ denotes the Hadamard (component-wise) product.
Substituting the expression for $\Delta \tilde X_k$ in the second 
equation of~(\ref{e-prox-kkt-2}) gives an equation $H\Delta y = g$
with
\BEQ \label{e-sdp-Hdef}
H_{ij} = \sum_{k=1}^l
\Tr (\widehat F_{ik} (S \circ \widehat F_{jk})), \qquad
g_i = r_i - \sum_{k=1}^l \Tr(\widehat F_{ik} (S\circ R_k)).
\qquad i,j = 1,\ldots,m.
\EEQ
The cost of this solution  method for the KKT 
system~(\ref{e-prox-kkt-1})--(\ref{e-prox-kkt-2}) is comparable
to the cost of solving the KKT systems in an interior-point method
applied to the conic optimization problem~(\ref{e-sdp-prox-qp}) 
without the quadratic term. 
The proximal operator can therefore be evaluated
at roughly the same cost as the cost of solving the converted
SDP~(\ref{e-sdp-converted-primal}) 
with the consistency constraints removed. 

\begin{figure}
\begin{center}
\begin{psfrags}
\psfrag{t}{}
\psfrag{x}[t][t]{matrix order ($p=|\beta_1|$)}
\psfrag{z}[b][b]{time (sec)}
\includegraphics[width=.5\linewidth]{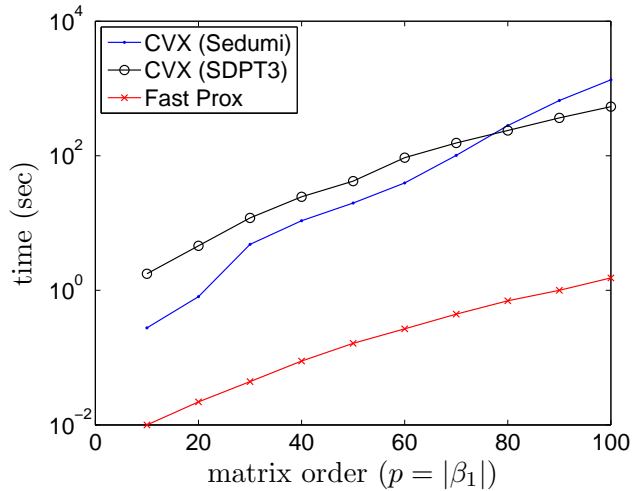}
\end{psfrags}
\caption{Time required for for a single proximal operator 
evaluation~(\ref{e-sdp-prox-qp}) on a dense subproblem with 
a single clique ($l=1$) of size $p=|\beta_1|$ and $m=p$ constraints (averaged over 10 trials). 
The CPU time of the general-purpose solvers SDPT3 and SEDUMI,
called via CVX,  is compared against a customized fast proximal 
operator.}
\label{f-proxspeed}
\end{center}
\end{figure}

To illustrate the value of this technique, we compare
in Figure~\ref{f-proxspeed} the time needed to solve
the semidefinite QP~(\ref{e-sdp-prox-qp}) using three 
methods: SEDUMI and SDPT3 called via 
CVX (version 2.0 beta)~\cite{GRB:07,GRB:08} in MATLAB ,
and an implementation of the algorithm described above
in CVXOPT~\cite{ADV:12c}.
The problems are dense and randomly generated with $l=1$ and
$m = p = |\beta_1|$.  
The figure shows CPU time versus the order $p$ of the matrix variable. 
(For details on the computing environment, see the beginning of
section~\ref{s-experiments}.)

\subsection{Correlative sparsity} \label{s-correlative-sdp}
The efficiency of the decomposition method depends crucially
on the cost of the proximal operator evaluations, which is
determined by the sparsity pattern of the Schur complement
matrix $H$~(\ref{e-sdp-Hdef}), \ie, the correlative
sparsity pattern of the reformulated problems.
Note that in general the scaled matrices $\widehat F_{ik}$ used
to assemble $H$ will be either completely dense 
(if $\tilde F_{ik} \neq 0$) or zero (if $\tilde F_{ik} =  0$).
Therefore $H_{ij} = 0$ if for each $k$ at least one of the coefficient 
matrices $\tilde F_{ik}$ and $\tilde F_{jk}$ is zero.
This rule characterizes the correlative sparsity pattern.

As pointed out in section~\ref{s-corr-sparse}, the correlative
sparsity can be enhanced by exploiting the flexibility in the
choice of parameters of the reformulated problem
(the matrices $\tilde F_{ik}$).  
The definition~(\ref{e-Fik}) is one possible choice, 
but any set of matrices that satisfy~(\ref{e-sdp-converted-dual})
can be used instead.
While the optimal choice is not clear in general, it is straightforward 
in the important special case when the index set $\{1,\ldots,m\}$ can 
be partitioned in $l$ sets $\nu_1$, \ldots, $\nu_l$, with the property
that if $i\in\nu_j$, then all the nonzero entries of $F_i$ belong to 
the principal submatrix $(F_i)_{\beta_j\beta_j}$.
In other words $F_i = \mathcal P_{\beta_j}(F_i)$ for
$i\in\nu_j$. 
In this case, a valid choice for the coefficient matrices $\tilde
F_{ik}$ is to take
\[
\tilde F_{ij}  = \mathcal E_{\beta_j}(F_i), \qquad
\tilde F_{ik}  = 0, \quad k\neq j,
\]
when $i\in\nu_j$.
With this choice, the matrix $H$ can be re-ordered to be block-diagonal 
with dense blocks $H_{\nu_i\nu_i}$.  
Moreover the QP~(\ref{e-sdp-prox-qp})
is separable and equivalent to $l$ independent subproblems
\[
\begin{array}{ll}
\mbox{minimize} & \Tr(C_k\tilde X_k) + 
 (\sigma/2) \|\tilde X_k - Z_k\|_F^2 \\
\mbox{subject to} 
 & \Tr(\tilde F_{ik} \tilde X_k) = b_i, \quad i \in \nu_k \\
 & \tilde X_k \succeq 0.
\end{array}
\]

\section{Numerical examples} \label{s-experiments}

In this section we present the results of numerical experiments with the
decomposition method applied to semidefinite programs.
First, we describe how steplength selection can significantly affect 
(and impair) convergence speed and show how a simple adaptive steplength 
scheme can make the method more robust. 
Then, we apply the decomposition method to an approximate Euclidean 
distance matrix completion problem, motivated by an application
in sensor network node localization, and illustrate the convergence 
behavior of the method in practice. 
The problem involves a sparse matrix variable whose sparsity 
pattern is characterized by the sensor network topology, 
and is interesting because in the converted form the problem 
has block-diagonal correlative sparsity regardless of the
network topology. 
Finally, we present extensive runtime results for a family of
problems with block-arrow aggregate sparsity and 
block-diagonal correlative sparsity. 
By comparing the CPU times required by general-purpose interior-point 
methods and the decomposition method, 
we are able to characterize the regime in which each method is 
more efficient. 

The decomposition method is implemented in Python (version 2.6.5), using the conic 
quadratic optimization solver of CVXOPT (version 1.1.5)~\cite{ADLV:12} 
for solving the conic QPs~(\ref{e-sdp-prox-qp})
in the evaluation of the proximal operators.
SEDUMI (version 1.1) \cite{Stu:99} and SDPT3 (version 4.0) in MATLAB (version 2011 b) are used as the 
general-purpose solver for the experiments in 
sections~\ref{s-blck-arrow} and~\ref{s-edm}.
The experiments are performed on 
an Intel Xeon CPU E31225 processor 
(4 cores, 3.10 GHz clock speed) and 8 GB RAM, 
running Ubuntu 10.04 (Lucid).

\subsection{Adaptive steplength selection} \label{s-step}
The first experiment illustrates the effect of the choice of 
the steplength parameter $\sigma_k$ and explains the motivation behind the 
adaptive strategy~(\ref{e-adapt-t}).
We pick a randomly generated SDP with a block-banded sparsity pattern 
$V$ of order $p=402$ with $l=50$ cliques of size $|\beta_k| = 10$.  
The cliques correspond to overlapping diagonal blocks of order $10$, 
with overlap of size $2$.  
The correlative sparsity pattern in the converted SDP has a block-arrow
structure with 50 diagonal blocks of size $10\times 10$
and 10 dense rows and columns at the end. 
The number of primal constraints and dual variables is $m=510$.

\begin{figure}
\begin{center}
\includegraphics[width=.32\linewidth]{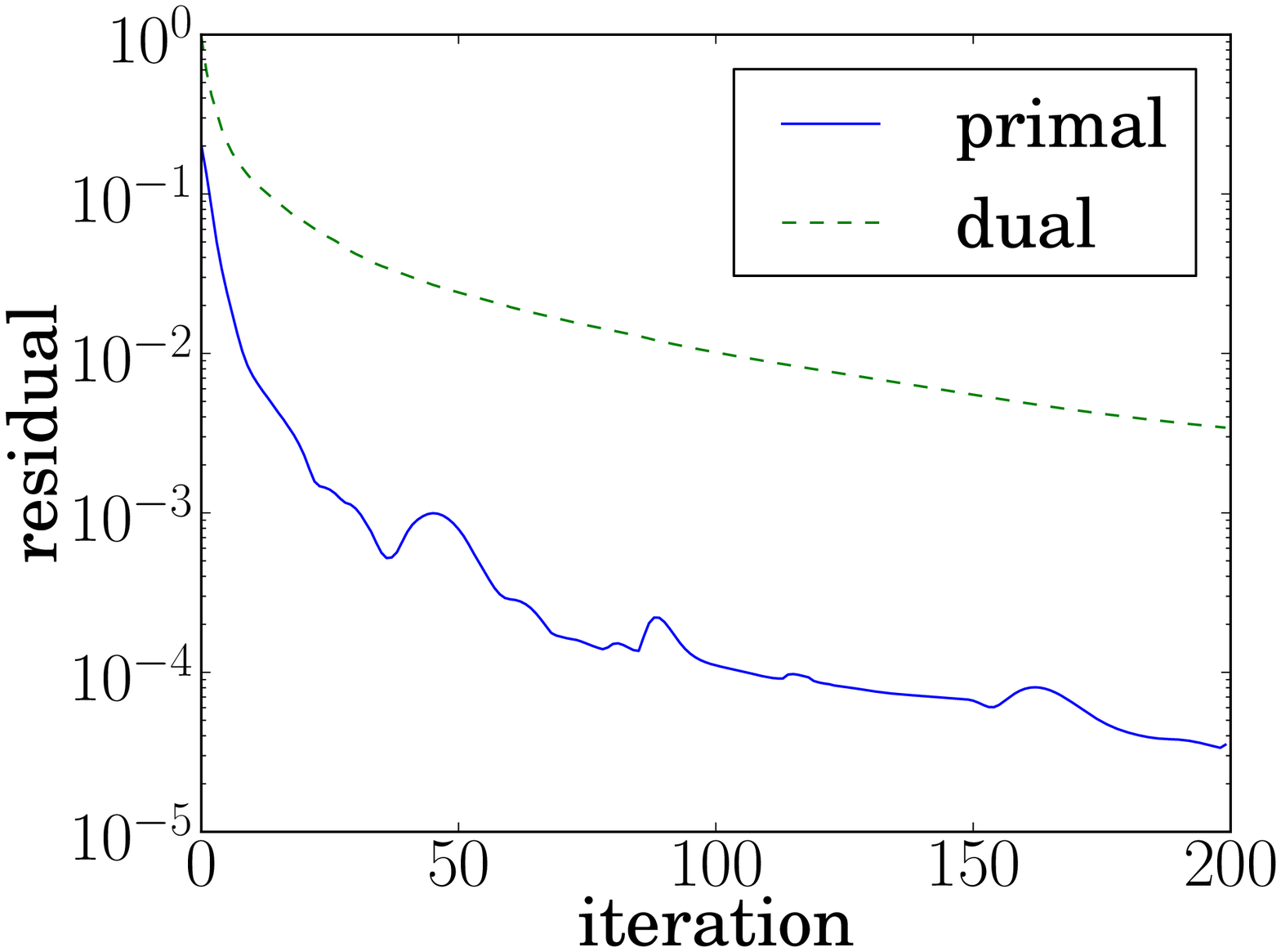}
\hspace*{\fill}
\includegraphics[width=.32\linewidth]{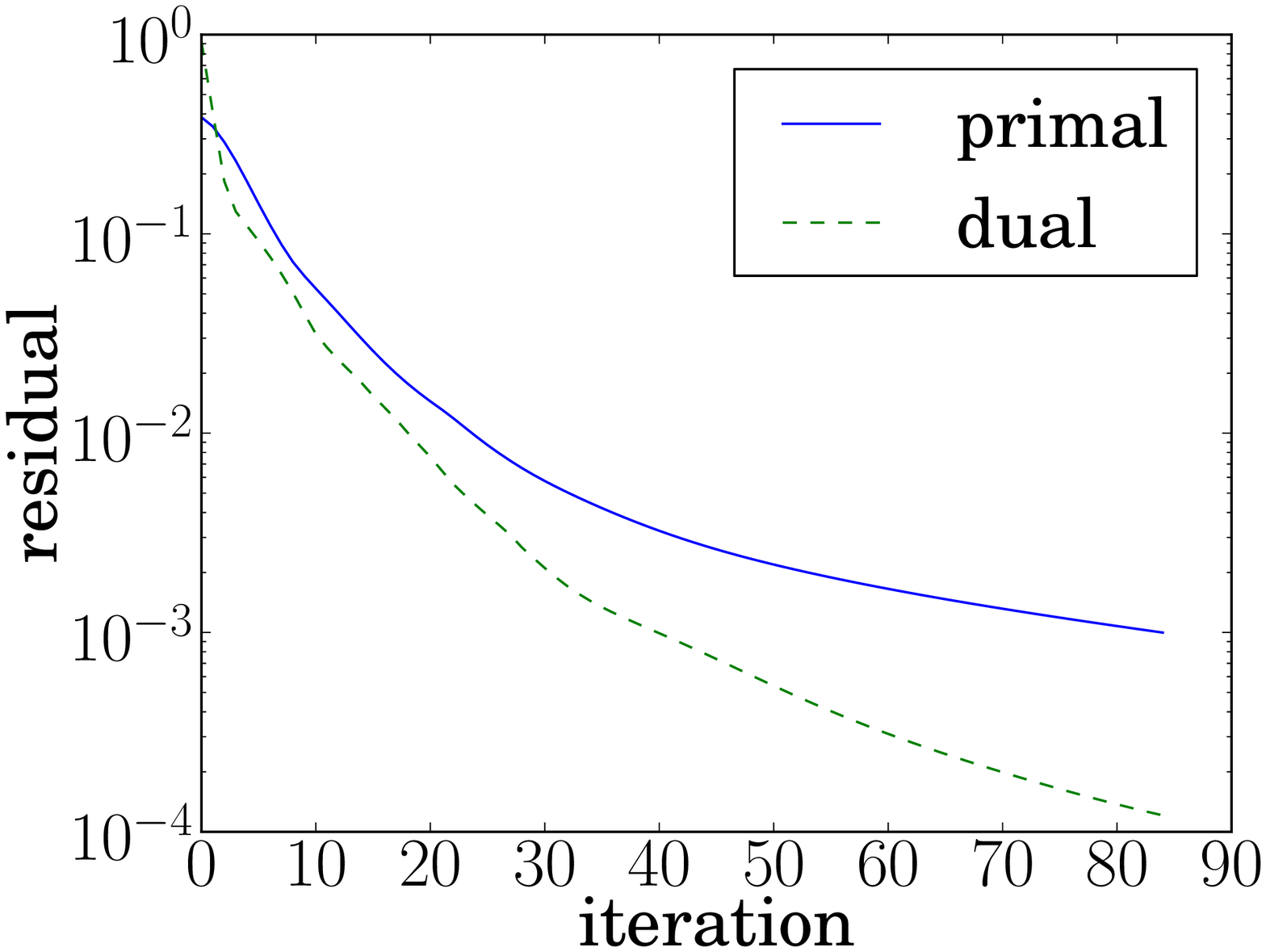}
\hspace*{\fill}
\includegraphics[width=.32\linewidth]{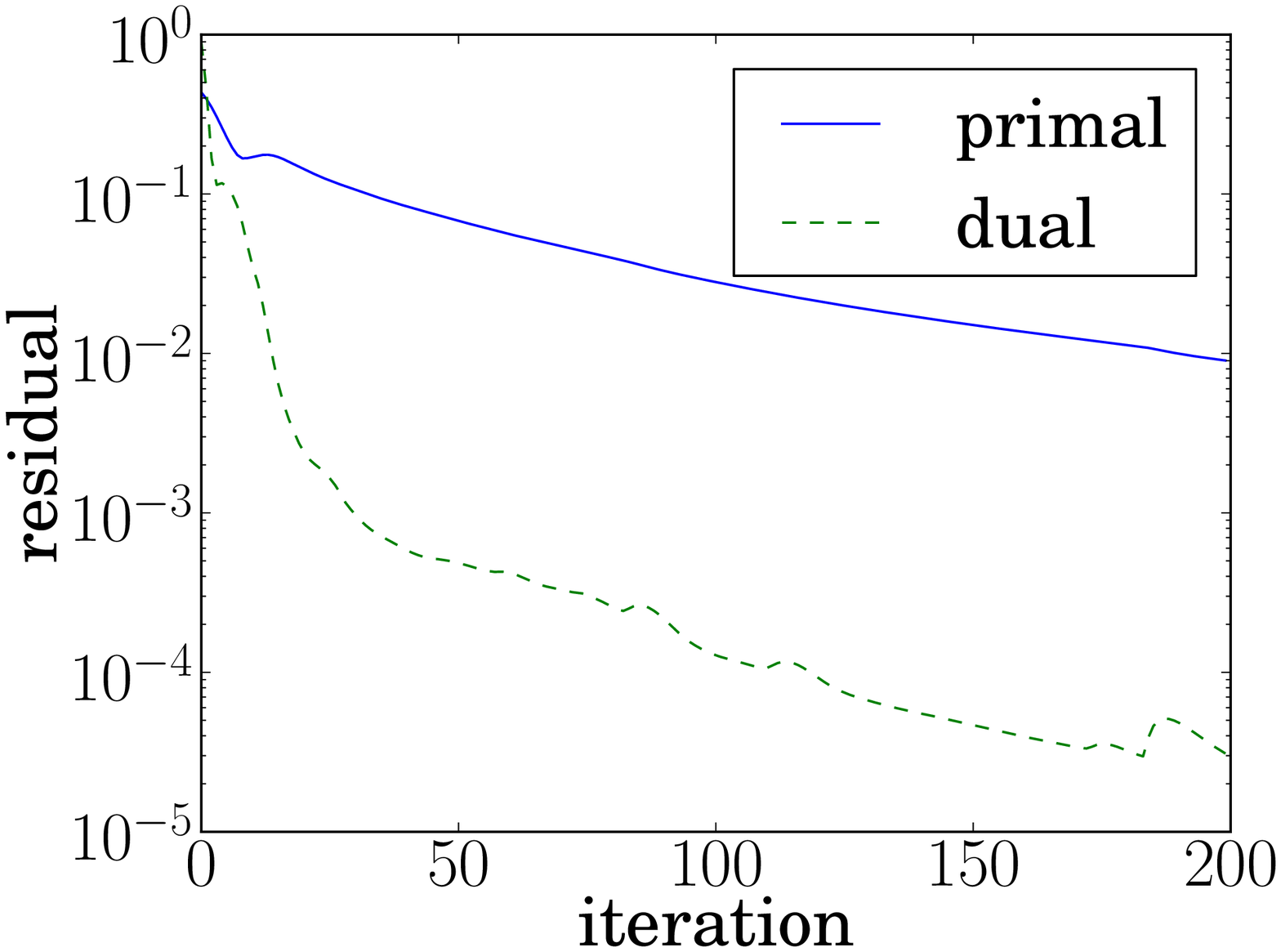}

\caption{Primal residual 
$\|r_{\mathrm p}^{(k)}\|_2/\|\tilde x^{(k)}\|_2$ 
and dual residual
$\|r_{\mathrm d}^{(k)}\|_2/\|v^{(k)}\|_2$ 
versus iteration number $k$ for 
three constant values of $\sigma_k$: 
$\sigma_k = 0.1$ (left), 
$\sigma_k = 0.01$ (middle), and  $\sigma_k = 0.001$ (right).} 
\label{f-sigma-res}
\end{center}
\end{figure}

Figure~\ref{f-sigma-res} shows the primal and dual residuals
$\|r_{\mathrm p}^{(k)}\|_2/\|\tilde x^{(k)}\|_2$ and
$\|r_{\mathrm d}^{(k)}\|_2/\|v^{(k)}\|_2$ 
for three constant values of the steplength parameter: $\sigma_k = 0.1$, 
$0.01$, and $0.001$.
As can be seen, the choice of $\sigma_k$ has a strong effect on
the speed of convergence.  
The figures suggest that when $\sigma_k$  is too large
(the steplength $1/\sigma_k$ is too small)
the dual residual decreases more slowly than the primal residual,
and when $\sigma_k$ is too small, the primal residual decreases more 
slowly.  For a good value of $\sigma_k$ in between, the 
two residuals decrease at about the same rate.

\begin{figure}
\begin{center}
\includegraphics[width=.32\linewidth]{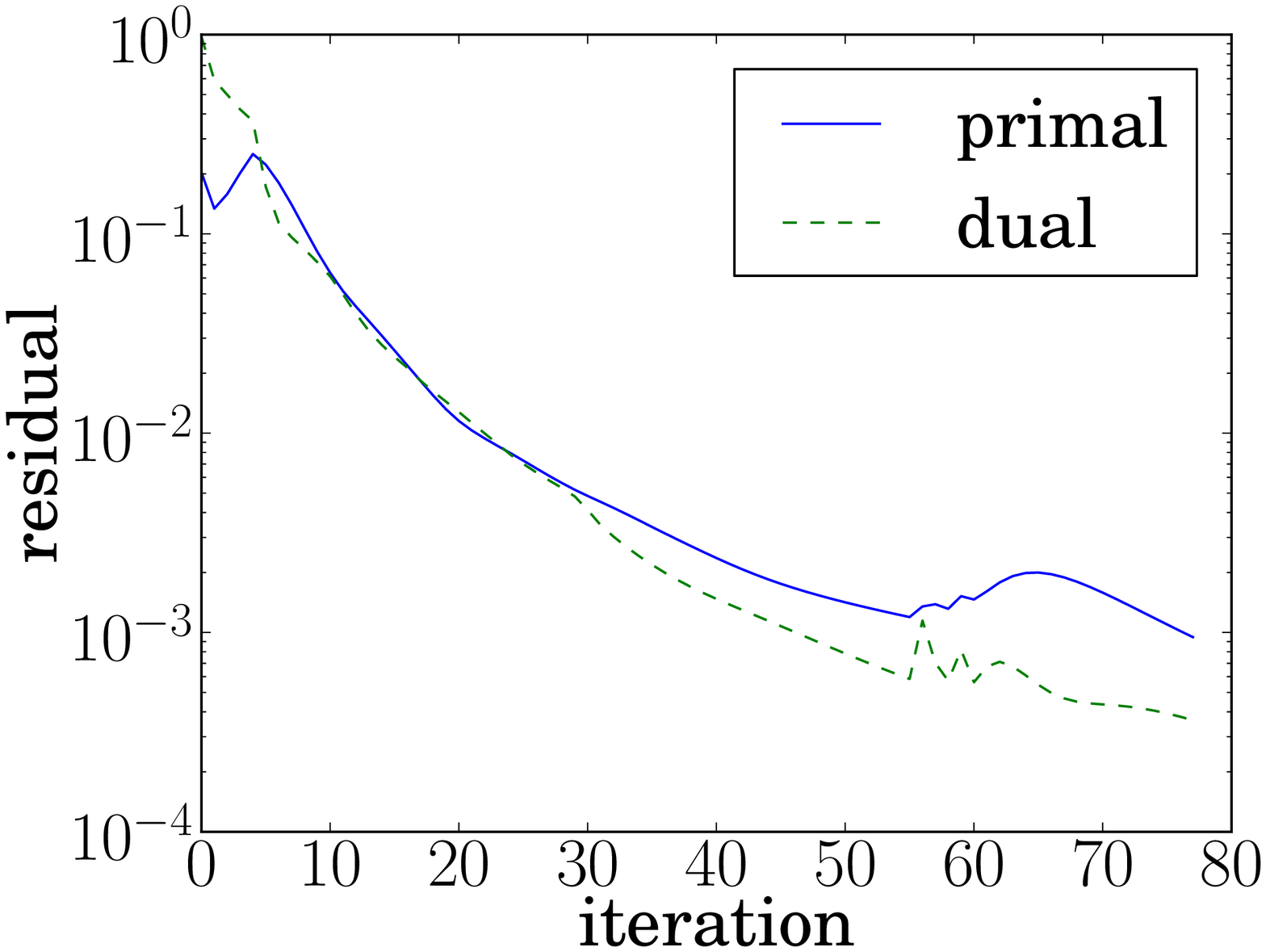}
\hspace*{\fill}
\includegraphics[width=.32\linewidth]{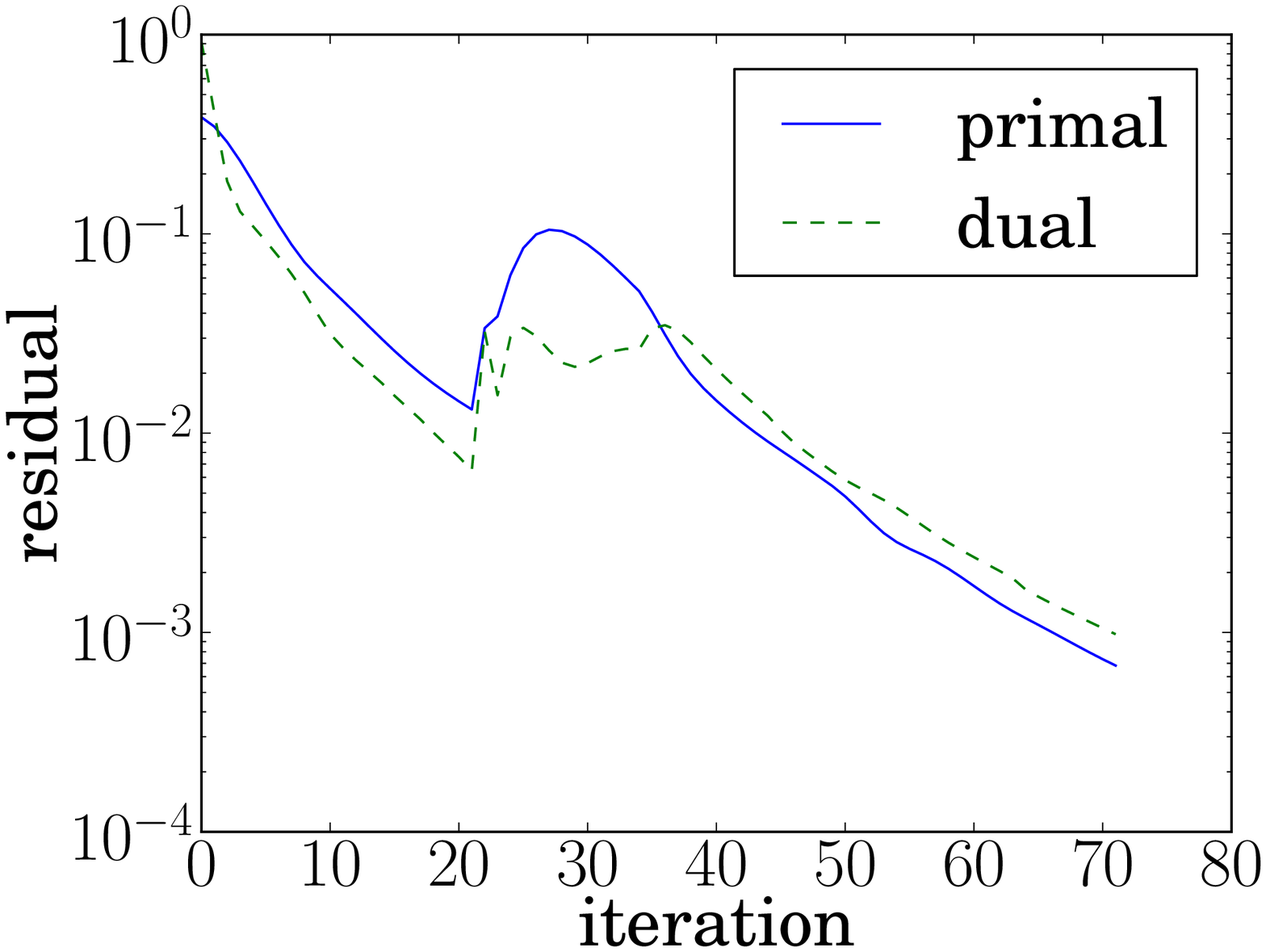}
\hspace*{\fill}
\includegraphics[width=.32\linewidth]{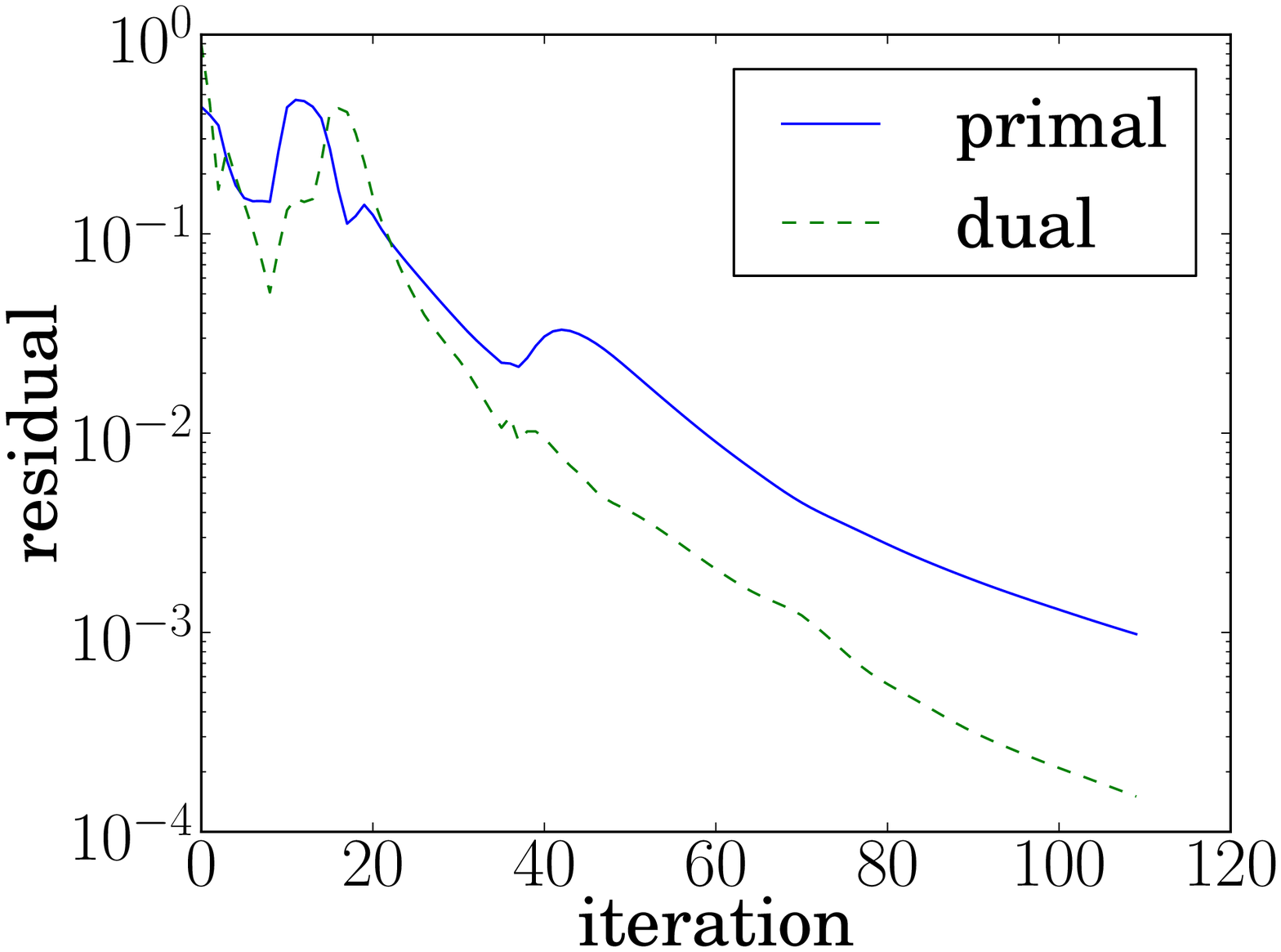}

\includegraphics[width=.32\linewidth]{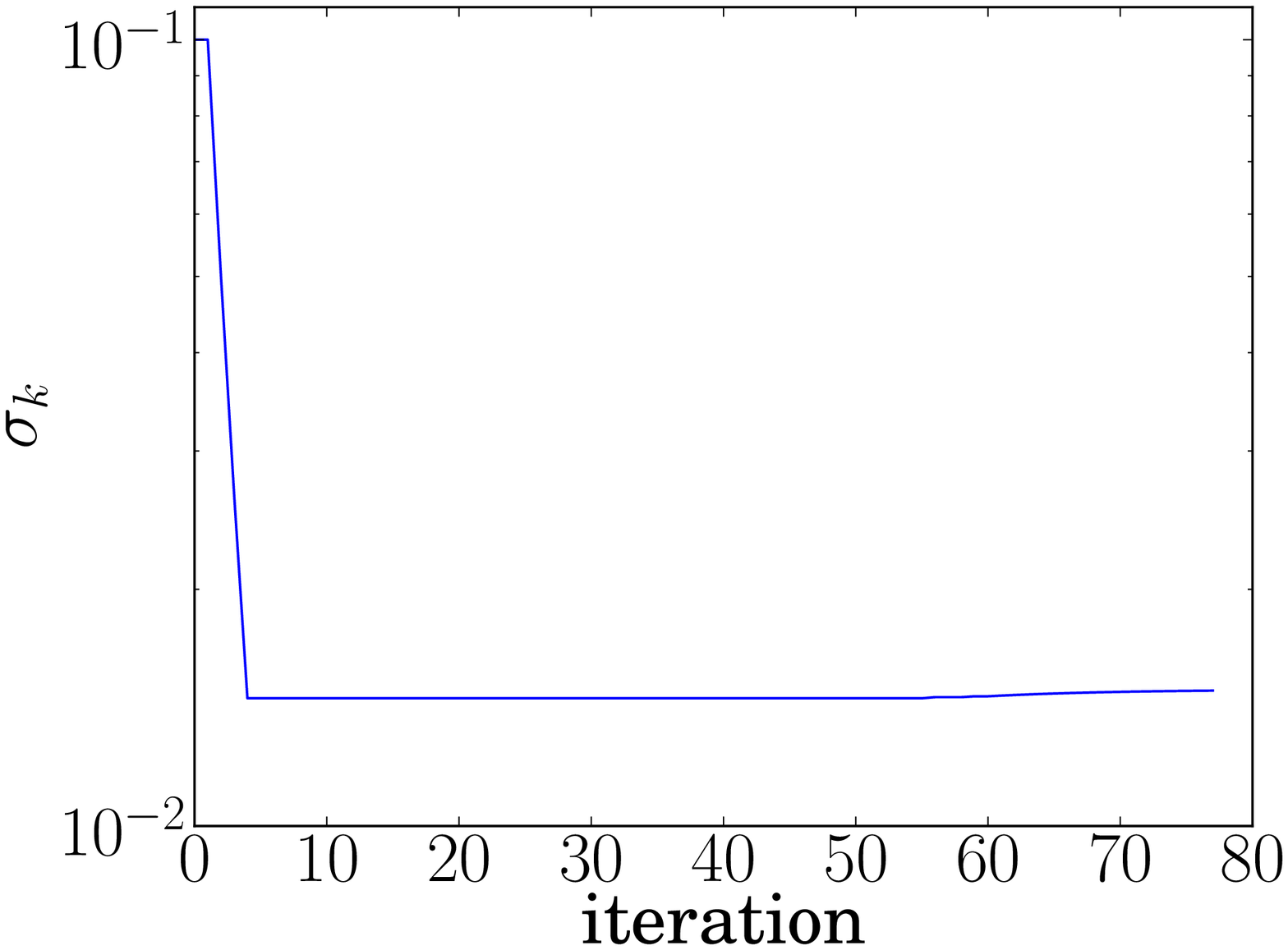}
\hspace*{\fill}
\includegraphics[width=.32\linewidth]{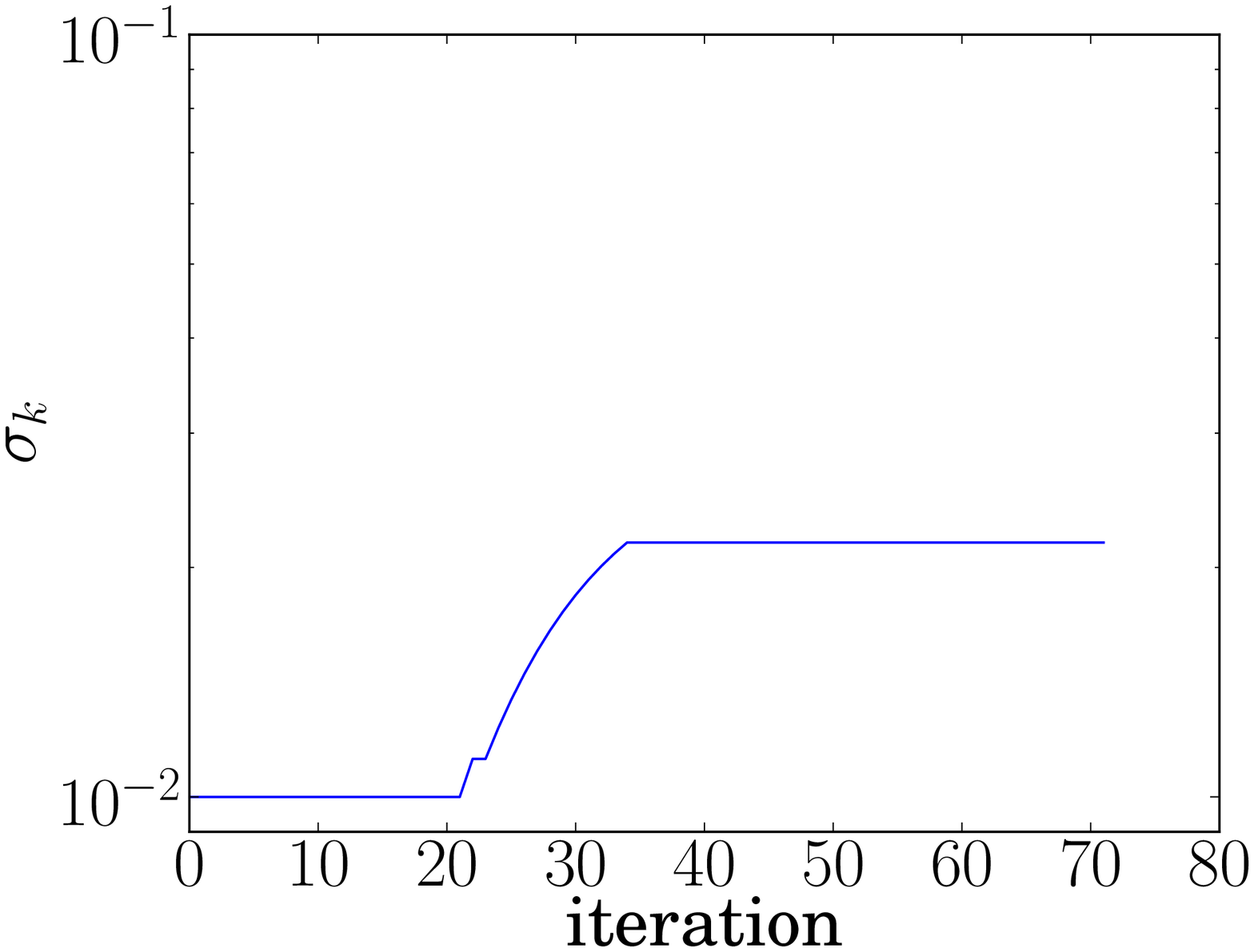}
\hspace*{\fill}
\includegraphics[width=.32\linewidth]{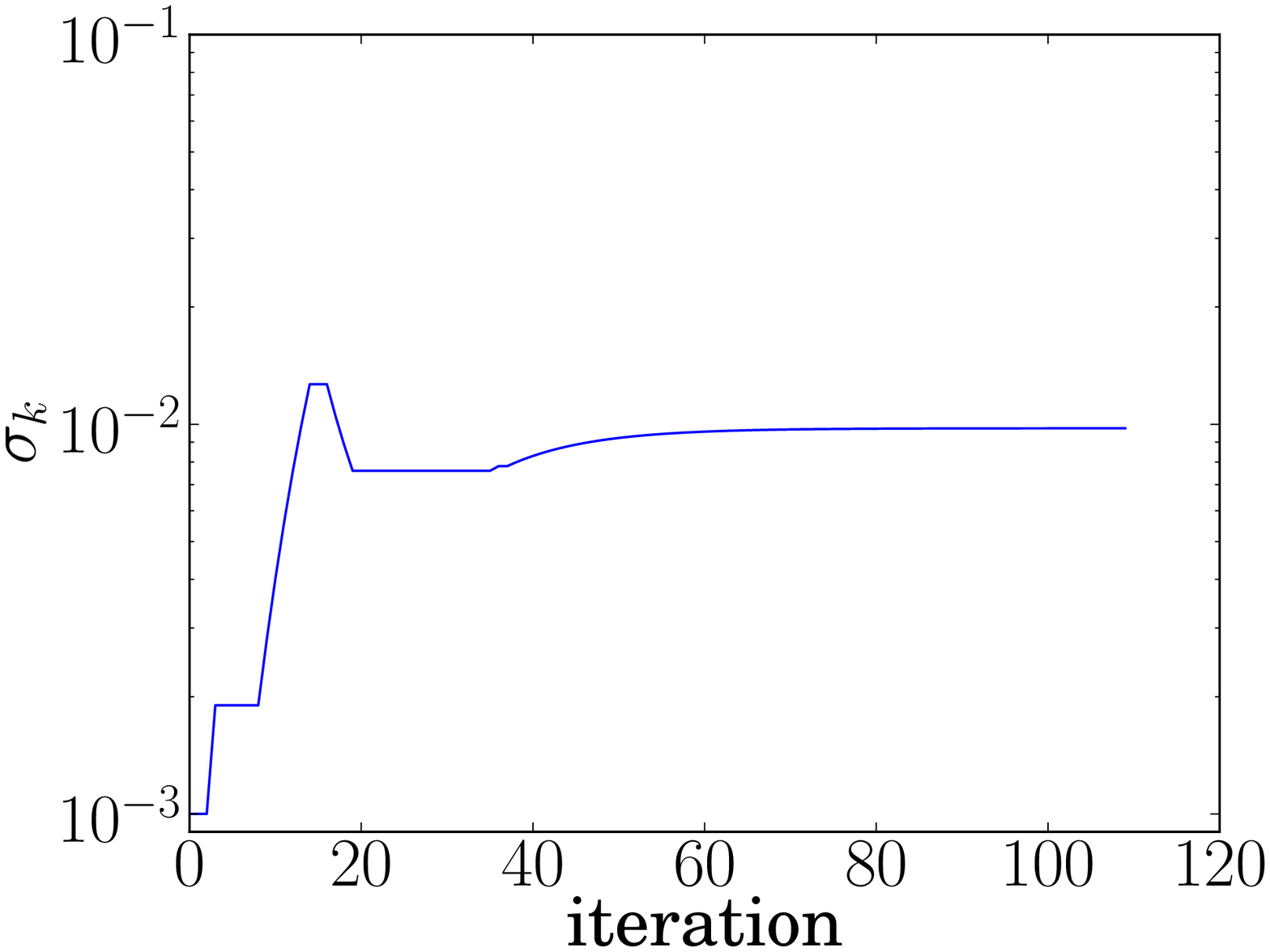}
\caption{
Primal and dual residuals versus iteration number with adaptive
selection of $\sigma_k$, starting with a value
$0.1$ (top left),  $0.01$ (middle), $0.001$ (right).
The graphs on the bottom row show the values of $\sigma_k$ during the 
three runs of the algorithm.}
\label{f-sigma-res2}
\label{f-sigma-t}
\end{center}
\end{figure}

This observation motivates  the adaptive strategy~(\ref{e-adapt-t}).
Figure~\ref{f-sigma-res2} shows the residuals if the
adaptive strategy is used, with $\mu = 2$, $\tau_k = 1+0.9^k$, and 
starting at three different values of $\sigma_k$  ($0.1$, 
$0.01$, $0.001$). 
Figure~\ref{f-sigma-t} shows the resulting values 
of $\sigma_k$ versus the iteration number $k$.

The convergence graphs indicate that a simple heuristic for adapting 
the steplength can improve the speed of convergence and make it less 
dependent on the initial steplength.
The specific convergence behavior depends on the parameter $\mu$ and 
the decay rule for $\tau_k$, but is much less sensitive to the choice
of $\sigma_0$. 
While in general the convergence with adaptive steplength is not 
faster than with a carefully tuned constant steplength,
the adaptive strategy is more robust than picking an arbitrary 
constant steplength.

\subsection{Approximate Euclidean distance matrix completion} 
\label{s-edm}
A Euclidean distance matrix (EDM) $D$ is a matrix with entries that 
can be expressed as squared pairwise distances 
$D_{ij} = \|x_i - x_j\|_2^2$ for some set of vectors $x_k$. 
In this section, we consider the problem of
fitting a Euclidean distance matrix to measurements 
$\widehat D_{ij}$ of a subset of its entries.
This and related problems arise in many applications,
including, for example, the sensor network node localization problem 
\cite{SoY:07,KKW:09,KWsnl}. 

Expanding the identity  in the definition of Euclidean distance matrix,
\[
  D_{ij} = \|x_i - x_j\|_2^2 = x_i^T x_i - 2x_i^T x_j + x_j^Tx_j,
\]
shows that a matrix $D$ is a Euclidean distance matrix
if and only if $D_{ij} = X_{ii} - 2X_{ij} + X_{jj}$
for a positive semidefinite matrix $X$   
(the Gram matrix with entries $X_{ij} = x_i^Tx_j$).
Furthermore, since $D$ only depends on the pairwise distances of the 
configuration points, we can arbitrarily place one of the points 
at the origin or, equivalently, set one row and column of $X$ to zero.
This gives an equivalent characterization: 
$D$ is a $(p+1)\times (p+1)$ Euclidean distance matrix if and only if
there exists a positive semidefinite matrix $X\in\symm^p$ such that
\[
   D_{ij} = \Tr(F_{ij}X), \quad 1 \leq i < j \leq p+1
\]
where
\[
  F_{ij} = \left\{\begin{array}{ll}
  (e_i-e_j) (e_i-e_j)^T &  1\leq i < j\leq p \\ 
  e_ie_i^T &  1 \leq i < j = p+1
  \end{array}\right.
\]
and $e_i$ denotes the $i$th  unit vector in $\reals^p$.

In the EDM approximation problem we are given a set of measurements 
$\widehat D_{ij}$ for entries $(i,j) \in W$ where 
\[
 W  \subseteq \{(i,j) \mid 1 \leq i < j \leq p+1\}.
\]
The problem of fitting a Euclidean distance matrix to the measurements
can be posed as 
\BEQ \label{e-edm-emb}
\begin{array}{ll}
\mbox{minimize} & 
 \sum\limits_{(i,j) \in W} |\Tr(F_{ij}X) - \widehat D_{ij}| \\*[3ex]
\mbox{subject to} & X \succeq 0,
\end{array}
\EEQ
with variable $X\in\symm^p$.  
(We choose the $\ell_1$-norm to measure the quality 
of the fit simply because the problem is more easily expressed
as a conic LP.)
Now let $V$ be a chordal sparsity pattern of order $p$
that includes the aggregate sparsity pattern of the matrices 
$F_{ij}$.  In other words, if $(i,j) \in W$ with $1\leq i < j \leq p$,
then $(i,j)$ is a nonzero in $V$.  
Moreover $V$ is chordal and includes all the diagonal entries in its
nonzeros. 
Such a pattern $V$ is called a \emph{chordal embedding} of $W$.
Then, without loss of generality, we can restrict
the variable $X$ in~(\ref{e-edm-emb}) to be a sparse matrix
in $\SV{p}{V}$  and we obtain the equivalent problem
\BEQ \label{e-edm-emb-c}
\begin{array}{ll}
\mbox{minimize} & \sum\limits_{(i,j) \in W} 
    |\Tr(F_{ij}X) - \widehat D_{ij}| \\*[3ex]
\mbox{subject to} & X \in\SVc{p}{V}.
\end{array}
\EEQ 
This problem is readily converted into a standard conic LP
of the form~(\ref{e-matrix-cone-LPs}), which can then be 
solved using the decomposition method of section~\ref{s-decmp-sdp}.
An interesting feature of this application is that the
correlative sparsity associated with the converted problem is block-diagonal. 

The conversion method and the block-diagonal correlative sparsity 
can also be explained directly in terms of the problem~(\ref{e-edm-emb-c}).
Suppose $V$ has $l$ cliques $\beta_k$, $k=1,\ldots, l$.
Suppose we partition the set $W$ in $l$ sets $W_k$ with the property that
if $(i,j) \in W_k$ and $1\leq i < j \leq p$, then 
$i,j\in\beta_k$, and if $(i,p+1) \in W_k$, then $i\in \beta_k$.
Then~(\ref{e-edm-emb-c}) is equivalent to
\BEQ \label{e-edm-emb-c-conv}
\begin{array}{ll}
\mbox{minimize} & \sum\limits_{k=1}^l \sum\limits_{(i,j) \in W_k} 
 |\Tr(F_{ij} \mathcal E_{\beta_k}^*(\tilde X_k)) - 
    \widehat D_{ij}| \\*[3ex]
\mbox{subject to} 
 & \mathcal E_{\eta_j} (\mathcal E_{\beta_k}^*(\tilde X_k)
   - \mathcal E_{\beta_j}^*(\tilde X_j)) = 0, \quad
   k=1,\ldots, l, \quad \beta_j\in\ch(\beta_k) \\*[1ex]
 & \tilde X_k \succeq 0, \quad k=1,\ldots,l,
\end{array}
\EEQ 
with variables $\tilde X_k \in \symm^{|\beta_k|}$, $k=1,\ldots,l$.
This problem can be solved using Spingarn's method.
At each iteration we alternate between 
projection on the subspace defined
by the consistency equations in~(\ref{e-edm-emb-c-conv})
and evaluation of a prox-operator, via the solution of 
\BEQ \label{e-edm-emb-c-prox}
\begin{array}{ll}
\mbox{minimize} & \sum\limits_{k=1}^l \sum\limits_{(i,j) \in W_k} 
 |\Tr(F_{ij} \mathcal E_{\beta_k}^*(\tilde X_k)) - 
    \widehat D_{ij}| 
 + (\sigma/2) \sum\limits_{k=1}^l \|\tilde X_k - Z_k\|_F^2 \\*[3ex]
\mbox{subject to} 
 & \tilde X_k \succeq 0, \quad k=1,\ldots,l.
\end{array}
\EEQ 
Note that this problem is separable because if $(i,j) \in W_k$, 
$F_{ij}$ is nonzero only in positions that are included
in $\beta_k\times \beta_k$.
The problems~(\ref{e-edm-emb-c-prox}) can be solved efficiently
via a straightforward modification of the interior-point method
described in section~\ref{s-prox-sdp}.

\begin{figure}
\begin{center}
\begin{psfrags}
\psfrag{t}{}
\includegraphics[width=.45\linewidth]{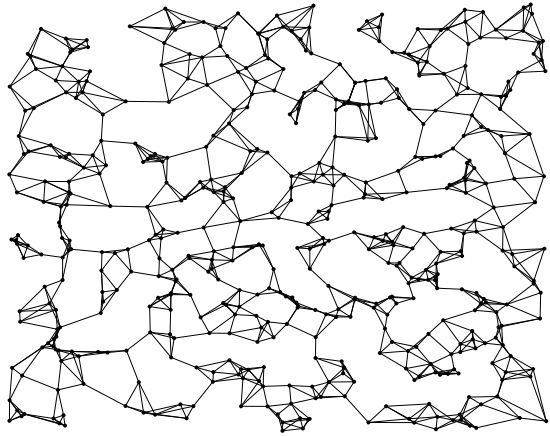}
\caption{
Nearest-neighbor network for a problem with $500$ nodes in two 
dimensions. 
Two nodes are connected if one of the two is among the 5 nearest 
neighbors of the other node.}
\label{f-edm-network}
\end{psfrags}
\end{center}
\end{figure}

We now illustrate the convergence of the decomposition method
on two randomly generated networks.
An example of a network topology is shown in 
Figure~\ref{f-edm-network} for a problem with 500 nodes. 
The network edges are assigned using the following rule:
a pair $(i,j)$ is in the sparsity pattern $W$ 
if one of the nodes is among the five nearest neighbors of the 
other node.

\begin{figure}
\begin{center}
\begin{psfrags}
\psfrag{t}{}
\includegraphics[width=.4\linewidth]{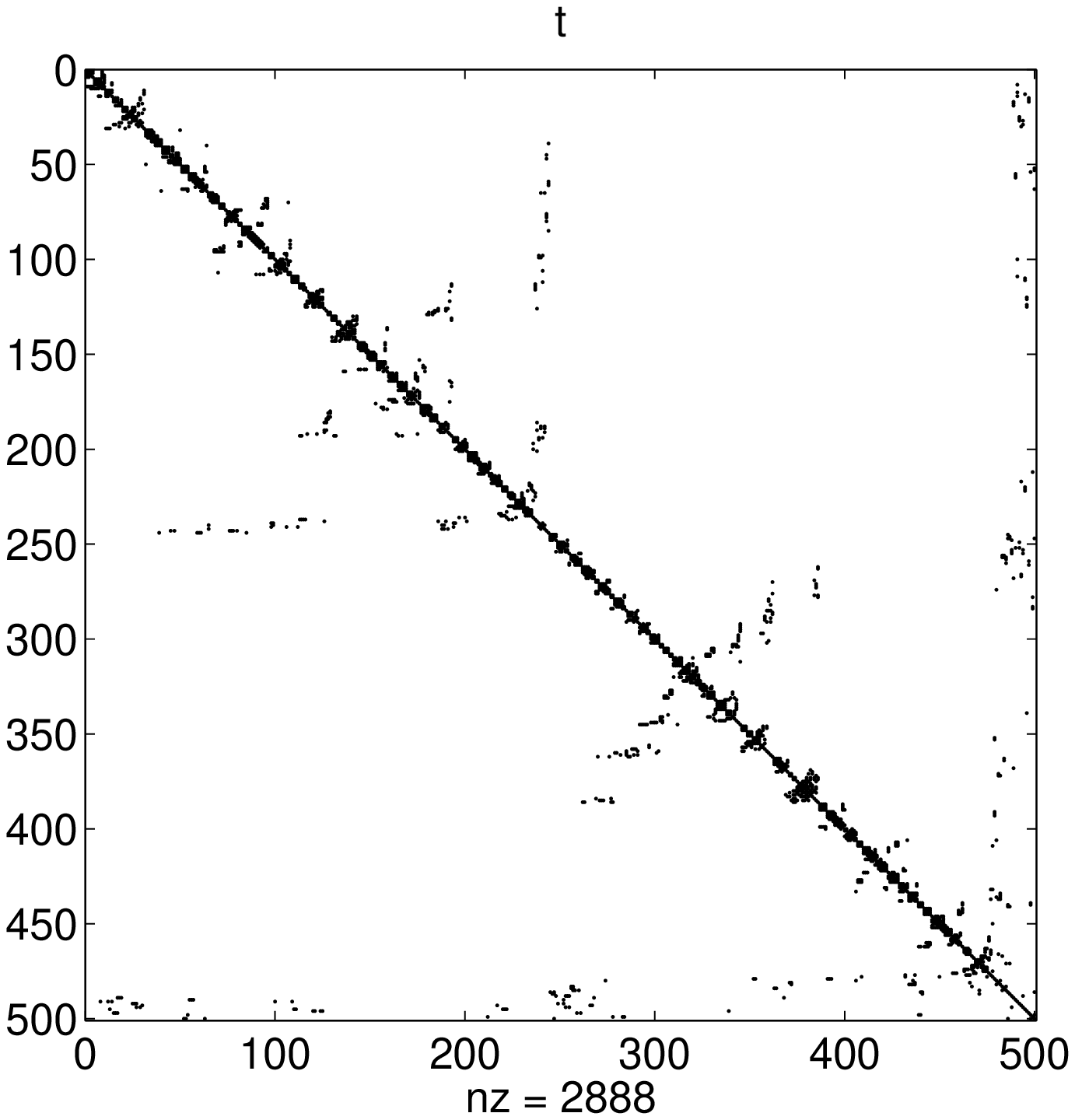}
\qquad
\includegraphics[width=.4\linewidth]{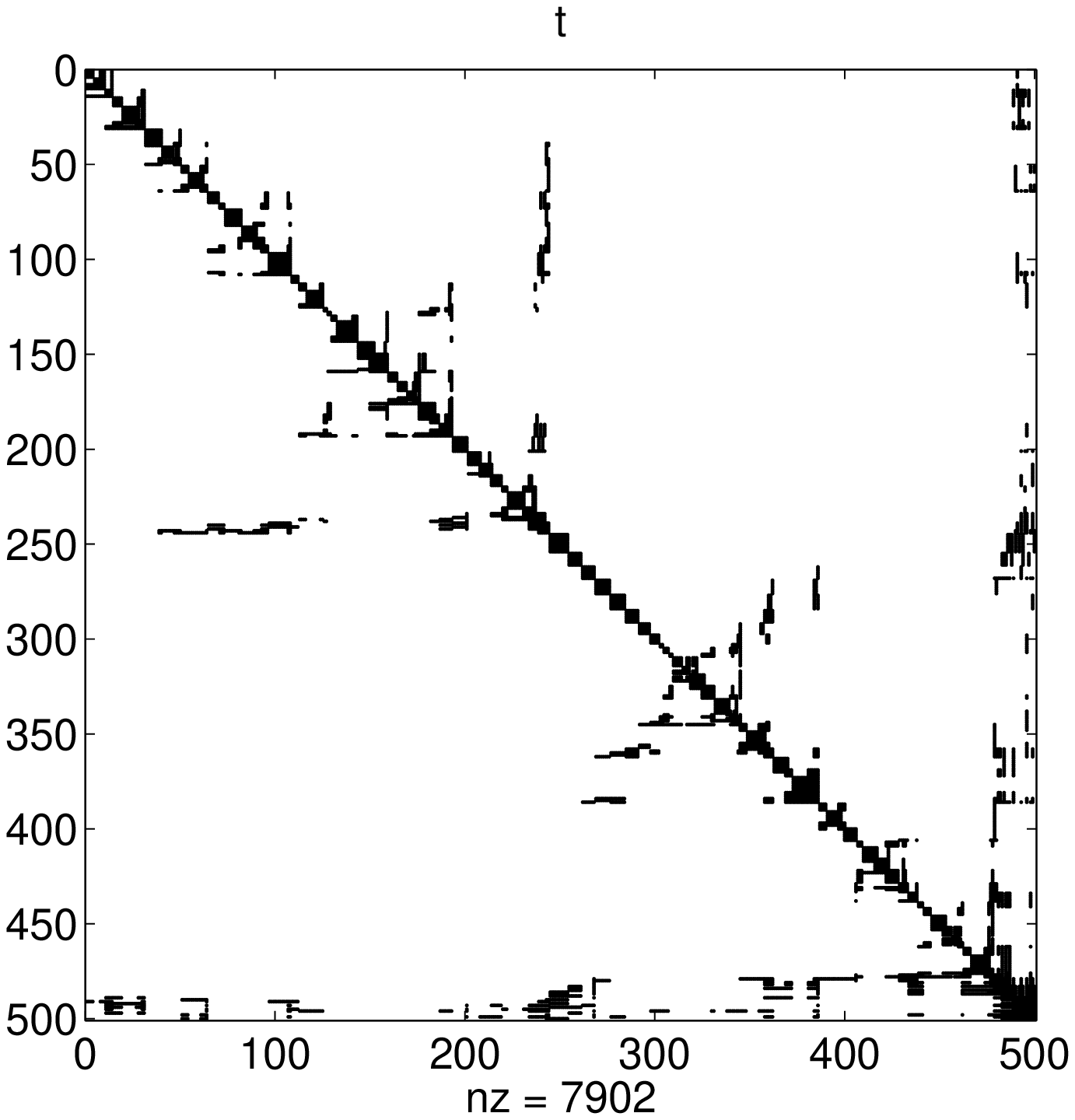}
\caption{Sparsity pattern for a network with 500 nodes after 
approximate minimum degree (AMD) reordering  and chordal embedding 
(left), and after clique merging (right). 
Before clique merging, there are 359 cliques with an average 
of 5 elements. 
After clique merging, there are 79 cliques with an average of 5 elements.
}
\label{f-edm-preprocessing}
\end{psfrags}
\end{center}
\end{figure}

To compute a chordal embedding $V$, we use an approximate minimum degree 
(AMD) reordering, which gives a permutation of the sparsity pattern 
that reduces fill-in (Figure \ref{f-edm-preprocessing}, left). 
Often, the resulting embedding contains many small cliques 
and for our purposes it is more efficient
to merge some neighboring cliques, 
using algorithms similar to those in \cite{CR:89,ReSc:09,HoSc:10}.
Specifically, 
traversing the tree in a topological order, 
we greedily merge clique $k$ with its parent if 
\[  
(|\beta_{\prnt(k)}| - |\eta_k|)(|\beta_k| - |\eta_k|) 
\leq t_{\mathrm{fill}} \quad \text{or} \quad 
\max(|\beta_k|-|\eta_{k}|, |\beta_{\prnt(k)}|-|\eta_{\prnt(k)}|) \leq t_{\mathrm{size}} 
 \] 
where $t_{\mathrm{fill}}$ is a threshold based on the amount of fill 
that results from merging clique $k$ with its parent, 
and $t_{\mathrm{size}}$ is a threshold based on the cardinality of the 
sets $\beta_{\prnt(k)}\setminus \eta_{\prnt(k)}$ and 
$\beta_{k}\setminus \eta_{k}$. 
In Figure~\ref{f-edm-preprocessing} (right) we show the result of this 
clique-merging technique 
using the values $t_\mathrm{fill} = t_\mathrm{size} = 5$.
This reduced the 359 original cliques with an average of 5 nodes 
each to 79 cliques with an average of 10 nodes.

\begin{figure}
\begin{center}
\begin{psfrags}
\psfrag{t}{}
\psfrag{x}[t][t]{iteration}
\psfrag{z}[b][b]{residual}
\includegraphics[width=.48\linewidth]{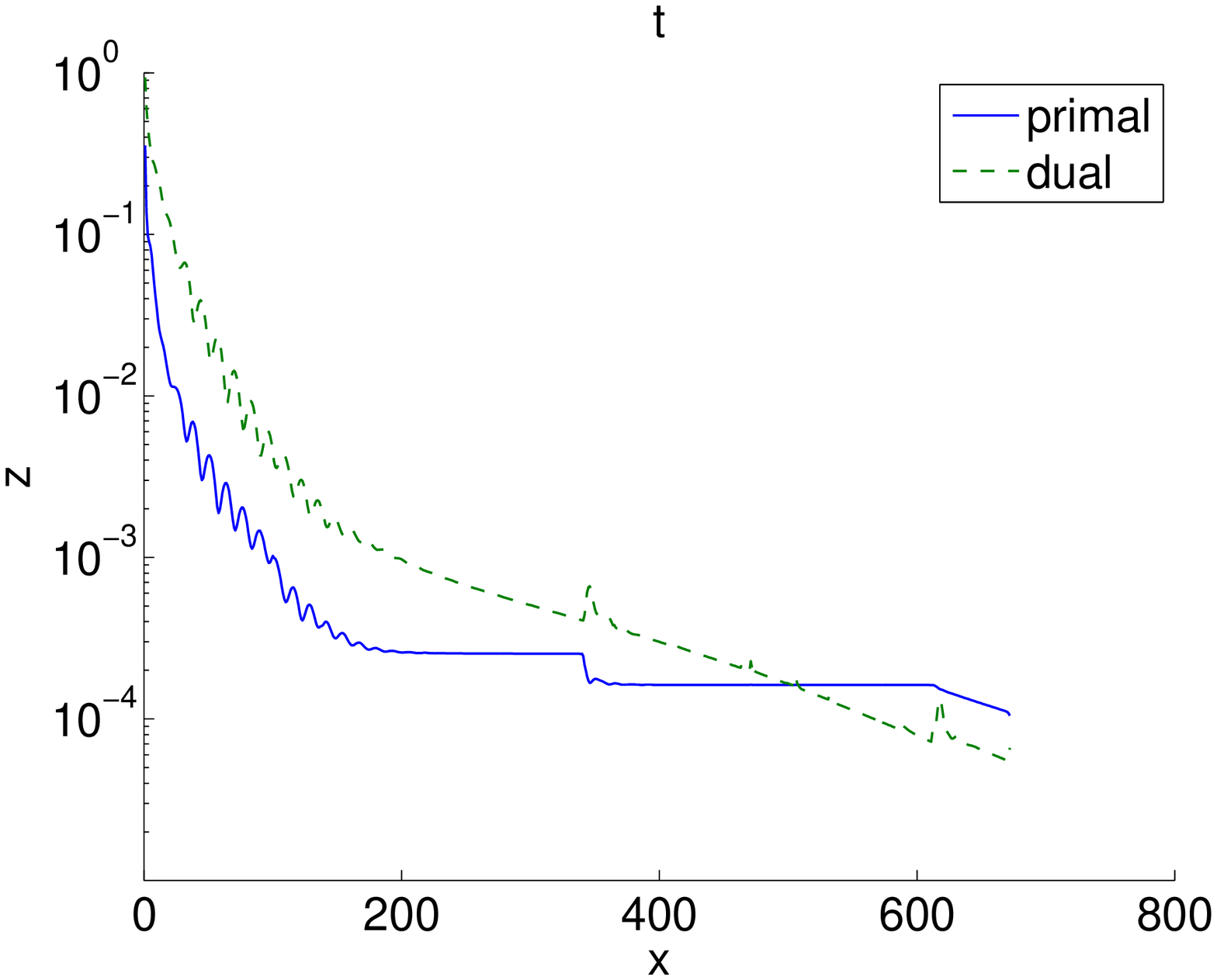}
\quad
\includegraphics[width=.48\linewidth]{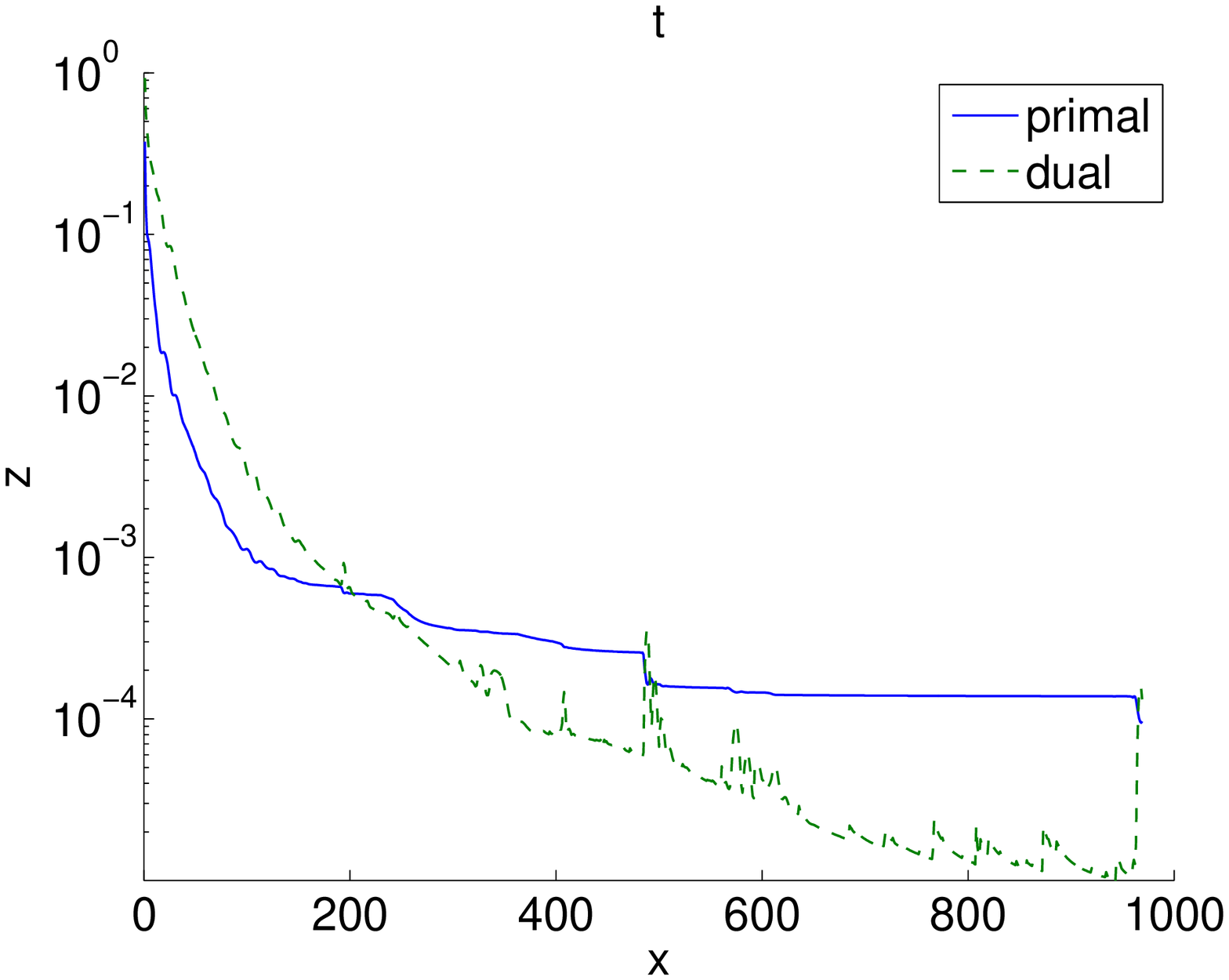}
\end{psfrags}
\caption{Relative primal and dual residuals versus iteration number 
for networks with 500 (upper left) and 2000 (right) nodes. 
For $n = 500$, there are 82 cliques, and for $n = 2000$, there 
are 310 cliques. A constant steplength parameter $\sigma_k = 5.0$ is used.
}
\label{f-edm-evolution}
\end{center}
\end{figure}

A typical convergence plot of the resulting problem is given in 
Figure~\ref{f-edm-evolution} for a network with 500 nodes (left) and 
2000 nodes (right). A constant value $\sigma_k = 5.0$
is used for the steplength parameter. 
The greedy clique merging strategy described above was used,
with the same threshold values.

\subsection{Block-arrow semidefinite programs} \label{s-blck-arrow}
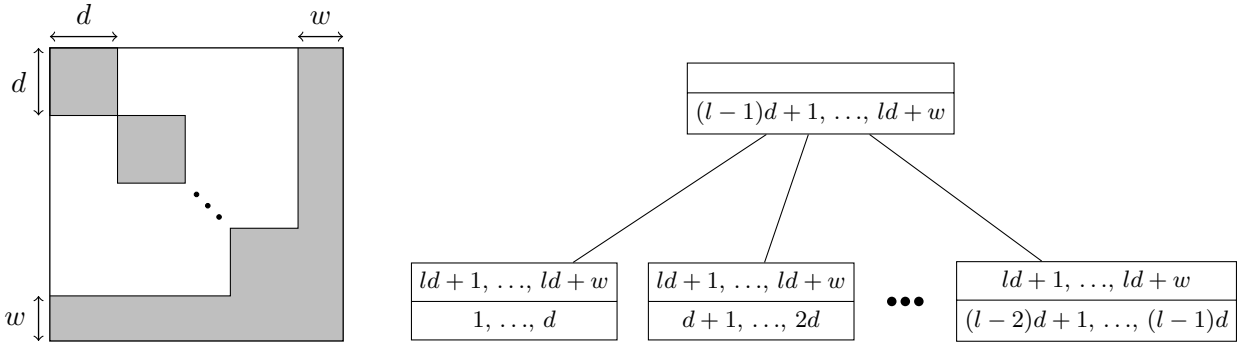
\begin{figure}
\begin{tikzpicture}[scale = 0.3]

\fill[line width = 2pt, color=gray!50] 
   (0,0) -- (0,-3) -- (3,-3) -- (3,0) -- (0, 0);
\draw (0,0) rectangle(3,-3);

\fill[line width = 2pt, color=gray!50] 
   (3,-3) -- (3,-6) -- (6,-6) -- (6,-3) -- (3, -3);
\draw (3,-3) rectangle(6,-6);

\fill[line width = 2pt, color=gray!50] 
   (8,-8) -- (8,-11) -- (11,-11) -- (11,-8) -- (8, -8);

\fill[line width = 2pt, color=gray!50] 
   (0,-11) -- (8,-11) -- (8,-8) -- (11,-8) -- (11, 0) -- (13,0) --
   (13,-13) -- (0, -13) -- (0,-11);
\draw (0,-11) -- (8,-11) -- (8,-8) -- (11,-8) -- (11,0);

\filldraw (6.5,-6.5) circle (.1);
\filldraw (7,-7) circle (.1);
\filldraw (7.5,-7.5) circle (.1);
\draw[line width = .5pt] (0,0) rectangle(13,-13);

\draw[<->, line width = .5pt] (-.5,0) -- (-.5,-3);
\node[anchor = east] at (-0.6,-1.5) {$d$};

\draw[<->, line width = .5pt] (-.5,-11) -- (-.5,-13);
\node[anchor = east] at (-0.6,-12) {$w$};

\draw[<->, line width = .5pt] (0,.5) -- (3,.5);
\node[anchor = south] at (1.5, 0.6) {$d$};

\draw[<->, line width = .5pt] (11, .5) -- (13, .5);
\node[anchor = south] at (12,0.6) {$w$};

\end{tikzpicture}
\hspace*{\fill}
\begin{tikzpicture}[scale=0.9, 
clique/.style = {rectangle split, rectangle split parts = 2, draw}]
\footnotesize
\node[clique, anchor=east](1) at (0,0){$ld+1$, $\ldots$, $ld+w$ 
  \nodepart{second} $1$, $\ldots$, $d$ };
\node[clique, anchor=east](2) at (3.5,0){$ld+1$, $\ldots$, $ld+w$ 
  \nodepart{second} $d+1$, $\ldots$, $2d$ };
\node[clique, anchor=west](3) at (5,0){$ld+1$, $\ldots$, $ld+w$ 
  \nodepart{second} $(l-2)d+1$, $\ldots$, $(l-1)d$ };
\foreach \x/\y in {4.05/0, 4.25/0, 4.45/0 } \fill (\x,\y) circle (2pt);
\node[clique, anchor=center](4) at (3,3){
  \nodepart{second} $(l-1)d+1$, $\ldots$, $ld+w$ };
\draw (1) -- (4);
\draw (2) -- (4);
\draw (3) -- (4);
\end{tikzpicture}
\caption{Block arrow pattern with $l$ cliques and corresponding clique 
tree.  The order of the matrix is $ld+w$.
The first $l$ diagonal blocks in the matrix have size $d$, the last
block column and block row have width $w$.
The cliques therefore have size $d+w$.
Each clique in the clique tree is partitioned in two sets:
the top row shows $\eta_k = \beta_k \cap \prnt(\beta_k)$;
the bottom row shows $\beta_k \setminus \eta_k$.  }
\label{f-blck-arrow}
\end{figure}

In the last experiment we compare the efficiency of the splitting
method with general-purpose SDP solvers.
We consider a family of randomly generated SDPs with 
a block-arrow aggregate sparsity pattern $V$ and a block-diagonal 
correlative sparsity pattern.  
The sparsity pattern $V$ is defined in Figure~\ref{f-blck-arrow}.
It consists of $l$ diagonal blocks of size $d\times d$, plus
$w$ dense final rows and columns.
We take the clique
\[
 \beta_l = \{ (l-1)d+1, \, \ldots, \, ld, \, ld+1, \, \ldots, \, ld+w\} 
\]
(with $\eta_l = \{\}$)  as root of the clique tree.
The other $l-1$ cliques $\beta_k$
and the intersections $\eta_k = \beta_k \cap \prnt(\beta_k)$ with
their parent cliques are 
\[
 \beta_k = \{(k-1)d+1, \, \ldots, \, kd\} \cup  \eta_k, \qquad
 \eta_k = \{ld+1, ld+2, \ldots, ld+w\}, \qquad k=1,\ldots,l-1. 
\]
We generate matrix cone LPs~(\ref{e-matrix-cone-LPs})
with $m=ls$ primal equality constraints, partitioned in $l$ sets 
\[
 \nu_k = \{(k-1)s+1, \, (k-1)s+2, \, \ldots, \, ks\},\qquad k=1,\ldots,l,
\]
of equal size $|\nu_k| =  s$.   
If $i\in \nu_k$, then the coefficient matrix
$F_i$ contains a dense $\beta_k\times \beta_k$ block, and is 
otherwise zero.  We will use the notation
\[
 (F_i)_{\beta_k\beta_k} = 
\left[\begin{array}{cc} A_i & B_i \\ B_i^T & D_i \end{array}\right],
\]
for the nonzero block of $F_i$ if $i\in \nu_k$.
The primal and dual SDPs can therefore be expressed as
\BEQ \label{e-block-arrow-unconverted}
\begin{array}[t]{ll}
\mbox{minimize} & \Tr(CX) \\ 
\mbox{subject to} & \mathcal A(X) = b \\ & X\succeq 0
\end{array} \qquad\qquad
\begin{array}[t]{ll}
\mbox{maximize} & b^Ty \\
\mbox{subject to} & \mathcal A^*(y) + S = C \\ & S\succeq 0
\end{array} 
\EEQ
with a linear mapping $\mathcal A: \symm^{(ld+w)\times (ld+w)} \rightarrow
\reals^{ls}$ defined as
\[
\mathcal A(X)_i = 
 \Tr \left(
 \left[\begin{array}{cc} A_i & B_i \\ B_i^T & D_i \end{array}\right]
 \left[\begin{array}{cc}
   X_{kk} & X_{k,l+1} \\ X_{l+1,k} & X_{l+1,l+1} 
   \end{array}\right] \right), \quad i \in \nu_k, 
   \quad k=1,\ldots, l,
\]
where $X_{ij}$ denotes the $i,j$ block of $X$.
(These blocks have dimensions $X_{ii} \in \symm^d$ for $i=1,\ldots,l$,
$X_{l+1,l+1} \in \symm^w$, $X_{l+1,i} \in \reals^{w\times d}$ for
$i=1,\ldots,l$.)
The adjoint $\mathcal A^*: \reals^{ls} \rightarrow 
\symm^{(ld+w)\times (ld+w)}$ is
\[
\mathcal A^*(y) = 
\left[\begin{array}{ccccc}
\sum\limits_{i\in\nu_1} y_i A_i & 0 & \cdots & 0 &
\sum\limits_{i\in\nu_1} y_i B_i \\
0 & \sum\limits_{i\in\nu_2} y_i A_i & \cdots & 0 & 
\sum\limits_{i\in\nu_2} y_i B_i \\
\vdots & \vdots & \ddots & \vdots & \vdots \\ 
0 & 0 & \cdots & \sum\limits_{i\in\nu_l} y_iA_i & 
\sum\limits_{i\in\nu_l} y_i B_i \\
\sum\limits_{i\in\nu_1}y_iB_i^T &
\sum\limits_{i\in\nu_2}y_iB_i^T & \cdots & 
\sum\limits_{i\in\nu_l}y_iB_i^T & 
\sum\limits_{i=1}^m y_i D_i
\end{array}\right].
\]

In the reformulated problem, the variable $X$ is replaced with 
$l$ matrices $\tilde X_k = X_{\beta_k\beta_k}$, \ie, defined as
\[
\tilde X_k = 
\left[\begin{array}{cc}
 (\tilde X_k)_{11} & (\tilde X_k)_{12} \\ 
 (\tilde X_k)_{21} & (\tilde X_k)_{22} 
 \end{array}\right]
 = \left[\begin{array}{cc}
   X_{kk} & X_{k,l+1} \\ X_{k,l+1}^T & X_{l+1,l+1}
 \end{array}\right], \quad k=1,\ldots,l,
\]
and the primal SDP is converted to
\BEQ \label{e-block-arrow-converted}
 \begin{array}{ll}
 \mbox{minimize} & \sum\limits_{k=1}^l \Tr(\tilde C_k \tilde X_k) \\
 \mbox{subject to} & 
 \Tr\left(\left[\begin{array}{cc}
   A_i &B_i \\ B_i^T & D_i \end{array}\right]
  \left[\begin{array}{cc}  
   (\tilde X_k)_{11} & (\tilde X_k)_{12} \\
   (\tilde X_k)_{21} & (\tilde X_k)_{22} \end{array}\right]\right)
 = b_i, \quad i\in\nu_k, \quad k=1,\ldots,l \\*[2ex]
 & (\tilde X_k)_{22} = (\tilde X_l)_{22}, \quad k=1,\ldots,l-1  \\
 & \tilde X_k \succeq 0, \quad k=1,\ldots,l
 \end{array}
\EEQ
where
\[
\tilde C_k = \left[\begin{array}{cc}
  C_{kk} & C_{k,l+1} \\ C_{k,l+1}^T & 0 \end{array}\right],
 \quad k=1,\ldots,l-1, \qquad
\tilde C_l = \left[\begin{array}{cc}
  C_{ll} & C_{l,l+1} \\ C_{l,l+1}^T & C_{l+1,l+1} \end{array}\right].
\]
With this choice of parameters, the correlative sparsity pattern 
of the converted SDP~(\ref{e-block-arrow-converted}) 
is block-diagonal, \ie, except for the consistency constraints
$(\tilde X_k)_{22} = (\tilde X_l)_{22}$ the problem is separable
with independent variables $\tilde X_k \in \symm^{p+w}$.
This allows us to compute the prox-operator by solving $l$ 
independent conic QPs.

\paragraph{Problem generation}
The problem data are randomly generated as follows.
First, the entries of $A_k$, $B_k$, $D_k$ are drawn independently
from a normal distribution $\mathcal N(0,1)$. 
A strictly primal feasible $X$ is constructed 
as $X = W + \alpha I$ where
$W\in\SV{p}{V}$ is randomly generated 
with i.i.d.\ entries from $\mathcal N(0,1)$ and $\alpha$ is chosen so that
$X_{\beta_k\beta_k} = W_{\beta_k\beta_k} + \alpha I \succ 0$ for
$k =1,\ldots,l$.
The right-hand side $b$ in the primal constraint is
computed as $b_i = \Tr(F_iX)$, $i=1,\ldots,m$.

Next, strictly dual feasible $y\in\reals^m$, $S \in \SV{p}{V}$ are 
constructed.  The vector $y$ has i.i.d.\ entries from $\mathcal
N(0,1)$ and $S$ is constructed as
$S = \sum_{k=1}^l \mathcal E^*_{\beta_k}(\tilde S_k)$, 
with $\tilde S_k = W_k + \alpha I$, $W_k\in\symm^{|\beta_k|}$ randomly 
generated with i.i.d.\ $\mathcal N(0,1)$ entries, and
$\alpha$ chosen so that $\tilde S_k \succ 0$.
Finally, the matrix $C$ is constructed as $C = S + \sum_i y_i F_i$.

\iffalse
\begin{itemize}
\item The elements of $B_i$ and the lower triangular elements of 
$A_i, D_i$ are drawn from a normal Gaussian distribution for 
$i\in \nu_k, \,k = 1,\hdots,l$.

\item A primal feasible $X\in \mathcal C$ is generated as 
\begin{eqnarray*}
X = \mathcal P_\mathcal V(W) + q \cdot 0.1 I
\end{eqnarray*}
where each lower triangular element of $W$ is drawn from a normal 
Gaussian distribution, and $q$ is the smallest integer such that 
$\mathcal E_{\beta_k}(X)\succeq 0$, $k = 1,\ldots,l$. for each 
clique $k$. 

\item A dual feasible $S\in \mathcal C^*$ is generated as
\begin{eqnarray*}
S = \sum_{k=1}^l\mathcal E_{\beta_k}^*(\xi_k)
\end{eqnarray*}
where
\begin{eqnarray*}
\xi_k = w_k + q \cdot 0.1 I
\end{eqnarray*}
where each lower-triangular element of $w_k$ is drawn from a normal 
Gaussian distribution, and $q$ is the smallest integer such that 
$\xi_k \succeq 0$. For large $p$ and $r$, this method is preferred to 
generating a dense $p(r+1)\times p(r+1)$ matrix $X$, which can quickly 
exceed memory limitations.

\item A dual feasible $y\in \reals^m$ is generated by drawing each $y_i$ from a normal Gaussian distribution.

\item A primal feasible  $b= \mathcal A(X)$ and a dual feasible $C = \mathcal A^*(y)+S$ is then directly computed.
\end{itemize}
\fi

\begin{figure}
\small
\begin{psfrags}
\psfrag{z}[b][b]{time (sec)}
\psfrag{t}{}
\psfrag{x}[t][t]{$w$ (arrow width)}
\includegraphics[width=.49\linewidth]{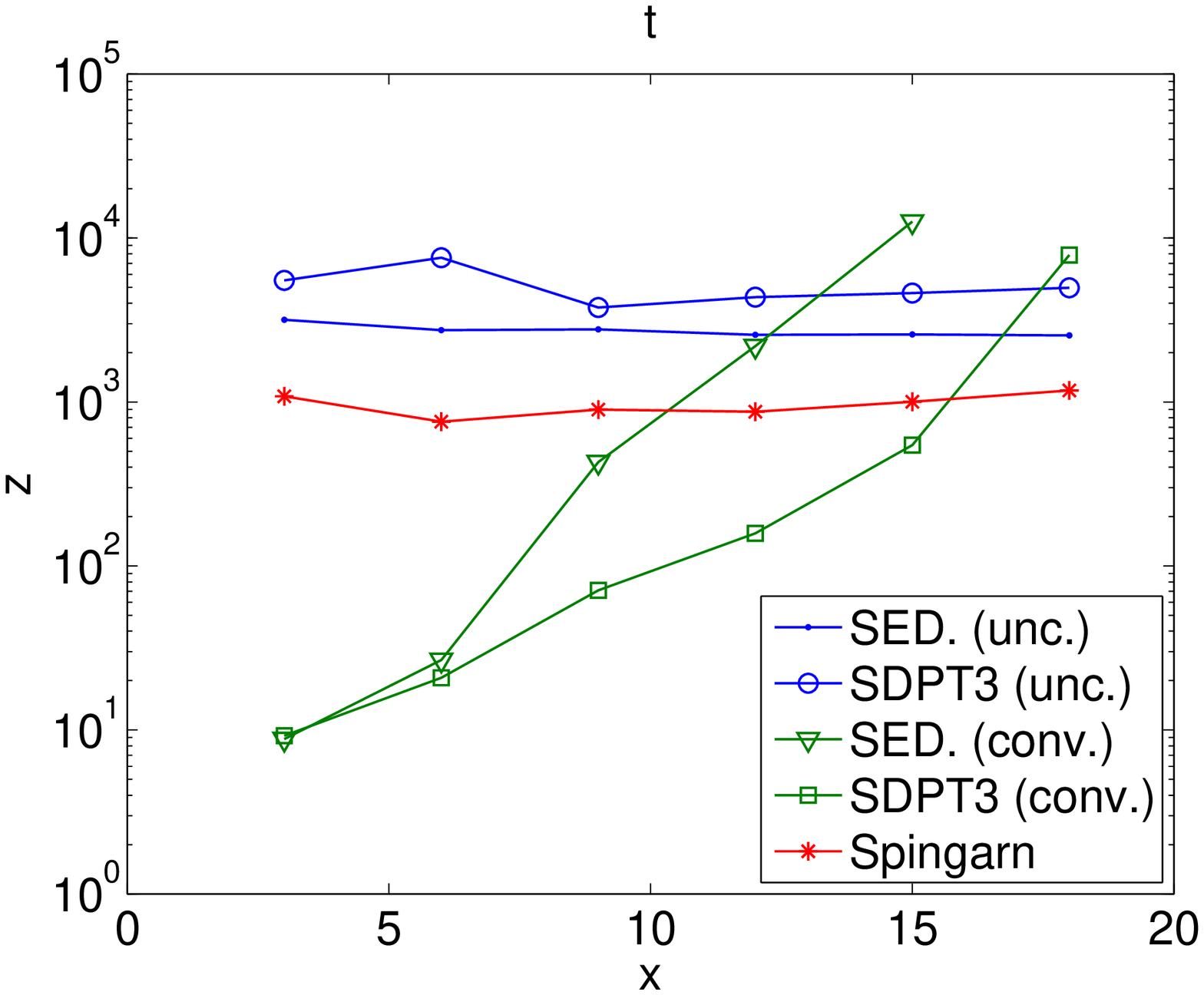}
\end{psfrags}
\hspace*{\fill}
\begin{psfrags}
\psfrag{z}[b][b]{time (sec)}
\psfrag{t}{}
\psfrag{x}[t][t]{$l$ (number of cliques)}
\includegraphics[width=.49\linewidth]{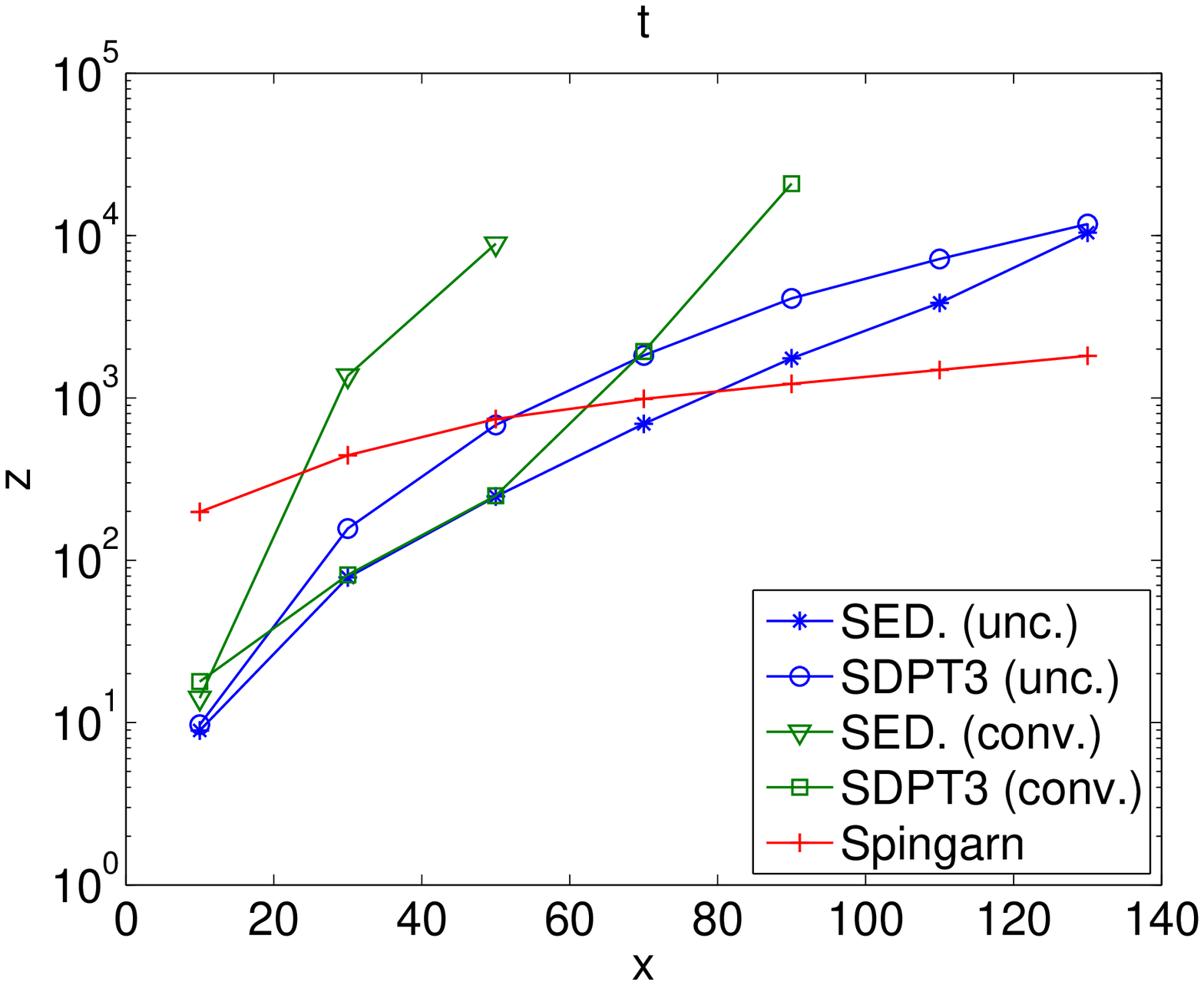}
\end{psfrags}
\caption{Solution time for randomly generated SDPs with block-arrow
sparsity patterns.  Times are reported for SEDUMI (SED.) and SDPT3
applied to the original (`unc.') and converted (`conv') SDPs,
and the Spingarn method applied to
the converted SDP.
The figure on the left shows the times as function of arrow width $w$,
for fixed dimensions $l=100$, $d=20$, $s =10$. 
The figure on the right shows the times versus number of cliques $l$,
for fixed dimensions $w=20$, $d=20$, $s =10$.
}
\label{fig:simple_blockarrow}
\end{figure}

\paragraph{Comparison with general-purpose SDP solvers}
In Figure~\ref{fig:simple_blockarrow} we compare the solution time
of Spingarn's algorithm with the general-purpose
interior-point solvers SEDUMI and SDPT3, applied to 
the unconverted and converted SDPs~(\ref{e-sdps}) 
and~(\ref{e-sdp-converted-primal}).
In the decomposition method we use a constant steplength parameter
$\sigma_k=5$ and relaxation parameter $\rho_k = 1.75$.
The stopping criterion is~(\ref{e-stopping}) with 
$\epsilon_\mathrm p = \epsilon_\mathrm d   = 10^{-4}$.
For each data point we report the average CPU time over 5 instances.  

To interpret the results, it is useful to consider the linear algebra
complexity per iteration of each method.
The unconverted SDP~(\ref{e-block-arrow-unconverted})
has a single matrix variable $X$ of order $p=ld+w$.  
The cost per iteration of an interior-point method is 
dominated by the cost of forming and solving the Schur complement
equation, which is dense and of size $m=sl$.
For the problem sizes used in the figures ($w$ small compared to 
$ld$) the cost of solving the Schur complement dominates the 
overall complexity.
This explains the nearly constant solution time in the first
figure (fixed $l$, $s$, $p$, varying $w$) and the increase with $l$
shown in the second figure.

The converted SDP~(\ref{e-block-arrow-converted})  
has $l$ variables $\tilde X_k$ of order $d+w$.
The Schur complement equation in an interior-point method has 
the general structure~(\ref{e-converted-schur}) with a leading
block-diagonal matrix ($l$ blocks of size $s\times s$)
augmented with a dense block row and block column of width
proportional to $lw^2$.  For small $w$, exploiting the block-diagonal 
structure in the Schur complement equation, allows one 
to solve the Schur complement equation very quickly and reduces
the cost per iteration to a fraction of a cost of solving the 
unconverted problem, despite the increased size of the problem.
However the advantage disappears with increasing $w$
(Figure \ref{fig:simple_blockarrow} left).    

The main step in each iteration of the Spingarn method 
applied to the converted problem is the evaluation of the
prox-operators via an interior-point method.  The Schur complement 
equations that arise in this computation are block-diagonal
($l$ blocks of order $s$) and therefore the cost of solving them is 
independent of $w$ and linear in $l$.
As an additional advantage, since the correlative sparsity 
pattern is block-diagonal, the proximal operator can be evaluated
by solving $l$ independent conic QPs that can be solved in parallel.   
This was not implemented in the experiment, but could reduce the solution 
time by a factor of roughly $l$.

\begin{figure}
\small
\centering
\psfrag{x}[t][t]{$l$ (number of cliques)}
\psfrag{z}[b][b]{number of iterations}
\psfrag{t}{}
\includegraphics[width=.49\linewidth]{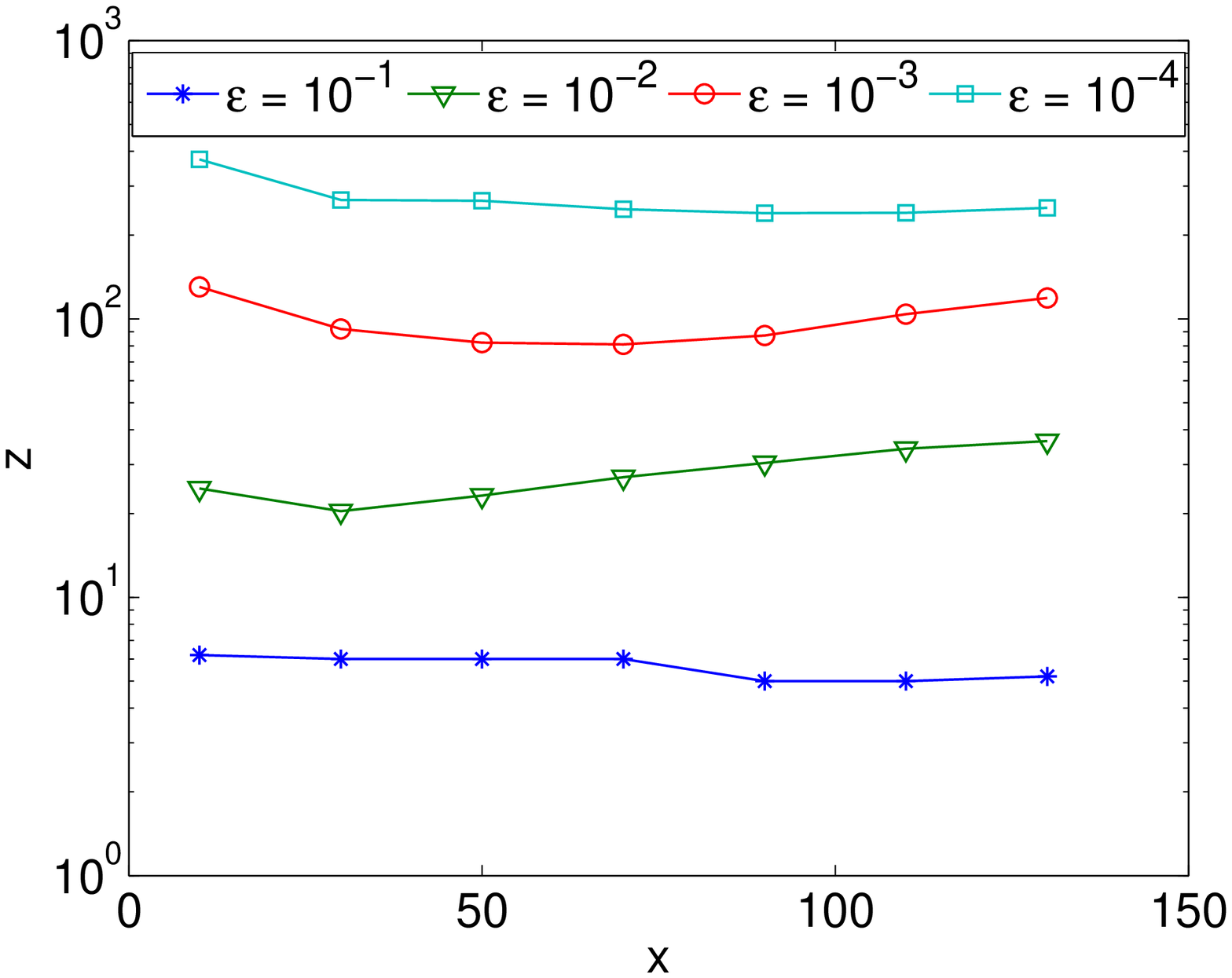}
\hspace*{\fill}
\includegraphics[width=.49\linewidth]{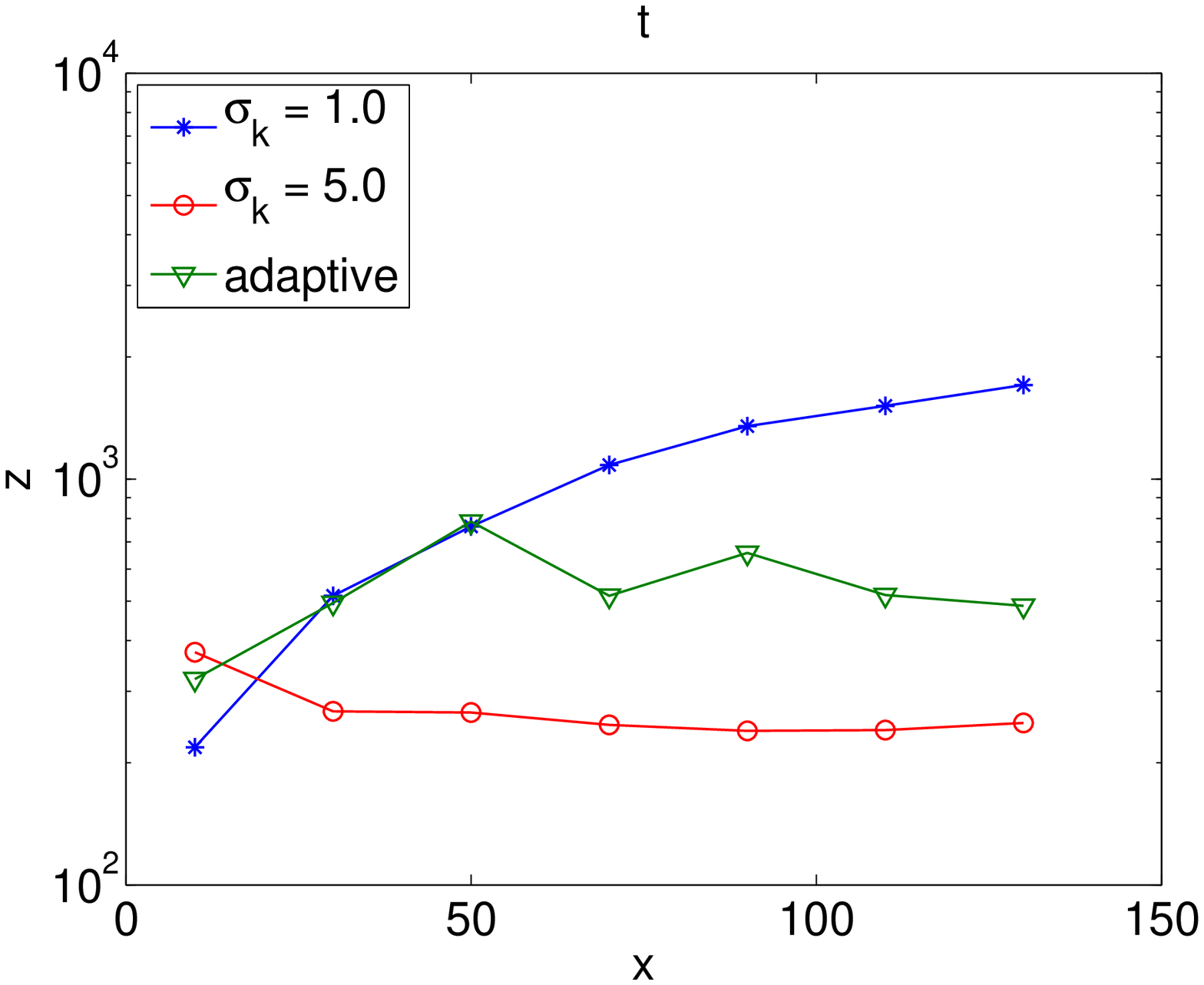}
\caption{\emph{Left.} Number of iterations for of Spingarn's
method based on the desired accuracy, for a problem instance with 
$d=20$, $w=20$, $s=10$, and using a fixed steplength parameter
$\sigma_k = 5$.
\emph{Right.} Number of iterations for the same problem with
$\epsilon = 10^{-4}$ and different choices of steplength.
}\label{f-accuracy}
\end{figure}

\paragraph{Accuracy and steplength selection}
The principal disadvantage of the splitting method, compared with
an interior-point method, is the more limited accuracy and the
higher sensitivity to the choice of algorithm parameters.
Figure~\ref{f-accuracy} (left) shows the number of iterations
versus $l$ for
different values of the tolerance $\epsilon$ used in the stopping
criterion.  
The right-hand plot shows the number of iterations versus $l$ for two
different constant values of the steplength parameter $\sigma_k$
($\sigma_k = 1.0$ and $\sigma_k=5.0$) and for an adaptively
adjusted steplength.

\section{Conclusions} \label{s-conclusion}
We have described a decomposition method that exploits partially separable
structure in linear conic optimization problems.  
The basic idea is straightforward: by replicating some of the variables, 
we reformulate the problem as an equivalent linear optimization 
problem with block-separable conic inequalities 
and an equality constraint that ensures that the replicated variables 
are consistent.
We can then apply Spingarn's method of partial inverses to this 
equality-constrained convex problem.
Spin\-garn's method is a generalized alternating projection method
for convex optimization over a subspace.
It alternates orthogonal projections on the subspace 
with the evaluation of the proximal operator of the cost function.
In the method described in the paper, these prox-operators are 
evaluated by an interior-point method for conic quadratic optimization.

When applied to sparse semidefinite programs, the reformulation
coincides with the \emph{clique conversion} methods which were
introduced in \cite{KKMY:11,FKMN:00} with the purpose of exploiting 
sparsity in interior-point methods for semidefinite programming. 
By solving the converted problems via a splitting algorithm 
instead of an interior-point algorithm 
we extend the applicability of the conversion methods to problems for 
which the converted problem is too large to handle by interior-point 
methods.
As a second advantage, if the correlative sparsity is block-diagonal,
the most expensive step of the decomposition algorithm 
(the evaluation of the proximal operator) becomes separable and can be  
parallelized.
The numerical experiments indicate that the approach is effective 
when a moderate accuracy (compared with interior-point methods)
is acceptable.   However the convergence can be quite slow and
strongly depends on the choice of steplength. 

A critical component in the decomposition algorithm for semidefinite
programming is the use of a customized interior-point method for
evaluating the proximal operators.  This technique allows us to
evaluate the proximal operator at roughly the same cost of solving
the reformulated SDP without the consistency constraints.
As a further improvement we hope to extend this technique to exploit
sparsity in the coefficient matrices of the reformulated problem, using
techniques developed for interior-point methods for sparse 
matrix cones~\cite{ADV:12}.

While sparse semidefinite programming provides the most important
application of our results, the techniques easily extend to other types
of partially separable cones.  In many of these extensions
partial separability does not require chordal structure 
(as it does in semidefinite programming).
As an example, second-order cone programs with partially separable
structure arise in machine learning problems involving 
sum-of-norm penalties that promote group sparsity \cite{BJMO:11}.

\newcommand{\etalchar}[1]{$^{#1}$}

\end{document}